\documentclass[11pt]{article}
\usepackage{amsmath,amsfonts,amssymb,amsthm,mathrsfs}
\usepackage[top=2.5cm,left=2.5cm,right=2.5cm,bottom=2.5cm]{geometry}
\usepackage[font=sf, labelfont={sf,bf}, margin=1cm]{caption}
\usepackage{graphicx,graphics}
\usepackage{epsfig}%pour les eps
\usepackage{latexsym}%encore des symboles
\usepackage[applemac]{inputenc}
\usepackage{ae,aecompl}
\usepackage[english]{babel}
\usepackage[colorlinks=true]{hyperref}
\usepackage{pstricks}
\usepackage{enumerate}

\newtheorem{theorem}{Theorem}[]
\newtheorem{theorem2}{Theorem*}[]

\newtheorem{proposition2}[theorem2]{Proposition*}
\newtheorem{corollary2}[theorem2]{Corollary*}
\newtheorem{remark}{Remark}[]
\newtheorem{definition}{Definition}[]
\newtheorem{proposition}[theorem]{Proposition}
\newtheorem{lemma}[theorem]{Lemma}
\newtheorem{corollary}[theorem]{Corollary}
\newtheorem{question}[]{Question}

\newcommand{\limLD}[3]{\underset{#3}{\mathsf{limLD}} \,\left(#1,#2\right)}

\title{\textsc{A glimpse of the conformal structure of random planar maps}}

\date{}

\author{Nicolas Curien\thanks{CNRS \& Université Paris 6, E-mail: nicolas.curien@gmail.com}}

\begin{document}
\maketitle
\begin{abstract}We present a way to study the conformal structure of random planar maps. The main idea is to explore the map along an SLE (Schramm--Loewner evolution) process of parameter $ \kappa = 6$ and to combine the locality property of the SLE$_{6}$ together with the spatial Markov property of the underlying lattice in order to get a non-trivial geometric information. We follow this path in the case of the conformal structure of random  triangulations with a boundary. 

Under a reasonable assumption called $(*)$ that we have unfortunately not been able to verify, we prove that the limit of uniformized random planar triangulations has a fractal boundary measure of Hausdorff dimension $\frac{1}{3}$ almost surely. This agrees with the physics KPZ predictions and represents a first step towards a rigorous understanding of the links between random planar maps and the Gaussian free field (GFF). \end{abstract}

%\noindent \emph{\textbf{AMS subject classifications:} }

%\medskip

\begin{figure}[!ht]
 \begin{center}
\begin{minipage}{6cm}\includegraphics[width=7cm]{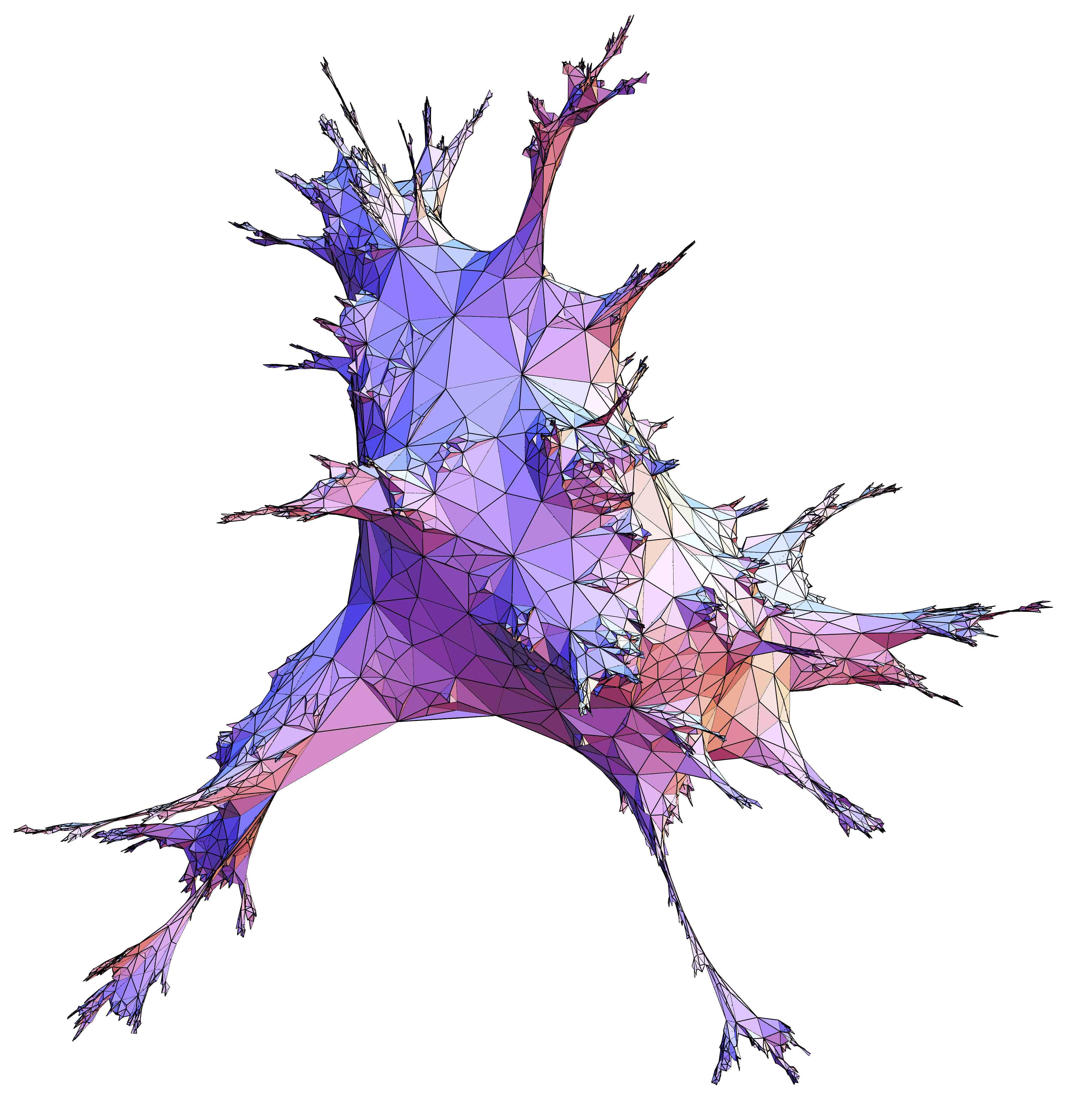}\end{minipage} $\quad$ \begin{minipage}{3cm} $\xrightarrow[]{ \mbox{``uniformization''}}$ \end{minipage} $\quad$  \begin{minipage}{6cm}\includegraphics[width=5.2cm]{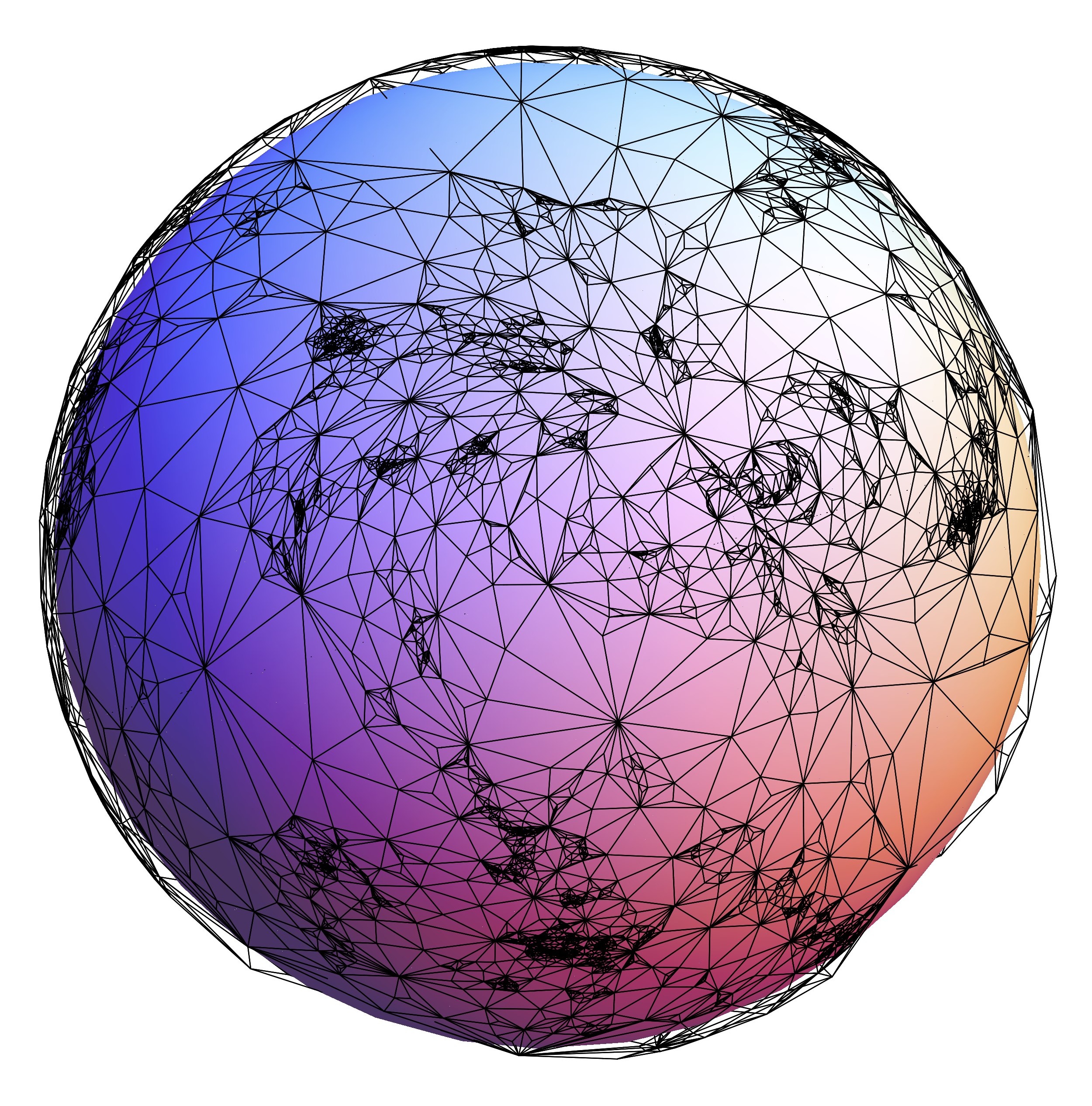} \end{minipage}
 \caption{A random triangulation embedded (not isometrically) in $ \mathbb{R}^3$ and an approximation of its uniformization on the two-dimensional sphere.}
 \end{center}
 \end{figure}

%%%%%%%%%%%%%%%%SECTION
\section*{Introduction} 
What does a typical random metric on the two-dimensional sphere look like? This concept plays a crucial role in the theory of two-dimensional quantum gravity where the famous KPZ relations (Knizhnik, Polyakov and Zamolodchikov  \cite{KPZ88}) are supposed to relate the dimensions of (some) sets under the random --or ``quantum''-- metric on the sphere $ \mathbb{S}_{2}$ to their dimensions with respect to the standard  Euclidean metric, see \cite{Gar12} for a smooth introduction. Nowadays, there are two mathematically rigorous approaches trying to make sense of  ``the random metric on $	\mathbb{S}_{2}$''.  

\paragraph{Random planar triangulations.} The first one is the theory of random planar triangulations (RPT) known as  ``dynamical triangulations'' in theoretical physics \cite{ADJ97}. The basic idea is to discretize a continuous surface into  finitely many triangles (or in any other basic tile) glued together: a triangulation that approximates the space. It seems natural to expect that such a discretization of  ``the random metric on $ \mathbb{S}_{2}$'' into $n$ triangles should yield a random triangulation  $T_{n}$ uniformly distributed over the set of all triangulations of $ \mathbb{S}_{2}$ with $n$ faces. 
 
Starting from this discrete model, Le Gall \cite{LG11} (see also Miermont \cite{Mie11} for the quadrangular case) has shown that after renormalizing the distances in $T_{n}$ by $n^{-1/4}$, the resulting random compact metric space indeed converges in distribution (for the Gromov--Hausdorff topology) towards a random compact metric space  called the Brownian map. This random metric space thus captures the metric properties of what a random metric on $ \mathbb{S}_{2}$ should be (in particular it is of Hausdorff dimension $4$ \cite{LG07}). However, although the Brownian map is known to be homeomorphic to the sphere  (see \cite{LGP08,Mie08}) the embedding is not canonically defined. The Brownian map cannot yet  be seen as $ \mathbb{S}_{2}$ endowed with a canonical random metric.

\paragraph{Gaussian free field.} The second approach is based on the Gaussian free field (GFF) which is a conformally invariant random distribution $h$ on the sphere. The ``random metric on $ \mathbb{S}_{2}$'' is then formally given by
 \begin{eqnarray} \label{eq:formalmetric} e^{\varrho h(z)} dz^2,\end{eqnarray} where $dz^2$ is the infinitesimal metric element on $ \mathbb{S}_{2}$ and $\varrho \geq 0$ is a parameter.  The last display would be easy to define if $h$ were a random smooth function, but unfortunately up to now, no rigorous construction is known to make sense of  \eqref{eq:formalmetric} (except in dimension one  \cite{BS11}), see \cite{MS13} for recent progress. Still, there are several equivalent ways to make sense of \eqref{eq:formalmetric} in terms of a random \emph{measure} and certain forms of the KPZ relations have been proved in this setup, see \cite{DS11,Kah85,RV11,RV13}.

\subsection*{Conformal structure of RPT} Though both paths have not succeeded in formally constructing a random metric living on $ \mathbb{S}_{2}$, we see that these approaches have different drawbacks: The RPT theory does yield a continuous metric but the embedding on the sphere is lacking, whereas in the GFF approach, the sphere (hence the embedding) is a built-in feature of the model but the random  metric seems hard to construct. However, the two theories are believed to eventually converge. This conjectured link has been made particularly clear (but remains unproven) by Duplantier \& Sheffield in \cite{DS11} and consists in understanding the \emph{conformal structure} of random planar maps (triangulations in this work) and to relate it to the GFF. For a nice exposition, see Garban's survey \cite{Gar12}. The goal of this work is to propose a possible way to rigorously begin this understanding.\medskip

\begin{figure}[!ht]
 \begin{center}
 \includegraphics[width=15cm]{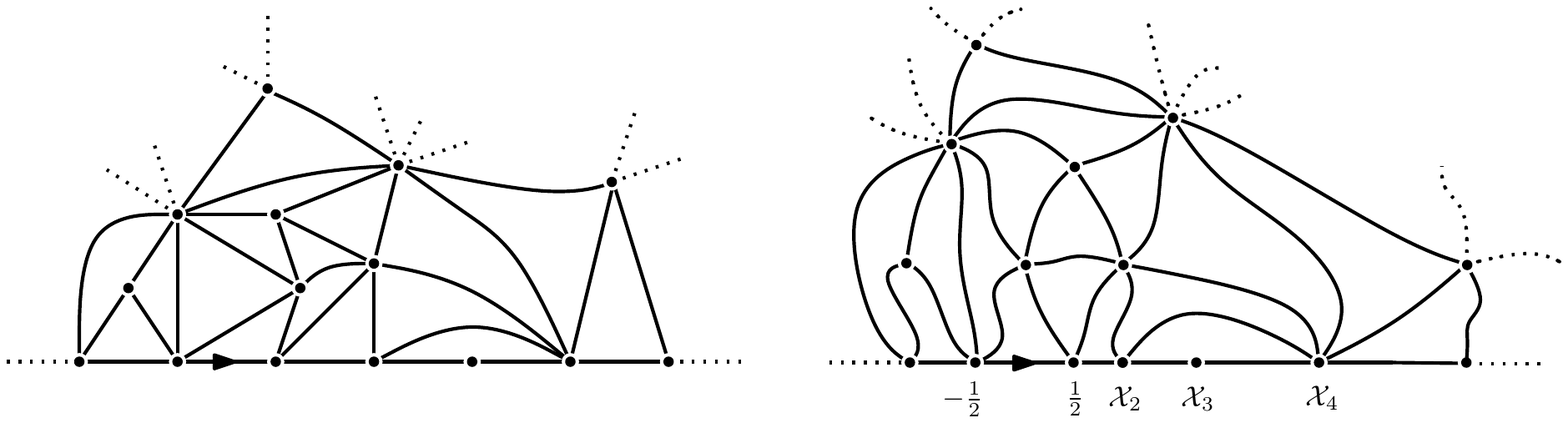}
 \caption{ \label{fig:uniformization} A part of a UIHPT --left-- and its uniformization --right-- (artistic drawing).}
 \end{center}
 \end{figure}

Formally, we focus here on the model of the uniform infinite half-planar triangulation (UIHPT) which is an infinite random triangulation $T_{\infty,\infty}$ with an infinite simple boundary obtained by Angel \cite{Ang05} as the local limit of triangulations with simple boundary whose size and perimeter both tend to infinity, see Section \ref{sec:UIHPT} for its definition and basics about planar maps. 
The UIHPT is also given with a distinguished oriented edge, called the root edge and oriented so that the infinite face is lying on its right, see Fig.\,\,\ref{fig:uniformization}.   From many respects, this model of random planar map is the simplest of all. The key property of this random lattice is its  particularly simple spatial Markov property which roughly says that after exploring a finite simply connected region of the map, then the remaining part is independent of the explored region and has the same law as the original lattice. See Section \ref{sec:peeling} for a precise statement. The spatial Markov property of random planar maps has been studied in details in \cite{AR13} and was at the core of many non-trivial results, see \cite{Ang03,ACpercopeel,BCsubdiffusive,MN13}. \medskip 

Our goal is to study the conformal structure of the boundary of the UIHPT. Formally, one consider $T_{\infty,\infty}$ as a random Riemann surface by seeing each triangle as an (Euclidean) equilateral triangle and gluing the charts along the edges and vertices of the map, see \cite{GR10} and Section \ref{sec:riemannconstruction} for details. Using the uniformization theorem, one can map the simply connected Riemann surface with a boundary obtained by the previous device onto the upper half plane $ \mathbb{H} = \{z \in \mathbb{C} : \Re(z) >0\}$. This mapping is unique provided that we fix the images of the origin and target of the root edge to be $-1/2$ and $1/2$ and send $\infty$ to $\infty$.

The conformal drawing  of the UIHPT (that is the image of the edges of $T_{\infty,\infty}$ by the above mapping) will be denoted by $ \mathscr{T}_{\infty,\infty}$ and we will commit an abuse of terminology when we will still speak about its vertices, edges and faces which are defined in an obvious way.  
For $k \geq 0$, the position of the $k$th vertex on the right of the origin of $ \mathscr{T}_{\infty,\infty}$ is denoted by $ \mathcal{X}_{k}$, in particular $ \mathcal{X}_{0}= -1/2$ and $ \mathcal{X}_{1}=1/2$.  For $n \geq 1$, we consider the random  probability measure $\mu_{n}$  on $[0,1]$ defined by
$$ \mu_{n}=  \frac{1}{n} \sum_{k=1}^n \delta_{\mathcal{X}_{k} / \mathcal{X}_{n}}. $$

\begin{theorem2} \label{thm:main}From any sequence of integers tending to $+\infty$ one can extract a subsequence $n_{k} \to \infty$ such that $\mu_{n_{k}}$ converges as $k \to \infty$ in distribution towards a random probability measure $\mu$  such that almost surely 
\begin{itemize}
\item $\mu$ is non-atomic,
\item $\mu$ has topological support equal to $[0,1]$,
\item the Hausdorff dimension\footnote{Recall that the dimension of a measure is the infimum of the dimensions of Borel sets of full mass} of $\mu$ is  $1/3$.
\end{itemize}
\end{theorem2}
\paragraph{The star condition.} We used the label Theorem* because our proof relies on an assumption  denoted by $(*)$ (see Section \ref{sec:introstar} for its definition)  that we strongly believe to hold, but have not been able to rigorously derive. Similarly, the results denoted by  Proposition*, Corollary*, Lemma* etc...  all rely on $(*)$.  Interesting on its own, the assumption $(*)$ is thus strongly motivated by the conditional results proved in this paper. See Section \ref{sec:comments} for a discussion and supports for  $(*)$. \medskip

The random measures $\mu_{n}$ are believed to converge (without the need to pass to a subsequence), and the candidate for the limiting random measure $\mu$ is defined as follows, see \cite{DS11,Gar12}. 
Let $h = \tilde{h}+h_{0}$ where $\tilde{h}$ is an instance of the mean zero Gaussian free field (GFF) on $ \mathbb{H}$ with zero boundary condition (see \cite[Section 3]{She10}) and $h_{0}(z) = - \varrho \log|z|$. We can define a random measure $\mu_\varrho$ on $ \mathbb{R}$ formally obtained as $$ e^{\varrho h(z)/2}dz,$$ where $dz$ is the Lebesgue measure on $ \mathbb{R}$. This random measure can be constructed using Kahane's theory of Gaussian multiplicative chaos or by means of regularization procedures, see \cite{RV10}, \cite[Section 6]{DS11} and \cite[Section 5]{RV13}. 
%  Since we (formally) have $ E[h(x)h(y)] = - 2 \log |x-y|$, the measure $\mu_{\varrho}$ can be constructed using Kahane's theory of Gaussian multiplicative chaos, see \cite{RV10} where $\mu_{\varrho}$ is called the log-normal multi-fractal random measure. An alternative equivalent construction of this measure using circle average regularization can be found in  \cite[Section 6]{DS11}. We refer the reader to the review \cite[Section 5]{RV13} for details of these constructions. For the particular parameter $\varrho = \sqrt{8/3}$ this measure is conjectured to describe the limit of the discrete counting measure on the boundary of a conformally uniformized map similar as the one considered in this work. More precisely, it is conjectured (see \cite{DS11,She10}) that:
\begin{question}[see {\cite[Conjecture 7.1]{DS11},\cite{She10}}] \label{ques:limit} Do the random measures $\mu_{n}$ converge in law towards the random  measure $$ \tilde{\mu}_{\varrho}=\frac{\mu_{\varrho}( \cdot \cap [0,1])}{\mu_{\varrho}[0,1]} \qquad \mbox{ with } \quad \varrho = \sqrt{8/3} \ ?$$\end{question}

Duplantier \& Sheffield \cite{DS11} and Rhodes \& Vargas \cite{RV11} recently showed that the KPZ relations derive from the analysis of the multi-fractal spectrum of the random measure $\mu_{\varrho}$. This analysis has  been undertaken first in \cite{BMKPZ} where it is shown that the dimension of $\mu_{\varrho}$ is $1/3$ when $\varrho = \sqrt{8/3}$. % (see also the proof of Theorem 6.1 in \cite{DS11} for a different approach). 
Hence Theorem* \ref{thm:main} strongly supports  Question \ref{ques:limit}.
% and can be seen as the first step towards the understanding of the full multi-fractal spectrum of $\mu$ and thus of the KPZ relations. 

%To our opinion, the value of this paper does not entirely lie in the statement of Theorem* \ref{thm:main} but rather on the strategy (that we believe is new) employed to study the conformal geometry of random planar maps. In particular, i
It might be the case that our approach actually yields a characterization of the subsequential-limits of the $\mu_{n}$'s which is shared with $\tilde{\mu}_{\rho}$ for $\rho = \sqrt{8/3}$, see Question \ref{ques:characterization} below. A positive answer to Question \ref{ques:characterization} and assumption $(*)$ would turn Question \ref{ques:limit} into a theorem.
%If true, then a positive answer to Question 1 would follow from $(*)$.  which could then be identified with the measure constructed from the exponential of the Gaussian free field in \cite[Section 6]{DS11} for the particular parameter $\varrho = \sqrt{8/3}$, see the last paragraph of this Introduction. We refer to Section \ref{sec:comments} for a discussion and for the statements of several conjectures about links between the GFF,  SLE and stable processes which we hope can directly be answered by specialists in the field.

%%%%%%%%%%%%%%%%%%%%SUBSECTION
\subsection*{Strategy} Our approach to investigate the conformal structure of random planar triangulations  is based on their exploration by an independent SLE$_{6}$ process. Recall that for $\kappa \geq 0$, the SLE$_{ \kappa}$   processes have been introduced by Schramm \cite{Sch98} in order to describe interfaces of conformally invariant models in two dimensions. See \cite{Law05,Wer04} for background. The SLE$_{6}$ process has a characteristic feature (that he shares with Brownian motion), which is called the locality property. The latter roughly means that its growth is locally defined and does not depend on the full curve, see Section \ref{sec:locality}. This property is one of the keys in  the determination of the Brownian intersection exponents by Lawler, Schramm and Werner \cite{LSW01} and is also central in this work.

Formally, the exploration of the UIHPT by an SLE$_{6}$ is defined as the exploration of $ \mathscr{T}_{\infty,\infty}$ by an independent chordal SLE$_{6}$ on $ \mathbb{H}$ started from $0$. \emph{A priori}, the exploration of the UIHPT thus depends on its  whole conformal structure since we formally need its uniformization $ \mathscr{T}_{\infty,\infty}$ to define it. However, the locality property of SLE$_{6}$ will imply that this exploration can in fact be performed by discovering the UIHPT ``step-by-step'' revealing only the parts necessary for the SLE$_{6}$ to displace. 

\begin{figure}[!ht]
 \begin{center}
 \includegraphics[width=14cm]{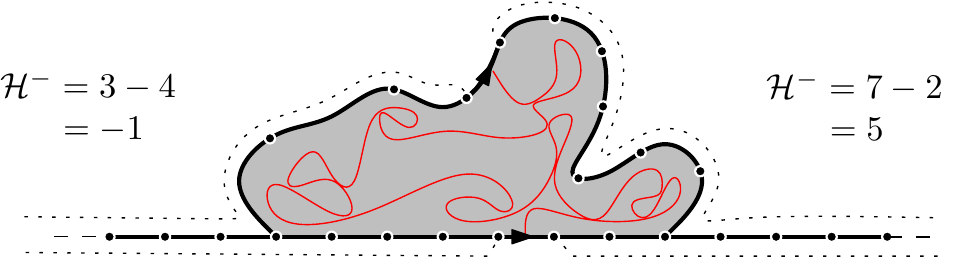}
 \caption{\label{fig:horo}Definition of the horodistances. After a piece of the UIHPT has been explored by the SLE process, one can define the variations of the boundary distances in the remaining triangulation from the edge on which the SLE is located. }
 \end{center}
 \end{figure}

This will show that the SLE$_{6}$ exploration of the UIHPT is \emph{Markovian}, in the sense that the submap discovered after some time (made of the triangles traversed by the SLE as well as the finite regions they enclose) is independent of the remaining of the map which is distributed as a standard UIHPT, see Section \ref{sec:markov}. 
Using Angel's peeling process (see \cite{Ang03,ACpercopeel} and Section \ref{sec:one-step}), we are able to understand the algebraic lengths of the boundary seen from $\pm \infty$ in the unexplored map. More precisely when the SLE is located on a boundary edge of the explored region, we can define two integer numbers $ \mathcal{H}^+$ and $ \mathcal{H}^-$ representing the variations of the boundary lengths towards $\pm\infty$ from this edge compared to the original boundary lengths from the root edge of the map, see Fig.\,\,\ref{fig:horo}.

In Theorem*\,\ref{thm:exploSLE} we show, under assumption $(*)$, that this \emph{horodistance} process  $(\mathcal{H}^{+}(i), \mathcal{H}^-(i))_{i\geq0}$ is mainly driven by the spatial Markov property of the map and converges (in distribution in the Skorokhod sense) after normalization by $n^{2/3}$ towards a pair $(S^{+},S^{-})$ made of independent standard $\frac{3}{2}$-stable spectrally negative Lévy processes, more precisely
\begin{eqnarray} \left(	\frac{ \mathcal{H}^{+}([nt])}{n^{2/3}}, \frac{  \mathcal{H}^{-}([nt])}{n^{2/3}}\right) & \xrightarrow[n\to\infty]{(d)} & 3^{-2/3} \cdot \big(S^{+}_{t},S^{-}_{t}\big) _{t \geq 0}.  \label{eq:introstable}\end{eqnarray}
The basic idea of the proof of Theorem* \ref{thm:main} is  to connect these horodistance processes to a geometric property, namely the fact that the SLE$_{6}$ bounces off $ \mathbb{R}_{+}$ and $ \mathbb{R}_{-}$,  \cite{BOW}. On an intuitive level at least, the times when $ \mathcal{H}^{-}$ (resp.\,\,$ \mathcal{H}^{+}$) reaches a new minimum correspond to the visits of $ \mathbb{R}_{-}$ (resp.\,\,$\mathbb{R}_{+}$) by the SLE$_{6}$ process, see Section \ref{sec:bouncing}.  We then compute, in two ways, the number $C( \varepsilon,n)$ of times the SLE$_{6}$ exploration of $ \mathscr{T}_{\infty,\infty}$ is alternatively bouncing of $ \mathbb{R}_{+}$ and $ \mathbb{R}_{-}$  between the point $ \mathcal{X}_{[ \varepsilon n]}$ and $ \mathcal{X}_{n}$. On the one hand, using the scaling limit of the horodistance process \eqref{eq:introstable} one is capable to compute $C( \varepsilon,n)$ (to be precise, its limit) in terms of interlaced minimal records of $S^+$ and $S^-$ (see Corollary \ref{cor:commuteStable})  and we find that as $ n\to \infty$
  \begin{eqnarray} \label{eq:intro1} C( \varepsilon,n) \quad \approx \quad   \frac{3 \sqrt{3}}{2\pi} |\log \varepsilon|. \end{eqnarray}
On the other hand, conditionally on  $ \mathscr{T}_{\infty,\infty}$ (and a fortiori on $( \mathcal{X}_{k})_{k \geq 0}$) it is known  (see   Corollary \ref{cor:commuteSLE} below or the related computation of Hongler and Smirnov \cite{HS11}) that the number of alternative commutings to $ \mathbb{R}_{+}$ and $ \mathbb{R}_{-}$ an SLE$_{6}$ is doing after having swallowed the point $  \mathcal{X}_{[\varepsilon n ]}$ until it swallows the point $ \mathcal{X}_{n}$ is roughly of order 
  \begin{eqnarray} \label{eq:intro2} C( \varepsilon,n) \quad \approx \quad \frac{ \sqrt{3}}{ 2 \pi} \log  \frac{ \mathcal{X}_{n}}{ \mathcal{X}_{[ \varepsilon n]}}.  \end{eqnarray} Equalizing \eqref{eq:intro1} and \eqref{eq:intro2} we find that $ \mathcal{ X}_{[ \varepsilon n]}/ \mathcal{X}_{n} \approx \varepsilon^{3}$ or  in terms of the limiting random measure $\mu$  that $\mu[0, \varepsilon^3] \approx \varepsilon$ or equivalently $\mu[0, \varepsilon] \approx \varepsilon^{1/3}$. This is the main idea of the proof of Theorem* \ref{thm:main} (iii).\medskip 
  
  % We conclude this introduction by pointing the conjectured link of our work with the Gaussian free field and the KPZ relations (at least in the case of ``pure gravity''). The interested reader should consult \cite{DS11,She10} for details,  references and additional open questions.

% More precisely, the formula is equivalent (roughly speaking) to the fact that if conditionally on $\tilde{\mu}_{ \varrho}$, we sample a point $X$ according to $ \tilde{\mu}_{\varrho}$ then we have (see the proof of Theorem 6.1 in \cite{DS11})  $$ \lim_{r \to 0^+} \frac{\log P( \mu_{\varrho} B_{r}(X)  \approx r^{\alpha})}{\log r} =  \frac{1}{2  \alpha \varrho}\left(\alpha - \frac{(2-\varrho^2/2)}{2}\right)^2.$$  hen $\varrho = \sqrt{8/3}$, the multi-fractal spectrum  vanishes for $\alpha = 1/3$ which ``implies'' that $\mu_{\varrho}$ almost surely has dimension $ \frac{1}{3}$: This  is precisely the content of Theorem* \ref{thm:main}. We can thus see our work as the first step towards the understanding of the full multi-fractal spectrum of $\mu_{\varrho}$ and thus of the KPZ relations. 

\medskip  The paper is organized as follows. In the first section we recall  the background on the UIHPT including its construction and the crucial spatial Markov property. The notion of Markovian exploration is introduced as well as  basics on the $ \frac{3}{2}$-stable process. The second section is devoted to the SLE$_{6}$ exploration of the UIHPT. We explain there why the locality property of the SLE resonates with the spatial Markov property of the underlying lattice and implies under assumption $(*)$ the convergence \eqref{eq:introstable}. In Section \ref{sec:bouncing}, we show how to translate \eqref{eq:introstable} into geometric information by studying the alternative bouncings of the SLE and interlaced minimal records of two independent stable processes. The proof of Theorem* \ref{thm:main} can be found in Section \ref{sec:proof}. The last section  contains conjectures, comments and possible extensions for future works.

\bigskip 
\noindent \textbf{Acknowledgments:} We are indebted to Jean Bertoin and Jean-Fran\c cois Le Gall for crucial advices on how to prove Proposition \ref{prop:commute} and \ref{prop:commuteSLE} respectively. Thanks also go to Wendelin Werner for a stimulating discussion and to Vincent Vargas for helpful comments on a first version of this work. We are grateful to the organizers of the  conference ``Planar statistical physics (2012)'' in les Diablerets where  this work started during a transit on the everlasting ski-lift ``Perche-Conche''.  Finally, we thank the anonymous referee for helpful comments.
\tableofcontents
 
%%%%%%%%%%%%%%%%SECTION
\section{Background on the half-plane UIPT}
 \label{sec:UIHPT}

%%%%%%%%%%%%%%%%%%%%SUBSECTION
\subsection{The UIHPT}
\label{sec:basics}
\paragraph{Triangulations.}
Recall that a planar map is a finite connected planar graph embedded in the sphere $ \mathbb{S}_{2}$ seen up to continuous deformations that preserve the orientation. There is a natural notion of vertex, edge and face in a planar map. The degree of a face is the number of half-edges surrounding the face.  As usual, all the maps considered in this work are rooted, that is, given with a distinguished oriented edge $ \vec{e}$ called the root of the map. A triangulation is a map whose faces are all of degree $3$. We will focus on type-II triangulations, that are triangulations without loops but possibly multiple edges.

A triangulation with a simple boundary is a planar map whose faces are all triangles except possibly the face on the right-hand side of the root edge called the external face which is bounded by a non-intersecting cycle (no pinch-points). In this work we only deal with simple boundaries and thus sometimes drop the adjective simple to lighten the writing. The perimeter of a triangulation with boundary is the degree of the external face, and a triangulation with a boundary of perimeter $p$ is also called a triangulation of the $p$-gon. The size of a triangulation with a boundary is its number of inner vertices (i.e.~not located on the boundary). By convention, the only  triangulation of the $2$-gon of size $0$ is made of a single oriented edge.

\begin{figure}[!ht]
 \begin{center}
 \includegraphics[width=7cm]{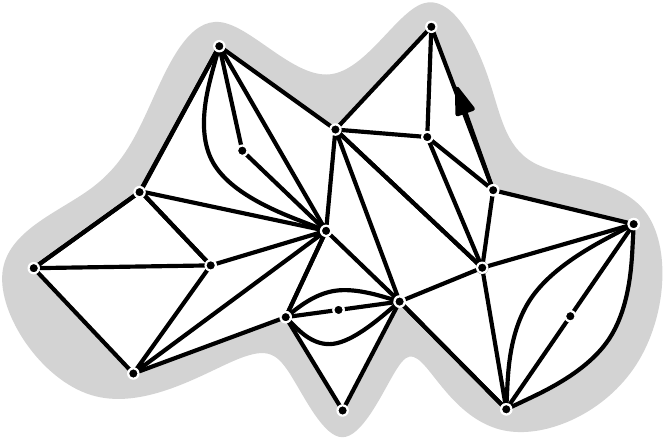}
 \caption{A triangulation of the $12$-gon of size $7$.}
 \end{center}
 \end{figure}
\paragraph{Local limits.}
Following \cite{AS03,BS01} we recall the local topology on the set of planar maps. If $m,m'$ are two rooted maps, the local distance between $m$ and $m'$ is 
 \begin{eqnarray*} \mathrm{d_{loc}}(m,m') &=& \big(1+ \sup\{r\geq 0 : B_{r}(m)=B_{r}(m')\}\big)^{-1}, \end{eqnarray*}
where $B_{r}(m)$ denotes the map formed by the vertices and edges of $m$ which are at graph distance smaller than or equal to $r$ from the origin of the root edge in $m$. The set of all finite rooted triangulations with boundary in not complete for this metric and we shall work in its completion obtained by adding infinite maps (see \cite{CMMinfini} for a detailed exposition in the quadrangular case). For any $p \geq 2$, we denote by $ T_{n,p}$ a random variable uniformly distributed on the set $ \mathcal{T}_{n,p}$ of all triangulations (of type II) of the $p$-gon having size $n$. The \textbf{U}niform \textbf{I}nfinite \textbf{H}alf-\textbf{P}lanar \textbf{T}riangulation (UIHPT) is obtained as a local limit of uniform triangulations with boundary by first letting their sizes tend to infinity and next sending their perimeters to infinity. More precisely, Angel \& Schramm \cite{AS03} and Angel \cite{Ang05} proved the following convergences in distribution for $ \mathrm{d_{loc}}$\begin{eqnarray*}
T_{n,p} \quad  \xrightarrow[n\to\infty]{(d)} \quad  T_{\infty,p} \quad \xrightarrow[p\to\infty]{(d)} \quad  {T}_{\infty,\infty}, \end{eqnarray*}
where $T_{\infty,p}$ is a random rooted infinite triangulation of the $p$-gon called the UIPT (for Uniform Infinite Planar Triangulation) of the $p$-gon and $T_{\infty,\infty}$ is the UIPT of the half-plane denoted by UIHPT (see \cite{CMboundary} for similar statements in the quadrangular case). This is the main character of this paper. 

The root edge of $T_{\infty,\infty}$ will always be denoted by $ \vec{e}$ (the external face is on its right) or $e$ if we consider the unoriented edge. The infinite simple boundary of $T_{\infty,\infty}$ can be identified with $ \mathbb{Z}$ by declaring that the root edge is $0 \to 1$. The UIHPT enjoys an invariance under re-rooting:  for any $k \in \mathbb{Z}$ the planar map obtained from $ T_{\infty,\infty}$ by re-rooting at the edge $k \to k+1$ is still distributed as the UIHPT. For this reason we might be loose on the precise location of the root edge in what follows.  

%%%%%%%%%%%%%%%%%%%%SUBSECTION
\subsection{One-step peeling of the UIHPT}
\label{sec:markov}
\label{sec:peeling}

One of the very nice features of the UIHPT is its spatial Markov property that can roughly be described as follows: Assume that we explore a simply connected region $ \mathcal{R}$ of $T_{\infty,\infty}$ that contains the root edge, then the exterior of $ \mathcal{R}$ is independent of $ \mathcal{R}$ and is distributed  as UIHPT. This describes the conditional laws of the different maps we obtain from $T_{\infty,\infty}$ after conditioning on the face that contains the root edge $e$. See \cite{Ang05,ACpercopeel} for details and proofs. \medskip 
%In order to make this statement precise we shall need some background. See \cite{ACpercopeel} for more details.

%\paragraph{Enumeration.}
First we recall the standard asymptotic  %We gather here several results about (asymptotic) enumeration of planar triangulations, see \cite{Ang05,Ang03} and the references therein for their derivations. For $n\geq 0$ and $p \geq2$ we have the following asymptotic equivalent for $ \# \mathcal{T}_{n,p}$:
$
   \# \mathcal{T}_{n,p}   \underset{n\to\infty}{\sim}  C_{p}\big({27}/{2}\big)^n n^{-5/2}$ for some $C_{p}>0$. So the series $\sum_{n\geq 0}  \#  \mathcal{T}_{n,p}(2/27)^n$ is finite and its sum is denoted by $Z_{p}$ (see \cite{Ang03} for exact expressions  of $Z_{p}$ and $C_{p}$). 
%  \begin{eqnarray} \label{zp}   Z_{p} &=& \frac{(2p-4)!}{(p-2)!p!}\left( \frac{9}{4}\right)^{p-1} \qquad \mbox{for } p \geq 2.\end{eqnarray}
  
 \begin{definition} \label{def:bolts} The free Boltzmann distribution of the $p$-gon is the probability measure on $ \cup_{n\geq 0}\mathcal{T}_{n,p}$ that assigns a weight $(2/27)^{n} Z_p^{-1}$ to each triangulation of the $p$-gon of size $n$.\end{definition}
 
\label{sec:one-step}

 Let $ T_{\infty,\infty}$ be a UIHPT. Assume that we reveal the face on the left of the root edge $ \vec{e}$, this operation is called the one-step peeling transition. Three (or two by symmetry) situations may appear depending on the ``form'' of the triangle revealed. Let us make a list of the possibilities and describe the probabilities and the conditional laws for each case. The set of forms is   \begin{eqnarray*} \mathsf{Forms} &:=& ( \mathrm{C}, 1) \cup \bigcup_{k \geq 1}\{( \mathrm{G},-k),(  \mathrm{D},-k)\}. \end{eqnarray*} 
 To help the reader remind the notation remember that ``$ \mathrm{C}$'' stands for center, ``$ \mathrm{G}$'' for gauche (left in French) and ``$  \mathrm{D}$'' for droite (right in French) and that the numbers $1$ or $-k$  represent the variation of the number of edges on the boundary. Here are all the possible cases: 
 
 \begin{itemize}
\item  The revealed triangle could simply be a triangle with a vertex lying in the interior of $ T_{\infty,\infty}$ (i.e. not on the boundary), see Fig.\,\ref{fig:peel1}. We say that the revealed triangle is of form $( \mathrm{C},1)$. This event happens with probability $ q_{1}$ where 
   \begin{eqnarray*} q_{1} & = & \frac{2}{3}.  \end{eqnarray*}   The remaining triangulation (in gray in Fig.\,\ref{fig:peel1}) denoted by $ \mathrm{Peel}( T_{\infty,\infty}; {e})$ is formed after removing the revealed triangle from $T_{\infty,\infty}$ and rooting the resulting map at the edge of the revealed triangle which is incident to the initial root vertex. Conditionally on this event  $ \mathrm{ Peel}( T_{\infty,\infty};e)$ has the same distribution as $T_{\infty,\infty}$.
\begin{figure}[!ht]
\begin{center}   \includegraphics[width=7cm]{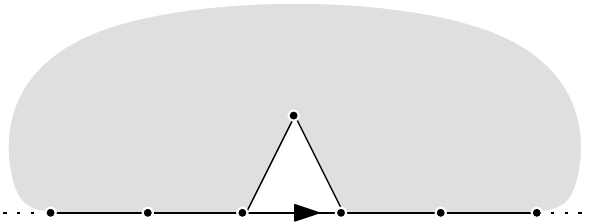} 
 \caption{ \label{fig:peel1}Case $ ( \mathrm{C},1)$.}
 \end{center}
 \end{figure}

 \item Otherwise, the revealed triangle has its three vertices lying on the boundary of $T_{\infty,\infty}$ and the third one is either $k \geq 1$ edges on the left of the root edge, in which case the triangle is said to be of form $( \mathrm{G},-k)$ or $k$ edges on the right of the root edge in which case the triangle is said to be of form $(  \mathrm{D},-k)$, see Fig.\,\ref{fig:peel2}. Note that $k >0$ because loops are not allowed since we are working with $2$-connected triangulations. By symmetry, these two events have the same probability $q_{-k}$ where
 \begin{eqnarray*} q_{-k} &=&\frac{ 
(2k-2)!}{ 
4^k (k-1)!(k+1)!}. \end{eqnarray*}
The revealed triangle thus encloses a triangulation with simple boundary of perimeter $k+1$ (the part in dark gray on Fig.\,\ref{fig:peel2}). Since $T_{\infty,\infty}$ has only one end, this enclosed part must be finite. The remaining infinite triangulation $ \mathrm{Peel}(T_{\infty,\infty}; {e})$ is formed by removing the revealed triangle and the enclosed triangulation from $T_{\infty,\infty}$ and rooting the resulting infinite triangulation with infinite boundary at the only edge adjacent to the revealed triangle. 

Then, conditionally on the fact that the revealed triangle has its third vertex lying $k$ edges away from the root edge, the enclosed triangulation and $ \mathrm{Peel}(T_{\infty,\infty}; {e})$ are independent, the first one follows a Boltzmann of the $k+1$-gon (see Definition \ref{def:bolts}) and $ \mathrm{Peel}(T_{\infty,\infty}; {e})$ is a UIHPT.

\begin{figure}[!ht]
 \begin{center}
\includegraphics[width=7cm]{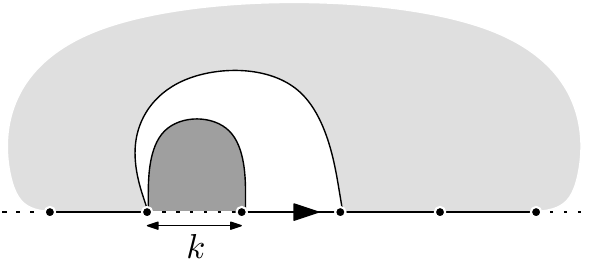}
 \caption{Case $( \mathrm{G},-k)$. \label{fig:peel2}}
 \end{center}
 \end{figure}
 \begin{remark}Conditionally on $k=1$, the enclosed triangulation can be of size $0$ with probability $Z_{2}^{-1} = 8/9$ in which case the revealed triangle is  glued on the boundary, see Fig.\,\ref{fig:swallowed} below.
\end{remark}
 
 \end{itemize}
 
After peeling the root edge, the triangle revealed may thus have two -if the form is $ ( \mathrm{C},1)$) or one (if the form is $( \mathrm{G},-k)$ or $(  \mathrm{D}, -k)$) edges which are part of the boundary of $ \mathrm{Peel}( T_{\infty,\infty}; e)$. These edges are called the \emph{exposed} edges as in \cite{ACpercopeel}. Also the edges of the boundary of $T_{\infty,\infty}$ except the peeled edge which are not part of the new boundary of $ \mathrm{Peel}(T_{\infty,\infty};e)$ are called the \emph{swallowed} edges. See Fig.\,\ref{fig:swallowed}. In the rest of the paper, we denote by $ \mathsf{F}$ a random variable over $ \mathsf{Forms}$ which has the law of the form of a one-step peeling of the UIHPT, that is  \begin{eqnarray} \label{def:F} P\big( \mathsf{F} = ( \mathrm{C},1)\big) = q_{1}, \quad \mbox{and} \quad P\big( \mathsf{F}=( \mathrm{G},-k)\big)=P\big( \mathsf{F}=( \mathrm{D},-k)\big) = q_{-k} \quad \mbox{ for }k \geq 1.  \end{eqnarray}  

\begin{figure}[!ht]
  \begin{center}
  \includegraphics[width=16cm]{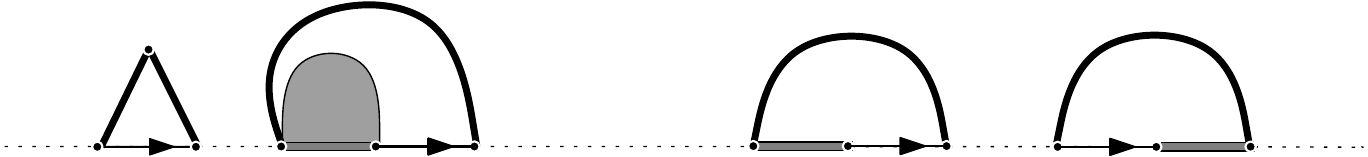}
  \caption{ \label{fig:swallowed}The exposed edges are in fat black lines and the swallowed
      ones are in fat gray lines. One the right, the two cases when the form is respectively $( \mathrm{G},-1)$ and $( \mathrm{D},-1)$ and the enclosed triangulation is of size $0$.}
  \end{center}
  \end{figure}
  
%%%%%%%%%%%%%%%%%%%%SUBSECTION
\subsection{Markovian exploration} \label{sec:horo}

 Let $T$ be an infinite triangulation with an infinite simple boundary such that $T$ has only one end. Extending what we have done in the case of the one-step peeling of $T_{\infty,\infty}$, for any non-oriented edge $a$ on the boundary of $T$ we denote by $ \mathrm{Peel}(T;a)$ the triangulation obtained from $T$ by removing the triangle adjacent to $a$ as well as the finite region it may enclose, rooted as in the preceding section. Similarly, define the form of the revealed triangle as before. We call this operation peeling the edge $a$ in $T$. \medskip 

An \emph{exploration} of $T_{\infty,\infty}$ is a sequence of nested subtriangulations\footnote{When we say a sequence of nested subtriangulations, we imagine that they are already given by nested embeddings. Indeed, in the case of presence of symmetries there could be many ways to see $T_{1}$ as a subtriangulation of $T_{0}$ etc... We do not intend to give a formal meaning to this and count on the intuition of the reader.} of $T_{\infty,\infty}$  \begin{eqnarray*}
\cdots \subset T_{n} \subset \cdots  \subset T_{2} \subset T_{1}\subset T_{0}=T_{\infty,\infty}  \end{eqnarray*} such that for any $i \geq 0$ the triangulation $T_{i+1}$ is obtained from $T_{i}$ by the peeling of an  edge %(necessarily unique)
 $a_{i}$ on the boundary of $T_{i}$. For each $i \geq 0$, we denote by $ K_{i}$ the ``complement'' triangulation of $T_{i}$ in $ T_{\infty,\infty}$ made of all the triangles peeled  at time $i$ as well as the finite regions they enclose. For definiteness, $K_{0}$ is the empty set. Alternatively, $K_{i}$ is obtained by cutting in $T_{\infty,\infty}$ along the boundary of $T_{i}$. This object is necessarily made of finitely many disjoint finite triangulations with simple boundary and will be called the ``known, explored or discovered'' part at time $i$ as opposed to $T_{i}$ which is the ``unknown, unexplored or undiscovered'' part. 

In this work, we further assume that $a_{0}$ is the root edge and that for $i \geq 0$ the edge $a_{i}$ to be peeled at time $i \geq 1$ is  located on the boundary of $K_{i}$ so that $(K_{i})_{i \geq 0}$ is a sequence of growing triangulations with simple boundary (there is a single growing component). \medskip  

Here comes the central notion introduced by Angel \cite{Ang03}:
\begin{definition}[Markovian exploration] \label{def:markovian} An exploration process is \emph{Markovian} if for every $i \geq 0$ the edge $a_{i}$ to peel at time $i$ is chosen using a (possibly random) algorithm that can use the knowledge of $K_{i}$ but does not depend on the unknown part $T_{i}$. 
\end{definition}

During a Markovian exploration of the UIHPT (also called a peeling process in \cite{ACpercopeel}) the peeling steps are iid. This has first been used by Angel in \cite{Ang03}, see also \cite[Proposition 4]{ACpercopeel}.

\begin{proposition} \label{prop:markoviid} During a \emph{Markovian} exploration of $T_{\infty,\infty}$ we have:
\begin{itemize}
\item For each $i \geq 0$, the half-planar triangulation $T_{i}$ is independent of $K_{i}$ and has the law of $T_{\infty,\infty}$,
\item The forms $(F_{i})_{i \geq 0}$ of the triangles revealed during the exploration are i.i.d.\,\,copies of $ \mathsf{F}$.
\end{itemize}
 \end{proposition}
%  \proof The proof is done by induction. Imagine that at step $i \geq 0$, the triangulation $T_{i}$ to be explored is a UIHPT independent of the explored part $K_{i}$. Since the edge $a_{i}$ to be peeled is chosen independently of $T_{i}$, by the re-rooting property of the UIHPT (see Section \ref{sec:basics}) we deduce that the triangulation $T_{i}$ re-rooted at $a_{i}$ still has the law of the UIPT of the half-plane. We then trigger a one-step peeling transition and reveal the face adjacent to $a_{i}$ as well as the finite region it encloses. By the last section, the form $F_{i}$ of this face is independent of $K_{i}$ and has the same distribution as $ \mathsf{F}$, besides $T_{i+1} = \mathrm{Peel}(T_{i};a_{i})$ is independent of $K_{i+1}$ and has the law of the UIHPT. \endproof

%%%%%%%%%%%%%%%%%%%%SUBVENTION
\subsection{Horodistances}
Let $ \mathcal{E}$ be an exploration process of the UIHPT. We will now keep track of the position of the peeling position with respect to $-\infty$ and $+\infty$ by using ``horodistances''. We use the notation introduced in the last sections where the dependance in $ \mathcal{E}$ is implicit. 

\begin{figure}[!ht]
 \begin{center}
 \includegraphics[width=15cm]{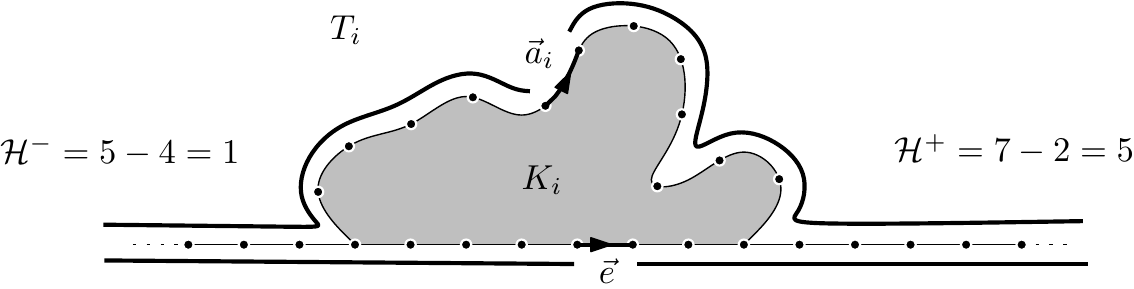}
 \caption{\label{fig:horodistances}Definition of the horodistances.}
 \end{center}
 \end{figure}

\paragraph{Definition.} Imagine that at step $i \geq 0$ we have discovered a subtriangulation  $K_{i} \subset T_{\infty,\infty}$ and that the next edge to peel is $ \vec{a}_{i}$ oriented such that the external face of $ T_{i}$ is on its right. We define two integer numbers $ \mathcal{H}^{-}(i)$ and $ \mathcal{H}^{+}(i)$ which represent the variations of the distances seen from $-\infty$ and $+\infty$ of the edge $ \vec{a}_{i}$ along the boundary. The definition should be clear on Fig.\,\,\ref{fig:horodistances}.

 Formally, denote $ \vec{a}_{i}^-$ the origin of $ \vec{a}_{i}$ and $ \vec{e}^-$ the origin of the root edge in $T_{\infty,\infty}$. Next consider the path $\gamma$ going from $ \vec{e}^{-}$ towards ``$-\infty$'' along the boundary of $T_{\infty,\infty}$ and $\gamma'$ the path going from $ \vec{a}^{-}_{i}$ towards ``$-\infty$'' along the boundary of $T_{i}$. Since $K_{i}$ is finite, these two paths eventually merge and $\gamma \backslash (\gamma\cap \gamma')$ as well as $\gamma' \backslash (\gamma\cap \gamma')$ are both finite. We  define $ \mathcal{H}^{-}(i)$ as the difference of the lengths of $\gamma'$ and $\gamma$ that is \begin{eqnarray*} \mathcal{H}^{-}({i}) &=& |\gamma' \backslash (\gamma\cap \gamma')| - |\gamma \backslash  (\gamma\cap \gamma')|.  \end{eqnarray*}
The quantity $ \mathcal{H}^{+}(i)$ is defined similarly using the other endpoints of $ a_{i}$ and $e$.

\paragraph{Splitting the variation.} One might think that during an exploration process of $T_{\infty,\infty}$ the horodistances from $\pm \infty$ are only ruled by the peeling forms $(F_{i})_{i \geq 0}$ of the exploration. This is not true since after peeling the edge $a_{i}$, the next edge $a_{i+1}$ to peel can be located anywhere on the boundary of $K_{i+1}$ and could thus introduce a change in the horodistances. For convenience, we will thus consider intermediate half-integer steps in the horodistance processes $ \mathcal{H}^\pm$ which take into account only the variation of the horodistances due to the peeling steps. 

Specifically, we introduce the following functions of forms: for every $f \in \mathsf{Forms}$
  \begin{eqnarray} \label{eq:defdeltapm} \Delta^-(f) &=&  \frac{1}{2} \cdot  \mathbf{1}_{f=( \mathrm{C},1)} - \sum_{k \geq 1} k \cdot  \mathbf{1}_{f = ( \mathrm{G},-k)}, \\
  \label{eq:defdeltapm2}
\Delta^+(f) &=&   \frac{1}{2} \cdot  \mathbf{1}_{f=( \mathrm{C},1)} - \sum_{k \geq 1} k \cdot  \mathbf{1}_{f = ( \mathrm{D},-k)}.  \end{eqnarray}
One can clearly recover $f$ using the pair $( \Delta^-(f), \Delta^+(f))$. Then for every $i \geq 0$ we set 
$$ \mathcal{H}^{\pm}\left( i + \frac{1}{2}\right) = \mathcal{H}^{\pm}(i) + \Delta^{\pm}( F_{i}). $$
In particular,  when the exploration is Markovian then $\mathcal{H}^{\pm} (i+\frac{1}{2})-\mathcal{H}^{\pm} (i)$ are i.i.d. of law $\Delta^{ \pm}( \mathsf{F})$. Geometrically if $ F_{i} \ne ( \mathrm{C},1)$ then $ \mathcal{H}^{\pm}(i+\frac{1}{2})$ corresponds to horodistances of the only edge of the revealed triangle in $T_{i+1}$, thus $ \mathcal{H}^\pm(i+1)$ would be equal to $ \mathcal{H}^\pm(i+\frac{1}{2})$ if the next edge to peel would be that one. However, when $ F_{i} = ( \mathrm{C},1)$ the quantities $ \mathcal{H}^{\pm}(i+\frac{1}{2})$ do not represent actual horodistances (since they are half-integers) but an ``imaginary horodistance'' of an edge sitting in between of the two edges of the revealed triangle in $T_{i+1}$.  The quantity $$ \eta^\pm(i) := \mathcal{H}^{\pm} (i+1)-\mathcal{H}^{\pm} \left(i+\frac{1}{2}\right)$$ thus corresponds to the difference of the new edge with respect to the ``predicted'' next edge to peel and heavily depends on the algorithm chosen for the exploration. When $F_{i} \ne ( \mathrm{C},1)$ then $\eta_{i}^\pm \in  \mathbb{Z}$ and $\eta_{i}^\pm \in \mathbb{Z} + \frac{1}{2}$ otherwise. Besides we always have 
 \begin{eqnarray} \label{eq:sum=0} \eta^+(i) + \eta^-(i) &=&0, \qquad \mbox{for every }i \geq 0.  \end{eqnarray}

\paragraph{Minimum process.} Finally, we will use an important geometric quantity that can be read from the horodistance process. Recall that in this work, we always peel on the boundary of the explored part so that $(K_{i})$ is a growing triangulation with simple boundary. For any $i \geq 0$ we introduce the infimum process of the horodistance along half-integer times :

 \begin{eqnarray*} \underline{\mathcal{H}}^+(i) = \inf \left\{ \mathcal{H}^+ \left(j + \frac{1}{2}\right) : j \leq i \right\} \wedge 0 \quad \mbox{ and }\quad \underline{\mathcal{H}}^-(i) = \inf \left\{ \mathcal{H}^- \left(j + \frac{1}{2}\right) : j \leq i \right\} \wedge 0.  \end{eqnarray*}
It is easy to see by induction that $ -\underline{\mathcal{H}}^+(i)$ (resp.\,\,$ -\underline{\mathcal{H}}^-(i)$) can be interpreted as the number of edges of $ T_{\infty,\infty}$ on the right (resp.\,\,left) of the root edge $ \vec{e}$ that have been swallowed in $K_{i}$ so far. For example on Fig.\,\,\ref{fig:horo} we have $ \underline{ \mathcal{H}}^-(i) = -4$ and $ \underline{ \mathcal{H}}^+(i)=-2$.  In particular, at time $i \geq 0$ the exploration process discovers a new triangle of form $( \mathrm{D}, \cdot)$ and such that the third vertex of this triangle is lying on the original boundary of $T_{\infty,\infty}$ if and only if we have 
   \begin{eqnarray} \mathcal{H}^+\left( i + \frac{1}{2}\right) &=& \underline{ \mathcal{H}}^+(i) \label{eq:hithor},  \end{eqnarray} and similarly for the left-hand side with ``$+$'' replaced by ``$-$''.
 
% \begin{figure}[!ht]
%  \begin{center}
%  \includegraphics[width=12cm]{horomin}
%  \caption{ \label{fig:minhoro} The running infimums (along half-integer times) of the horodistance processes have a geometrical meaning. A new minimum of $ \mathcal{H}^+$ corresponds to a visit of the right boundary by the exploration process.}
%  \end{center}
%  \end{figure}

%%%%%%%%%%%%%%%%%%%%SUBSECTION
\subsection{The spectrally negative $ \frac{3}{2}$-stable process}
 Using the exact expression of  the probabilities $(q_{k})$ defined in Section \ref{sec:one-step} one sees that the random variables $ \Delta^\pm(\mathsf{F})$ defined by \eqref{eq:defdeltapm} and \eqref{eq:defdeltapm2} are bounded above by $1/2$ and satisfy 
   \begin{eqnarray} E[ \Delta^\pm( \mathsf{F})] = 0 \quad \mbox{ and } \quad P( \Delta^\pm( \mathsf{F}) =-n)  \underset{n \to \infty}{\sim} \frac{n^{-5/2}}{4 \sqrt{ \pi}}. \label{eq:domstable}  \end{eqnarray}
In other words,   $ \Delta^+( \mathsf{F})$ and $ \Delta^-( \mathsf{F})$ are both in the domain of attraction of the totally asymmetric stable random variable of parameter $ \frac{3}{2}$. Let us recall some basic facts   about the standard $ \frac{3}{2}$-stable spectrally negative Lévy  process  (with no positive jumps) with no drift which will be denoted by $(S_{t} :t \geq 0)$ and simply referred to as the ``$3/2$-stable process'' in the rest of this paper. We refer to \cite{Ber96} for details.  By standard we mean that the process $S$ satisfies $E[\exp(\lambda S_{t}) ] = \exp( t \lambda^{3/2})$ for all $\lambda>0$ or equivalently its Lévy measure is given by   \begin{eqnarray*} \Pi( dx) & = & \frac{3}{4 \sqrt{ \pi}}|x|^{-5/2} dx \mathbf{1}_{x <0}.  \end{eqnarray*} 
This process enjoys  the scaling property with parameter $3/2$ that is $(S_{t} : t \geq 0) = (\lambda ^{-2/3} S_{\lambda t} : t \geq 0)$ in distribution for any $\lambda >0$. 

\begin{figure}[!ht]
 \begin{center}
 \includegraphics[width=6cm]{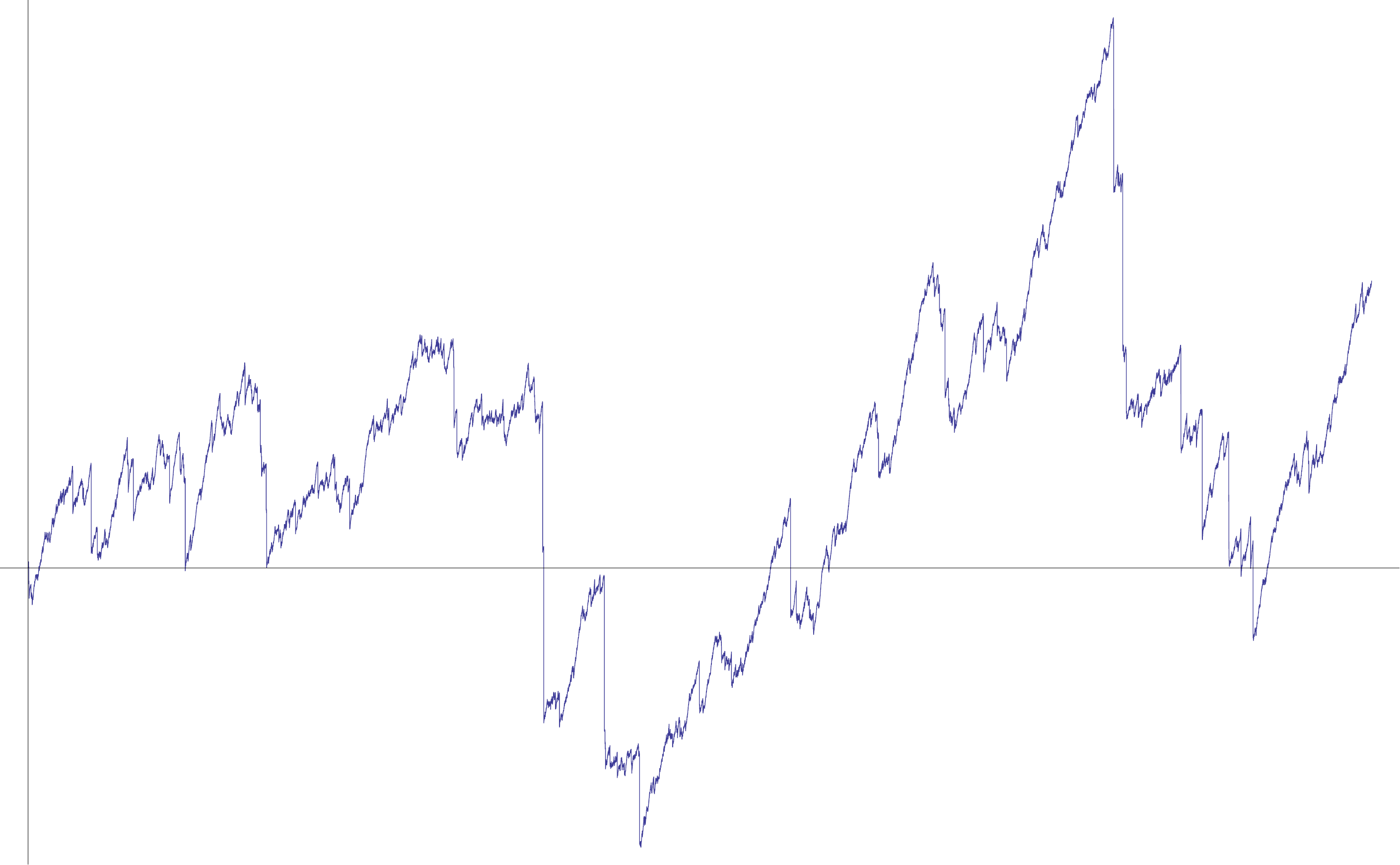} \hspace{1cm}\includegraphics[width=6cm]{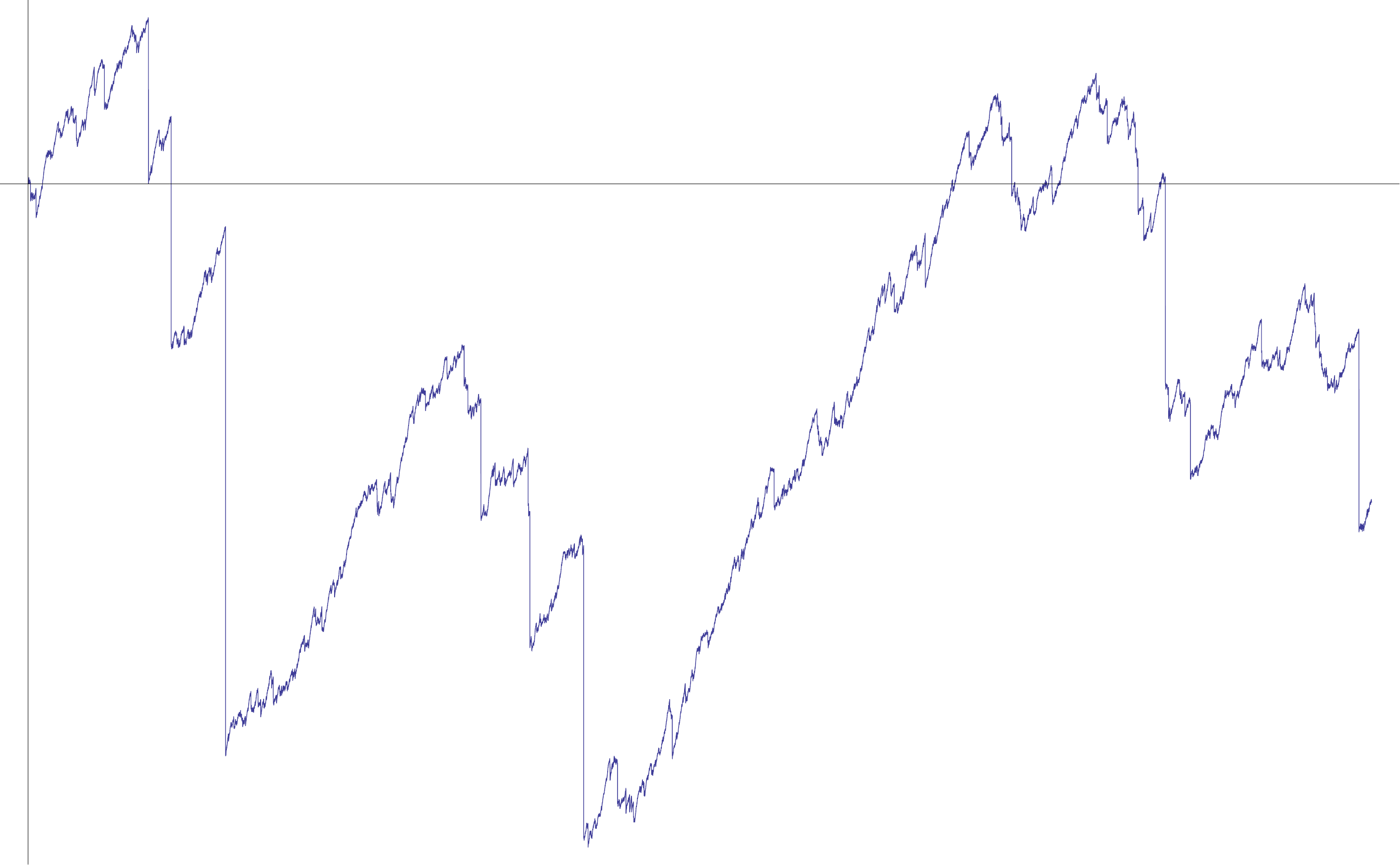}
 \caption{Two (approximated) samples of the process $S$.}
 \end{center}
 \end{figure}

The process $S$ will appear in this work as the scaling limit of discrete walks. Recall that if $\xi$ is a centered probability distribution over $ \mathbb{R}$ with increments bounded from above and such that $P( \xi \leq {-k}) \sim c k^{-3/2}$ as $k \to \infty
$, if $X_{1}, \ldots , X_{n}, \ldots$ are i.i.d. copies of $\xi$  with cumulative sum $Y_{n} = X_{1} + \cdots + X_{n}$ then we have the following convergence in distribution in the sense of Skorokhod 

  \begin{eqnarray} \label{eq:conv3-2}  \left(\frac{Y_{[nt]}}{(Kn)^{2/3}} \right)_{t \geq 0} &  \xrightarrow[n\to\infty]{(d)} & ({S}_{t})_{t \geq 0},  \end{eqnarray} with $K = 2c \sqrt{ \pi}$ and where $[x]$ denotes the largest integer less than or equal to $x$, see \cite{JS03}.\bigskip

  \begin{proposition}  \label{prop:deuxstables} If $(F_{i})_{i \geq0}$ are i.i.d.\ random variables distributed as $ \mathsf{F}$ then we have the following convergence in the sense of Skorokhod
   \begin{eqnarray*} n^{-2/3}  \cdot \left( \sum_{i=0}^{ [nt]} \Delta^+(F_{i}),\sum_{i=0}^{ [nt]} \Delta^-(F_{i}) \right)_{t \geq 0} &  \xrightarrow[n\to\infty]{(d)} & 3^{-2/3} \cdot \big(S^{+}_{t}, S^{-}_{t}\big)_{t \geq 0},  \end{eqnarray*}where $S^{+}$ and $S^{-}$ are independent standard $ \frac{3}{2}$-stable processes with no positive jumps. \end{proposition}
   \proof Although the variables $\Delta^+(F_{i})$ and $\Delta^-(F_{i})$ are not exactly independent, this is more or less an easy consequence of \eqref{eq:conv3-2}. Let us provide the details. To gain independence we Poissonize time. More precisely, we give us a Poisson clock of parameter $1$ and at each time $(s_{i}, i \geq 1)$ the clock rings, we sample a form $F_{{i}}$ according to $ \mathsf{F}$. Equivalently, every form $f \in \mathsf{Forms}$ appears with an independent Poisson clock of parameter $P( \mathsf{F}=f)$. For $t\geq 0$ let
   $$G_{t} = \sum_{\begin{subarray}{c}F_{{i}} = ( \mathrm{G},-k)\\ s_{i}\leq t\end{subarray}}-k \qquad \mbox{and}  \qquad D_{t} = \sum_{\begin{subarray}{c}F_{{i}} = (\mathrm{D},-k)\\ s_{i}\leq t \end{subarray}}-k $$ respectively be the sums of the (negative) jumps of left and right forms. We also put $C_{t}$ for  the number of centered forms  $(\mathrm{C},1)$ appeared before time $t$ (which is thus a Poisson variable of parameter $2t/3$). Then we have
    \begin{eqnarray} \label{eq:trivail} \left(\sum_{s_{i}\leq t} \Delta^-(F_{i}),\sum_{s_{i} \leq t} \Delta^+(F_{i}) \right)_{t \geq 0} = \left( G_{t} , D_{t}\right)_{t \geq 0} + \frac{1}{2}(C_{t},C_{t})_{t \geq 0}.  \end{eqnarray}
On the one hand, by Donsker's theorem, we have $$\left(\frac{C_{nt}- \frac{2}{3}nt}{\sqrt{n}}\right)_{t \geq 0}  \quad  \xrightarrow[n\to\infty]{(d)} \quad  (B_{t})_{t \geq 0},$$ where $B$ is a (multiple of a) Brownian motion. On the other hand, since $G_{t}$ and $D_{t}$ are now independent, by \eqref{eq:domstable} and the fact that $\sum_{k\geq1} kq_{-k} = 1/3$ we have
  \begin{eqnarray} \label{eq:poisson}\left(\frac{G_{nt} + \frac{nt}{3}}{n^{2/3}},\frac{D_{nt} + \frac{nt}{3}}{n^{2/3}}\right)_{t \geq 0} \quad \xrightarrow[n\to\infty]{(d)}  \quad  3^{-2/3}\cdot  (S^-_{t},S^+_{t})_{t \geq 0}, \end{eqnarray}  in distribution for the Skorokhod topology where $S^-$ and $S^+$ are independent standard $\frac{3}{2}$-stable processes with no positive jumps. Remark now that the last display holds if we replace $ \frac{nt}{3}$ by $ C_{nt}/2$ since the $ \sqrt{n}$ fluctuations of $C_{nt}$ around $2n/3$ are crushed by the  $n^{2/3}$ renormalization. Using \eqref{eq:trivail} and a standard depoissonization argument, this implies the proposition. \endproof

%%%%%%%%%%%%%%%%SECTION
\section{SLE$_{6}$ on the half-plane UIPT}
The goal of this section is to explain how to discover a half-plane UIPT using an SLE$_{6}$ process and to prove that, under hypothesis $(*)$, this exploration is the continuous limit of the discrete critical percolation interface in an appropriate sense. To help the reader digest our argument, we first recall the results of Angel \cite{Ang05,Ang03} on site percolation interface in $T_{\infty,\infty}$ using our formalism. We refer to \cite{ACpercopeel} for more details.

%%%%%%%%%%%%%%%%%%%%SUBSECTION
\subsection{Percolation exploration} \label{sec:perco}
Let $ T_{\infty,\infty}$ be the half-plane UIPT. Conditionally on $T_{\infty,\infty}$ we color each vertex of the triangulation independently white or black with equal probability, except for the vertices of the boundary : color in white those on the right of the root edge and in black those on the left. See Fig.\,\ref{fig:exploration} below. 
\begin{figure}[!ht]
 \begin{center}
 \includegraphics[width=16cm]{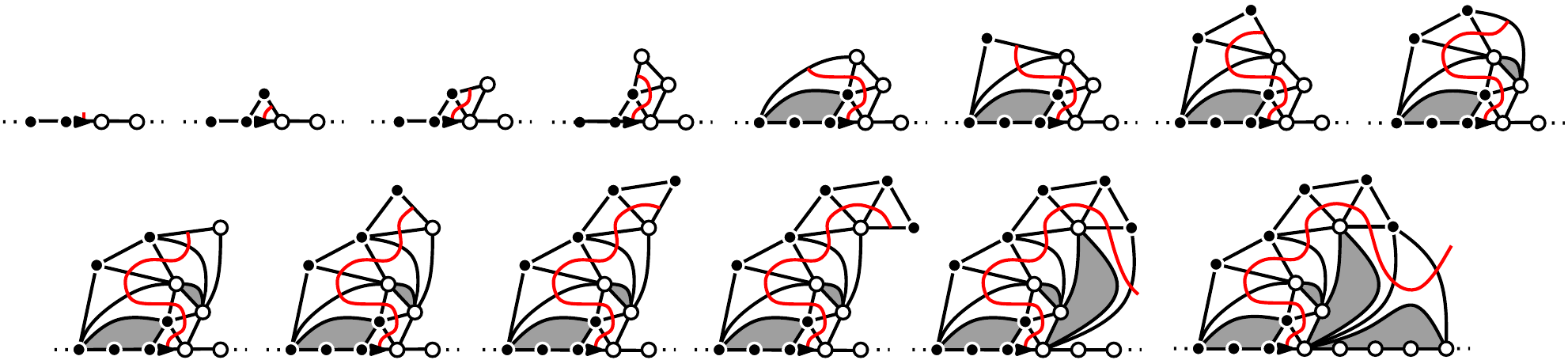}
 \caption{ \label{fig:exploration}Site percolation exploration on the UIHPT.}
 \end{center}
 \end{figure}
 
 It is possible to use the spatial Markov property of the UIHPT in order to discover step-by-step the percolation interface: at each step we reveal the triangle of the current boundary that lies between the black and the white component. If this triangle discovers a new vertex then reveal its color as well. It is easy to see  that if this algorithm has been used since the beginning, then at each step there is a unique edge $ \bullet - \circ$ located on the current boundary, see Fig.\,\ref{fig:exploration} above. This defines a \emph{Markovian} exploration of the UIHPT, see \cite{Ang05,Ang03,ACpercopeel}. If we denote by $ \mathcal{H}_{ \mathrm{Perc}}^{-}(n)$ and $ \mathcal{H}^{+}_\mathrm{Perc}(n)$ the horodistances of the edge $a_{n}$ at the $n$th step of peeling  then we have 
 
 \begin{theorem}[Angel \cite{Ang05}] \label{thm:exploperco}	We have the following convergence in distribution  \begin{eqnarray*}\left( \frac{\mathcal{H}_\mathrm{Perc}^{-}([nt])}{n^{2/3}} , \frac{\mathcal{H}_\mathrm{Perc}^{+}([nt])}{n^{2/3}} \right) & \xrightarrow[n\to\infty]{(d)} & 3^{-2/3} \cdot  (S^-_{t}, S^+_{t})_{t \geq 0}  \end{eqnarray*} where $S^-$ and $S^+$ are independent standard $ \frac{3}{2}$-stable processes with no positive jumps.\end{theorem}
 \proof[Proof (Sketch)] For every $i \geq 0$, if the revealed face at time $i$ is of the form $( \mathrm{C},1)$ then we let $\epsilon_{i} =\frac{1}{2}$ when the revealed vertex is white and $\epsilon_{i}=-\frac{1}{2}$ when it is black.   We set $ \epsilon_{i}=0$ otherwise. The description of the exploration process shows that for every $ i \geq 0$ we have 
$$  \begin{array} {rcr} 
 \eta_{ \mathrm{Perc}}^-(i)=\mathcal{H}^{-}_{ \mathrm{Perc}}(i+1)-\mathcal{H}^{-}_{ \mathrm{Perc}}(i+ \frac{1}{2})&=&  -\epsilon_{i},\\ \ \\
 \eta_{ \mathrm{Perc}}^+(i)=\mathcal{H}^{+}_{ \mathrm{Perc}}(i+1)-\mathcal{H}^{+}_{ \mathrm{Perc}}(i + \frac{1}{2})&=& + \epsilon_{i}.\end{array}$$
Indeed, when the revealed face is not of the form $( \mathrm{C},1)$ then $\epsilon_{i}=0$ and the next edge to peel is necessarily the unique edge of the revealed triangle belonging to the new infinite boundary (this edge is easily seen to be of type $ \bullet - \circ$). However, when the revealed face is of form $ ( \mathrm{C},1)$ then it has two edges belonging to the new infinite boundary and the next edge to peel is either the ``left'' edge of the revealed triangle if $ \epsilon_{i} = \frac{1}{2}$ or  its ``right'' edge if $ \epsilon_{i} = - \frac{1}{2}$.  Since the variables $ \epsilon_{i}$ are i.i.d. bounded centered variables we deduce that 
$$ \frac{\epsilon_{1} + \cdots + \epsilon_{[nt]}}{\sqrt{n}}  \quad \xrightarrow[n\to\infty]{(d)} \quad  (B_{t})_{t \geq 0},$$ where $B$ is  a multiple of a Brownian motion. On the other hand, since the exploration is Markovian the increments of the horodistances between $i$ and $i+\frac{1}{2}$ are independent copies of  $( \Delta^{-}( \mathsf{F}), \Delta^{+}( \mathsf{F}))$. We can thus combine the last display with Proposition \ref{prop:deuxstables} to get the desired result (notice again that the $\sqrt{n}$ scaling of the Bernoulli variables is crushed by the $n^{2/3}$ renormalization as in the proof of Proposition \ref{prop:deuxstables}). \qedhere

%%%%%%%%%%%%%%%%%%%%SUBSECTION
\subsection{The Riemann surface construction}
In this section we show how to associate with the UIHPT a Riemann surface that we will later use to define SLE processes on $T_{\infty,\infty}$. We follow the presentation of \cite{GR10} where the authors showed that the Riemann surface associated to the UIPT (of the full plane) is conformally equivalent to $ \mathbb{C}$. \bigskip

\label{sec:riemannconstruction}
We associate with any locally finite triangulation $T$ a Riemann surface $ \mathbf{T}$ by considering each triangle of the map as a standard Euclidean equilateral triangle endowed with its distance and use the combinatorics of the map to glue the triangles between each other. Formally, we first construct a topological space by gluing triangles according to the pattern of the map; this topological space is then endowed with a Riemann surface structure using the following coordinate charts:
\begin{itemize}
\item for any point located in the interior of a triangle or on a boundary edge we simply see this triangle as a standard equilateral triangle (whose sides have length $1$) in the complex plane and use the identity map,
\item if the point belongs to an interior edge, then place the two adjacent triangles (there must be two different triangles since we are considering type II triangulations) next to each other in the complex plane and use again the identity map,
\item if the point is located on an interior vertex with $d\geq 2$ adjacent equilateral triangles $t_{1}, t_{2}, \ldots , t_{d}$ arranged in cyclic order then we use the map $z \mapsto z^{6/d}$ as coordinate, that is, the point $z= r e^{ \mathrm{i} t}$ for $t \in [0,\pi/3]$ and $r < 1/2$ belonging to the triangle $t_{j}$ is sent to $\big(r\exp(\mathrm{i}(t+(j-1)\pi/3))\big)^{6/d}$.
\item If the point is a vertex on the boundary we modify the above chart by using $z \mapsto z^{3/d}$.
\end{itemize}
It is easy to check that the coordinate changes are analytic and thus this atlas does define a Riemann surface structure (in fact a  Euclidean surface with conical singularities at vertices of degree different from $6$), see \cite{GR10} for details.  \bigskip

In the case of the UIHPT we obtain a (random) simply connected Riemann surface with a boundary denoted by $\mathbf{T}_{\infty,\infty}$. By the uniformization theorem, this surface can be mapped onto the upper half-plane $ \mathbb{H} = \{ z \in \mathbb{C} , \Re(z)>0\}$, i.e.\,\,there exists a (random) bi-holomorphic function  $  \phi_{ \mathbf{T}_{\infty,\infty}} : \mathbf{T}_{\infty,\infty} \to \mathbb{H}$. This map is unique provided that we fix the images of three points : the origin of the root edge is sent to $-1/2$, its target to $1/2$ and the infinity of $ \mathbf{T}_{\infty,\infty}$ is sent to the infinity of $ \mathbb{H}$. The image of the edges of $T_{\infty,\infty}$ in $ \mathbf{T}_{\infty,\infty}$ under this conformal map is thus a canonical proper embedding of $T_{\infty,\infty}$ in $ \mathbb{H}$ and is denoted by $ \mathscr{T}_{\infty,\infty}$, see Fig.\,\ref{fig:uniformization}. \medskip

Once we have constructed this canonical representation of the UIHPT, one can consider various stochastic processes on it. For example we can \emph{define} a Brownian motion (up to time parametrization) moving over $ T_{\infty,\infty}$ (more precisely over $ \mathbf{T}_{\infty,\infty}$) as the pre-image under $\phi_{ \mathbf{T}_{\infty,\infty}}$ of a standard reflected Brownian motion on $ \mathbb{H}$. The goal of the next subsection is to study one very special random process over $ \mathbf{T}_{\infty,\infty}$  : the SLE process of parameter $\kappa=6$.

\begin{remark} They are various ways to construct a canonical embedding of a planar map, see \cite{Bef08}. However, we work here with Riemann's uniformization because it is well-suited to define and use the SLE$_{6}$ exploration (see below).
\end{remark}

%%%%%%%%%%%%%%%%%%%%SUBSECTION
\subsection{SLE$_{6}$ exploration}
\label{sec:defSLE}
 We recall the definition of the chordal $ \mathrm{SLE}_{6}$ in the upper half-plane. The reader is referred to \cite{Law05,Wer04} for details and proofs. Let $B_{t}$ be a standard linear Brownian motion and consider the flow of conformal mappings obtained by  solving the following PDE: \begin{eqnarray} \label{def:SLE6} \partial_{t} g_{t}(z) = \frac{2}{g_{t}(z) - \sqrt{6} B_{t}}, \qquad g_{0}(z) = z.  \end{eqnarray}
For each $t \geq 0$, the function $g_{t}$ maps a certain simply connected domain $H_{t} \subset \mathbb{H}$ onto the upper half-plane $ \mathbb{H}	$. Furthermore, it is by now classical that $\mathbb{H}\backslash H_{t}$ can be represented as the hull of a random curve $ \gamma :  \mathbb{R}_{+} \to \mathbb{H}$ starting from $0$, that is   \begin{eqnarray*} H_{t} = \mbox{infinite open component of } \mathbb{H} \backslash \gamma([0,t]).  \end{eqnarray*}
This curve is called the Schramm--Loewner curve of parameter $\kappa=6$ and abbreviated by $ \mathrm{SLE}_{6}$. %It enjoys the following scaling property: $ (\gamma_{t} : t\geq 0) = \lambda^{-1/2}( \gamma_{\lambda t} : t \geq 0)$ in distribution for all $\lambda >0$. 
For $\kappa > 4$, the $ \mathrm{SLE}_{6}$ is not a simple curve (it touches itself) and furthermore bounces on the real axis infinitely many often (this will be crucial in the sequel). 

\begin{figure}[!ht]
 \begin{center}
 \includegraphics[width=10cm]{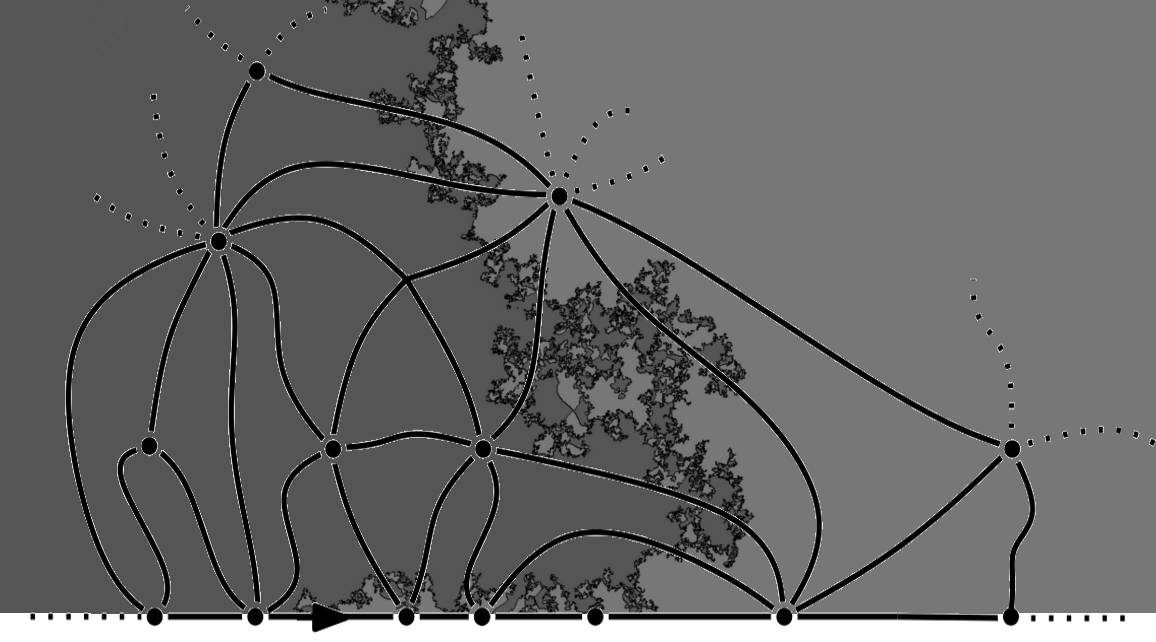}
 \caption{Simulation of an SLE$_{6}$ in the half-plane, courtesy of Vincent Beffara. On top of it, the uniformization of a UIHPT.}
 \end{center}
 \end{figure}

Independently of $T_{\infty,\infty}$, consider a standard SLE$_{6}$ curve  $(\gamma_{t})_{t \geq 0}$ on $ \mathbb{H}$ started from $0$. We \emph{define} the SLE$_{6}$ on the (Riemann surface associated to the) half-plane UIPT as the path  \begin{eqnarray} \label{eq:defSLE}  \big( \phi_{ \mathbf{T}_{\infty,\infty}}^{-1}( \gamma_{t}) \big)_{t \geq 0}.  \end{eqnarray} Although this process runs over the Riemann surface $ \mathbf{T}_{\infty,\infty}$, one will abuse notation and say that the SLE$_{6}$ explores the UIHPT itself and that $\gamma$ is running directly over $T_{\infty,\infty}$. One can thus make sense of the discrete notion of face, edge or point of $T_{\infty,\infty}$ visited by the SLE$_{6}$. \bigskip 

\begin{center} \textit{In the following, ``exploration of the UIHPT'' will always refer to the above SLE$_{6}$ exploration.}\end{center}

Let us begin with a few remarks concerning this process. Since the points are polar sets for the SLE$_{6}$ on $ \mathbb{H}$, it follows  that the curve $\gamma$ on the UIHPT almost surely does not visit the vertices of $T_{\infty,\infty}$ (recall that the root edge is uniformized onto $[- \frac{1}{2}, \frac{1}{2}]$). The SLE$_{6}$ defines an (\textit{a priori} non-Markovian) exploration of the half-plane UIPT:  For any $t \geq 0$, we denote by $ \mathrm{Hull}(t)$ the subtriangulation of $T_{\infty,\infty}$ obtained as the union of all the faces visited by the curve $\gamma$  before time $t$ as well as the finite regions they enclose. The growing subtriangulations $ \{\mathrm{Hull}(t)\}_{t \uparrow}$ are then naturally associated with an (\textit{a priori} non-Markovian) exploration process of $T_{\infty,\infty}$. After forgetting the continuous time parametrization, we denote by $$(a_{i})_{i\geq 0}, \quad (K_{i})_{i \geq 0}, \quad (T_{i})_{i\geq 0}, \quad   (\mathcal{H}^{-}(i),\mathcal{H}^{+}(i))_{i \geq 0}$$ the sequences of peeled edges,  explored and remaining parts, and horodistances in this exploration. Note that we used only the curve $\gamma$ up to time parametrization to define this exploration. Finally, we denote by $ \mathcal{F}_{n}$ the $\sigma$-field generated by the knowledge of the  part $K_{n}$ at (the discrete) time $n$ as well as the curve $\gamma$ restricted up to the first visit of a face not in $K_{n}$ (whose tip is thus located on $a_{n}$).

%%%%%%%%%%%%%%%%%%%%SUBSECTION
\subsection{Locality of SLE$_{6}$ and the spatial Markov property}
\label{sec:locality}

Remark that one could have considered other explorations  (on the Riemann surface) of $T_{\infty,\infty}$ using different SLE$_{\kappa}$ curves by mimicking Definition \eqref{eq:defSLE}. However, the SLE of parameter $ \kappa=6$ plays a very special role since it defines a \emph{Markovian} exploration in the sense of Definition \ref{def:markovian}. Let us explain this crucial point in more details. \medskip

A characteristic property of the SLE$_{6}$ process that it shares with Brownian motion is the locality property, see \cite[Section 6.3]{Law05}.  This property is reminiscent of the fact that SLE$_{6}$ is the scaling limit of site percolation interface on the triangular lattice \cite{Smi01} and loosely speaking means that the SLE$_{6}$ curve does not feel the boundary of the domain it explores until it   touches it. A key consequence for us is the following : 
\begin{quote}\emph{Although the definition of the SLE$_{6}$ over $T_{\infty,\infty}$ given via \eqref{eq:defSLE} \emph{a priori} depends on the Riemann uniformization of the UIHPT,  the locality property enables us to define the curve $\gamma$ (up to time reparametrization) running over $ \mathbf{T}_{\infty,\infty}$  by discovering the UIHPT ``step-by-step'' and revealing only the parts necessary for the SLE$_{6}$ to displace.} \end{quote}

\begin{figure}[!ht]
 \begin{center}
 \includegraphics[width=16cm]{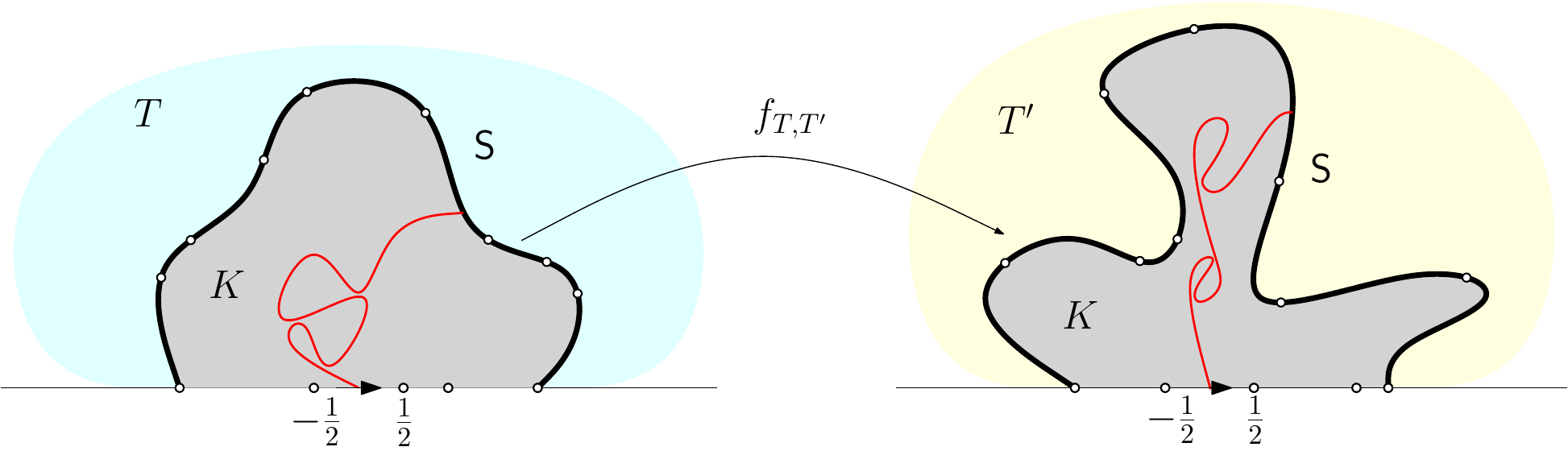}
 \caption{ \label{fig:locality} Rephrasing the locality property in our context.}
 \end{center}
 \end{figure}

More precisely, fix a finite triangulation $K$ with a simple boundary having a distinguished segment $ \mathsf{S}$ of boundary edges not containing the root edge; and  let $T$ be an infinite triangulation with an infinite boundary. We consider the triangulation $K + T$ obtained by gluing $T$ on the segment of $K$ and keeping the root of $K$. After uniformizing this triangulation onto $ \mathbb{H}$ as in Section \ref{sec:defSLE} we consider an independent SLE$_6$ curve $\gamma$ running on $K+T$. We denote by $ \gamma|_{{K}}$ the curve $\gamma$ seen up to time-reparametrization stopped at the first hitting of an edge of $ \mathsf{S}$, see Fig.~\ref{fig:locality}.
\begin{lemma} The law of $\gamma|_{ {K}}$ does not depend on $T$. In other words, the evolution of the SLE$_{6}$ inside $K$ can be performed without requiring the information outside $K$. 
\end{lemma}
\proof Consider two infinite triangulations with infinite boundary $T$ and $T'$.  After forming the two gluings of $K$ with $T$ and $T'$ along $ \mathsf{S}$, uniformize these two maps onto $ \mathbb{H}$ by sending the root edge to $[-1/2,1/2]$ and $\infty$ to $\infty$, see Fig.\,\,\ref{fig:locality}. 
In these uniformizations, the images of $ {K}$ thus form two different $\mathbb{H}$-neighborhoods %\footnote{an open subset $ \mathcal{N} \subset \mathbb{H}$ is a $ \mathbb{H}$-neighborhood of $x$ if $B(x, \varepsilon) \cap \mathbb{H} \subset \mathcal{N}$ for some $ \varepsilon>0$, see  \cite[Section 4.6]{Law05}}
 $N_{1}$ and $N_{2}$ of the origin (in gray in Fig.\,\,\ref{fig:locality}), see \cite[Chapter 6.3]{Law05}. The composition of the uniformizing maps thus yields a locally real\footnote{a univalent function $\phi : \mathcal{N} \to \mathbb{H}$ is locally real at $x_0$ if $\phi(z) = a_0+a_1(z-x_0)+a_2(z-x_0)^2+...$ locally around $x_0$ with $a_0,a_1,a_2,... \in \mathbb{R}$, see \cite[Section 4.6]{Law05}} conformal transformation $f_{T,T'}$ sending $N_{1}$ to $N_{2}$. The locality property of the SLE$_{6}$ \cite[Theorem 6.13]{Law05} precisely tells us that the image of an SLE$_{6}$ curve in $N_{1}$ has the law of an SLE$_{6}$ in $N_{2}$.  Otherwise said, the image of the SLE$_{6}$ running on the uniformization of $K+T$ and stopped when touching (the image of) $ \mathsf{S}$, once pushed by $f_{T,T'}$, is  an  SLE$_{6}$ exploring $K+T'$  stopped when touching $ \mathsf{S}$. The statement of the lemma follows.
\endproof

A repetitive use (left to the reader) of the last lemma shows that the edge to peel at time $i \geq 0$ is  independent of the remaining part $T_{i}$ and so:
\begin{corollary} \label{prop:markovianSLE6} The exploration process of $T_{\infty,\infty}$ induced by the SLE$_{6}$ is {Markovian}.
\end{corollary}

  \endproof

%%%%%%%%%%%%%%%%%%%%SUBSECTION
\subsection{The $(*)$ property} \label{sec:introstar}
The goal of this section is to introduce the  condition $(*)$ under which the scaling limit of the horodistance processes is known. Although we have been unable to prove it, we will try to convince the reader it is true, see Section \ref{sec:comments}.  Recall the notation
 \begin{eqnarray*}  
  \mathcal{H}^{-}(i+1)-  \mathcal{H}^{-}(i+\frac{1}{2})=  \eta^-_{i},\qquad \mbox{and} \qquad
\mathcal{H}^{+}(i+1)- \mathcal{H}^{+}(i+\frac{1}{2})=  \eta^+_{i}.\end{eqnarray*}
Our hypothesis $(*)$ on which most of the interesting results of this paper rely is 
\begin{center}\fbox{$(*) \qquad \qquad  \displaystyle \sup_{t \in [0,1]} \frac{\eta^+_{1}+\eta^+_{2}+ \cdots+ \eta^+_{[tn]}}{n^{2/3}}  \quad \xrightarrow[n\to\infty]{(P)} \quad  0$,} \end{center}
where $(P)$ denotes convergence in probability. 

\begin{theorem2} \label{thm:exploSLE} We have the following convergence in distribution   \begin{eqnarray*} \left(	\frac{ \mathcal{H}^{+}([2nt]/2)}{n^{2/3}}, \frac{  \mathcal{H}^{-}([2nt]/2)}{n^{2/3}}\right) & \xrightarrow[n\to\infty]{(d)} & 3^{-2/3} \cdot (S^+_{t},S^-_{t}) _{t \geq 0}  \end{eqnarray*} in the Skorokhod sense where $(S^+,S^-)$ is a pair of independent standard $ \frac{3}{2}$-stable processes. \end{theorem2}
\begin{remark} \label{rem:percoexplo}Theorem* \ref{thm:exploSLE} has to be compared with Theorem \ref{thm:exploperco}. In a weak sense, it says that the SLE$_{6}$ indeed is the scaling limit of critical percolation interfaces in the UIHPT (see \cite{Smi01} for the regular case), at least from the horodistances point of view. See Section \ref{sec:comments} for comments and open questions. \end{remark}

\proof[Proof*.] By Corollary \ref{prop:markovianSLE6} the SLE$_{6}$ exploration is Markovian and so by Proposition \ref{prop:markoviid} the variations of the horodistances between $i$ and $i+\frac{1}{2}$ are independent and distributed as $ ( \Delta^{-}(\mathsf{F}), \Delta^{+}( \mathsf{F}))$.  By Proposition \ref{prop:deuxstables}, we thus have 
 \begin{eqnarray*} n^{-2/3}\cdot\left( \sum_{i=0}^{[nt]} \mathcal{H}^{-}(i+ \frac{1}{2})- \mathcal{H}^-(i),\sum_{i=0}^{[nt]} \mathcal{H}^{+}(i+ \frac{1}{2})- \mathcal{H}^+(i)\right)_{t \geq 0} &\xrightarrow[n\to\infty]{(d)}& 3^{-2/3} \cdot (S^+_t,S^-_t)_{t \geq 0}. \end{eqnarray*} The condition $(*)$ then precisely entails that the increments between $i+1/2$ and $i+1$ cannot perturb the scaling limit. More precisely, since $ \eta_{i}^+= - \eta_{i}^-$ by \eqref{eq:sum=0}, condition $(*)$ implies
\begin{eqnarray*} n^{-2/3}\cdot\left( \sum_{i=0}^{[nt]} \mathcal{H}^{-}(i+1)- \mathcal{H}^-(i+ \frac{1}{2}),\sum_{i=0}^{[nt]} \mathcal{H}^{+}(i+ 1)- \mathcal{H}^+(i+\frac{1}{2})\right)_{t \geq 0} &\xrightarrow[n\to\infty]{}& 0, \end{eqnarray*} in probability for the Skorhokhod topology. Combining the last two displays yields to Theorem* \ref{thm:exploSLE} when $[2nt]/2$ is replaced by $[nt]$. To get the full statement, notice  that condition $(*)$ together with \eqref{eq:sum=0} also implies that $n^{-2/3} \cdot \sup_{i \leq n} \eta^\pm_i \to 0$ in probability (see also Proposition \ref{prop:boundetai} below for a stronger statement not depending on $(*)$). \endproof 

\begin{remark} At this point, the cautious reader may wonder why we have not chosen to explore the UIHPT using a Brownian motion instead of an SLE$_6$. Indeed, Brownian motion also enjoys the locality property and hence produces a Markovian exploration. The problem is that, contrary to the SLE$_6$, from time to time two consecutive peeling points for the Brownian motion may be far apart (in terms of horodistance): this occurs when the Brownian motion dive deep into the explored part so that the next peeling point is almost uncorrelated with the preceding one. Clearly, the analogous of condition $(*)$ for the Brownian exploration of $T_{\infty,\infty}$ does not hold and understanding the behavior of the horodistances, even on a heuristic level, is a very difficult problem.
\end{remark}

%%%%%%%%%%%%%%%%%%%%SUBSECTION
\subsection{Tail bound for the $\eta_{i}^\pm$}
Although the collective behavior of the $\eta_{i}^\pm$ is the content of the condition $(*)$ and remains conjectural, one can establish almost exponential bounds on the tails of the $\eta_{i}^\pm$. \bigskip

\begin{proposition}[Bounds for the $\eta_{i}^\pm$]  \label{prop:boundetai}  For $i \geq 0$, denote by $ \mathsf{D}_{i}$ the maximal degree of a vertex in $K_{i}$ within distance $2$ of its exposed boundary (that is the boundary in common with $T_{i}$). There exist some constants $c_{1},c_{2}>0$ such that for every $i \geq 0$ and every $k \geq 1$
\begin{enumerate}[(i)] 
\item $P(  \mathsf{D}_{i} \geq k) \leq i c_{1} \exp({-c_{2} k})$,
\item  conditionally on $ \mathcal{F}_{i}$, we have 
$ P( |\eta_{i}^\pm| \geq k) \leq c_{1}\exp({-c_{2}k  \mathsf{D}_{i}^{-4} }).$
\item Consequently we have for every $ \varepsilon>0$ $$ \lim_{n \to \infty} \frac{\sup_{i \leq n}|\eta_{i}^\pm|}{\log^{5+ \varepsilon} n} = 0. $$
\end{enumerate}
\end{proposition}

\begin{figure}[!ht]
 \begin{center}
 \includegraphics[width=12cm]{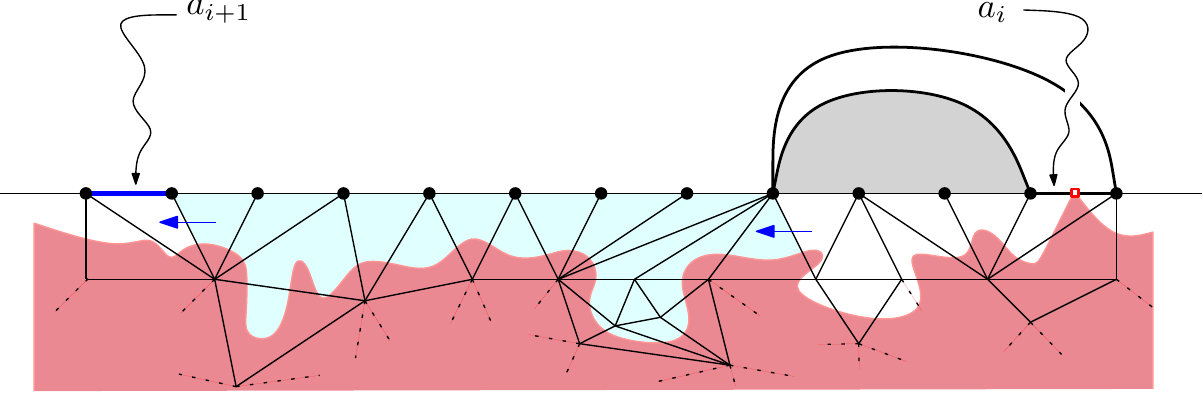}
 \caption{ \label{fig:etai}Illustration of the proof. The hull of the SLE is represented in red. To get to the desired next edge to peel, the SLE has to cross the channel and avoid the boundary edges.}
 \end{center}
 \end{figure}

\proof[Proof of Proposition \ref{prop:boundetai}] $(i)$. This statement should not be surprising for experts since it is more-or-less folklore that the maximal degree in a random triangulation is logarithmic in its size, see \cite[Lemma 4.2]{AS03} and \cite[Proposition 12]{BCsubdiffusive} for similar statements. However, we give a full proof for completeness. First of all, an easy adaptation of \cite[Lemma 4.2]{AS03} to the case of the UIHPT shows that the degree of the origin of the root edge in $T_{\infty,\infty}$ has an exponential tail. Actually, a slight generalization of it (left to the reader) shows that the maximal degree of a vertex within distance $3$ of the root edge $ {e}$ of $ T_{\infty,\infty}$ (that is  of one of its extremities)  also has an exponential tail, namely there exist $c_{1},c_{2}>0$ such that for every $k \geq 0$
  \begin{eqnarray} P\Big(\exists v \in T_{\infty,\infty}, \mathrm{d}_{ \mathrm{gr}}^{T_{\infty,\infty}}(v, {e})\leq 3 : \mathrm{deg}(v) \geq k\Big) \leq c_{1} e^{-c_{2}k},   \label{eq:expo4} \end{eqnarray} where the notation $ \mathrm{d}_{ \mathrm{gr}}^G$ stands for the graph distance in the graph $G$.   Next, for $i \geq 1$  if $v$ is a vertex of $K_{i}$ within distance $2$ of its exposed boundary, then if $j = \inf \{ n \leq i : v \in K_{n}\}$ is the first time at which $v$ is discovered (note that since $K_{0}= \varnothing$ we have $j\geq 1$), then $v$ is necessarily a vertex of $T_{j-1}$ and an easy geometric argument shows that 
  
  $$ \mathrm{d}_{ \mathrm{gr}}^{T_{j-1}}(v,a_{j-1}) \leq 3.$$  Recall from Proposition \ref{prop:markoviid} that for a Markovian exploration process, for every $j \geq 0$ the unexplored part $T_{j}$ rooted at $a_{j}$ is distributed as a standard UIHPT. By the union bound and  \eqref{eq:expo4} we thus have
 \begin{eqnarray*}P\big(  \mathsf{D}_{i} \geq k \big) &\leq & P\Big( \exists v \in T_{\infty,\infty}, \mathrm{d} ^{T_{j}}_{ \mathrm{gr}}(v,a_{j}) \leq 3 \mbox{ for some } 0 \leq j \leq i : \mathrm{deg}(v) \geq k \Big)\\ &\leq& i \  P \Big( \exists v \in T_{\infty,\infty}, \mathrm{d}_{ \mathrm{gr}}^{T_{\infty,\infty}}(v, {e})\leq 3 : \mathrm{deg}(v) \geq k \Big)\\
 & \underset{ \eqref{eq:expo4}}{\leq} & i c_{1}  e^{-c_{2} k}.
  \end{eqnarray*} 
$(ii)$ We only give  a detailed sketch and leave the precise details to the careful reader. Imagine the situation just after having peeled the $i$th edge. Let $k \geq 1$ (large) and let us evaluate the probability that the next edge to peel is the $k$th edge on the left of the root of the triangulation $T_{i+1}$. Since all the triangles that share an edge with the boundary of $T_{i}$ have been visited by the SLE$_{6}$ process, that means that the curve $\gamma$ has to travel in a narrow region towards the left to finally exit at the desired edge while bumping on its past and without touching any edge on the boundary of $T_{i+1}$ during its journey. See Fig.\,\,\ref{fig:etai}.
 
\begin{figure}[!ht]
  \begin{center}
  \includegraphics[width=14cm]{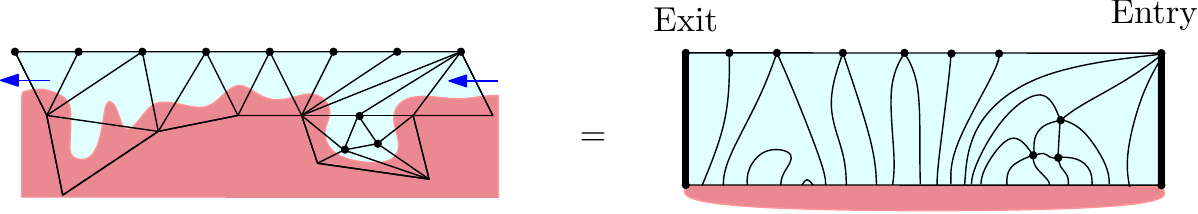}
  \caption{ \label{fig:beauchannel}The channel and its associated conformally equivalent rectangle (artistic representation).}
  \end{center}
  \end{figure}

 In particular, we can define a ``\emph{channel}'' (in light blue on Fig.\,\, \ref{fig:etai} and \ref{fig:beauchannel}) as being the region separating the target edge from the current position of the SLE with two edges playing the role of the entry and exit of the channel, see Fig.\,\ref{fig:beauchannel}. To show the bound of the proposition, we will prove that the probability that an SLE$_{6}$ crosses the channel without touching the above boundary edges is very low.  This is intuitively clear since the latter is a narrow and long path (when $k$ is large), but what really matters is its \emph{conformal} width. 

 More precisely, we consider the Riemann surface $ \boldsymbol{\mathcal{C}}$ associated to the channel made by the  parts of the triangles that are not contained in the hull of the SLE$_{6}$, see Fig.\ref{fig:beauchannel}. By the uniformization theorem, we can map  $ \boldsymbol{\mathcal{C}}$ onto a rectangle where the vertical sides correspond to the entry and exit of the channel. Then by standard properties, the probability that an SLE$_{6}$ process crosses such a rectangle without touching its above boundary is at most $c_{1}\exp( -c_{2} {L})$ where $c_{1},c_{2}>0$ and $ {L}$ is the ratio (which does not depend on the uniformization) of the horizontal length by the vertical length of the rectangle also called the extremal length or conformal moduli. The statement of the proposition thus reduces to show that the extremal length of the channel is at least   \begin{eqnarray} \label{eq:g2}  L \geq c \frac{k}{ \mathsf{D}_{i}^4},  \end{eqnarray} for some constant $c>0$. 

\begin{figure}[!ht]
  \begin{center}
  \includegraphics[width=7cm]{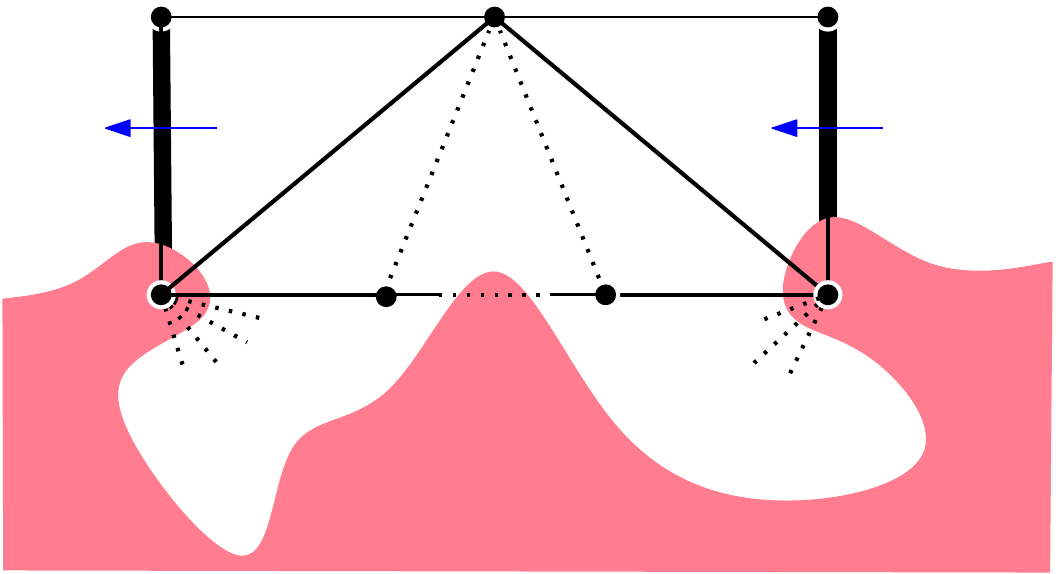}
  \caption{ \label{fig:block} A  block  requires a minimal $\rho_{0}$-length to be crossed.}
  \end{center}
  \end{figure}

  For this we use the definition of the extremal length of the channel $ \boldsymbol{\mathcal{C}}$ which is seen as a gluing of parts of equilateral triangles and thus endowed with the locally Euclidean metric and measure. If $\rho :  \boldsymbol{\mathcal{C}} \to \mathbb{R}_{+}$ is a positive function (also called ``metric'') we let $ \mathrm{Area}(\rho)$ be the integral of $\rho^2$ with respect to the Lebesgue measure on $ \boldsymbol{\mathcal{C}}$. Also, if $\Gamma$ is a smooth path going from the Entry to the Exit of the channel, we define the $\rho$-length of $\Gamma$ as
  $$ \mathrm{Length}_{\rho}( \Gamma) = \int_{\Gamma}  |ds| \ \rho,$$
  where $|ds|$ denotes the Euclidean element of length. With this piece of notation, the extremal length $L$ of $ \boldsymbol{\mathcal{C}}$ is expressed as (see \cite[Chapter 4]{Ahl73})
 \begin{eqnarray} \label{eq:eldef} 
  L = \sup_{ \rho}\inf_{\Gamma} \frac{ \big(\mathrm{Length}_{\rho}( \Gamma)\big)^2}{ \mathrm{Area}(\rho)},  \end{eqnarray} where the supremum is taken over all ``metrics'' $ \rho :  \boldsymbol{\mathcal{C}} \to \mathbb{R}_{+}$ and the infimum runs over all rectifiable paths joining Entry to Exit in the channel. To show  \eqref{eq:g2} we consider a particular metric $\rho_{0}$ defined as follows:  the function $\rho_{0}$ is constant and equals to $1$ on every (part of) triangle of $\boldsymbol{\mathcal{C}}$ which contains a vertex at  combinatorial distance less than $1$ from  the above boundary of the channel. Otherwise $\rho_{0}=0$ on the rest of the channel. Because $ \mathsf{D}_{i}$ is the maximum vertex degree within distance $2$ of the exposed boundary of $K_{i}$ we have 
   \begin{eqnarray} \mathrm{Area}(\rho_{0}) \leq  \frac{ \sqrt{3}}{4} \cdot (k+1) \cdot \mathsf{D}_{i}^2.  \label{eq:area} \end{eqnarray}
We now have to bound from below the $\rho_{0}$-length of a smooth path crossing $ \boldsymbol{\mathcal{C}}$. To do so, we will identify a combinatorial pattern in the channel that requires a minimal $\rho_{0}$-length to be traversed. First notice that  all the combinatorial triangles adjacent to the above boundary of the channel are either pointing upwards $\Delta$ or downwards $\nabla$. A   \emph{block} is a sequence $ \nabla, \Delta,\ldots, \Delta, \nabla$ together with the triangles ``grafted'' on the bottom of the upwards triangles. See Fig.\,\ref{fig:block}. 
  
As already mentioned, all the downwards triangles $\nabla$ contain a piece of the curve $\gamma$ for otherwise they would not have  been discovered. An easy geometrical argument shows that the $\rho_{0}$-length needed to cross a block is bounded from below by some universal constant $c'>0$. Hence, the minimal $\rho_{0}$-length of a curve $\Gamma$ crossing the channel is at least $c'$ times the number $B$ of blocks of this channel. However, it is easy to see that 
$$ B \geq \lfloor \frac{k}{10 \mathsf{D}_{i}} \rfloor,$$ which combined with \eqref{eq:area} and the definition \eqref{eq:eldef} of $L$ finishes the proof of the $(ii)$.

\noindent $(iii)$ For the final item we have using $(i)$ and $(ii)$
 \begin{eqnarray*} P( |\eta_{i}^+| \geq \log^{5+ \varepsilon} i)& \leq& P( \mathsf{D}_{i} \geq \log^{1 +{ \varepsilon}/{5}} i) + P(|\eta_{i}^+| \geq \log^{5 + \varepsilon} i \mid  \mathsf{D}_{i} \leq \log^{1 + \varepsilon/5} i) \\ & \leq& i c_{1} e^{-c_{2} \log^{1 + \varepsilon/5} i} + c_{1}e^{-c_{2} \log^{1 + \varepsilon/5} i}.  \end{eqnarray*}
The right-hand side is obviously summable in $i\geq 1$ and so an application of Borel--Cantelli's lemma finishes the proof of the proposition.\qedhere \bigskip

%%%%%%%%%%%%%%%%SECTION
\section{Bouncing off the walls}
\label{sec:bouncing}
The basic idea of Theorem* \ref{thm:main} is the following: When the horodistance $ \mathcal{H}^{+}$ (resp.\,\,$ \mathcal{H}^{-}$) reaches a new minimum value, this geometrically corresponds to a visit of $ \mathbb{R}_{+}$ (resp.\,$ \mathbb{R}_{-}$) by the SLE curve $\gamma$.  This heuristic is not exact on a discrete level but becomes true in the limit (see Proposition* \ref{prop:limitbouncing}). This enables us to relate the number of alternative visits to $ \mathbb{R}_{+}$ and $ \mathbb{R}_{-}$ by the curve $\gamma$ in terms of alternative minimal records of $S^+$ and $ S^-$ (Proposition \ref{prop:commute}).

%%%%%%%%%%%%%%%%%%%%SUBSECTION
\subsection{Discrete bouncing} \label{sec:discretebouncing}
For any $n \geq 0$, we introduce the first time $\tau_{+}(n)$ after $n$ such that the peeling of the edge $a_{\tau^{+}(n)}$ discovers a triangle of form $( \mathrm{D},\cdot)$ whose third vertex is lying on the original boundary of $T_{\infty,\infty}$.  Equivalently, using \eqref{eq:hithor} we have
 \begin{eqnarray} \label{eq:deftauplusn} \tau^+(n) &=& \inf \left\{ k \geq n : \mathcal{H}^+ \left( k + \frac{1}{2}\right) = \underline{ \mathcal{H}}^+(k)\right\}. \end{eqnarray}  The quantity $\tau^{-}(n)$ is defined by similar means. Thanks to Theorem* \ref{thm:exploSLE} and since $\liminf S^+ = \liminf S^- = -\infty$ we have $ \tau^+(n) < \infty$ and $ \tau^-(n) < \infty$ almost surely for every $n \geq 0$.
 
A peeling time $n$ is \emph{good} if the tip of the $SLE_{6}$ is located in the middle third of the edge to be peeled.

\begin{lemma}[Discrete bouncing]  \label{lem:discetebouncing}

There exists some constant $c>0$ such that on the event $\{ \mathcal{H}^+(n+1/2) = \underline{ \mathcal{H}}^+(n)\}$ and $n$ being a good peeling time, then conditionally on $ \mathcal{F}_{n}$ there is a probability at least $c$ that $\gamma$ touches $ \mathbb{R}_{+}$ within the next two peeling steps.  \end{lemma}

Obviously, a similar lemma holds when ``$+$'' is replaced by ``$-$''. 

\proof[Proof (Sketch).] Conditionally on the event considered, there is a probability bounded away from $0$ that the next peeling edge is good and is the left-most edge of the revealing triangle at time $n$ and that furthermore,  the peeling of that edge discovers a triangle ``glued'' on the boundary as in the following picture (see Remark 1). It is then easy to see that on this event the SLE$_{6}$ can touch $ \mathbb{R}_{+}$ with  a probability bounded away from $0$.%Fix $n\geq 0$ and suppose that $\tau^+(n) = N < \infty$. At time $N$ the peeled triangle of form $(\mathrm{D},\cdot)$ has thus  ``swallowed'' all the exposed boundary of $K_{N}$ on the right of $a_{N}$. If $(F_{i})_{i \geq 0}$ denote the forms of the triangles peeled, the event we now consider  is the following (see Fig.\,\ref{fig:hitRplus} and recall Fig.\,\ref{fig:swallowed}): $$ E_{N} = \left\{ \begin{array}{l} \eta^+_{N} = 0,  \eta^+_{N+1} = 0, F_{N+1} =  ( \mathrm{D},-1), F_{N+2} =  ( \mathrm{G},-1), \\ \mbox{and the associated enclosed triangulations of perimeter $2$ }\\ \mbox{discovered at time $N+1$ and $N+2$ are empty}\end{array}\right\}.$$
 
\begin{figure}[!ht]
 \begin{center}
 \includegraphics[width=15cm]{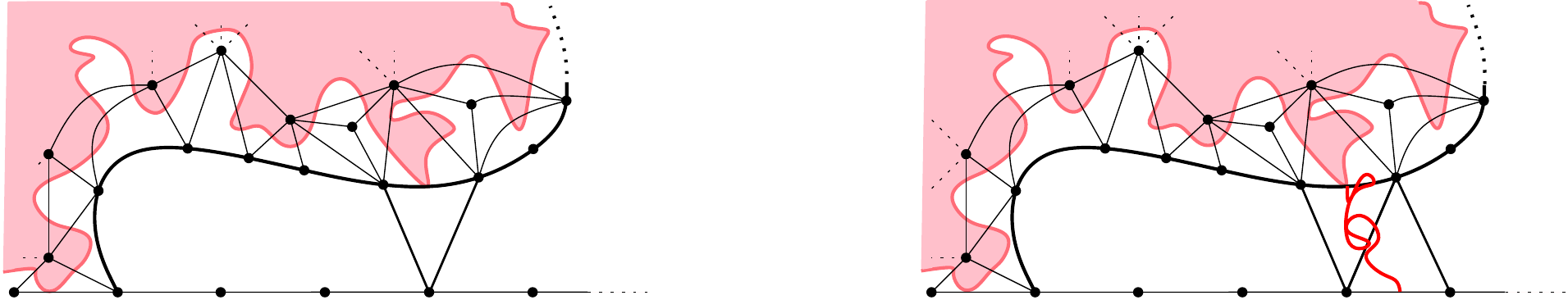}
 \caption{ \label{fig:hitRplus}Illustration of the proof.}
 \end{center}
 \end{figure}
% In words, on this event the SLE$_{6}$ first exists  towards the only exposed edge of the triangle discovered at time $N$ and discovers a triangle ``glued'' to the explored part on the right, then exists towards the exposed edge of this triangle and discovers a triangle glued to the explored part on  the left. It should be clear that conditionally on $ \mathcal{F}_{N}$, the last event has a probability bounded away from $0$ and that conditionally on $E_{N}$ the path $\gamma$ touches the real axis between time $N$ and $ N+3$ with positive probability. We leave the details to the reader. 
\qedhere

\paragraph{Commutings.}

We will now describe the limit as $n \to \infty$ of the random times $\tau^+(n)$ using the scaling limit of the horodistance processes given by Theorem* \ref{thm:exploSLE}. Recall that $S^+$ and $S^-$ are two independent  standard $ \frac{3}{2}$-stable processes with only negative jumps. We denote by $$\underline{S}^+_{t} = \inf\{ S^+_{u} : 0 \leq u \leq t\} \quad   \mbox{and} \quad  \underline{S}^-_{t} = \inf\{ S^-_{u} : 0 \leq u \leq t\}$$ be the running infimum processes of $S^+$ and $S^-$. For every $t \geq 0$ introduce 
 \begin{eqnarray*}  \xi^{+}(t) = \inf\big\{ u \geq t : S^{+}_{u} = \underline{S}^{+}_{u}\big\}, \end{eqnarray*}
and put a similar definition for $\xi^{-}(t)$. By standard properties of the spectrally negative $\frac{3}{2}$-stable process, for every $t >0$ we have $ \xi^+(t)>t$ almost surely. Furthermore the time  ${\xi^+(t)}$ a.s.~corresponds to a jump of the process which reaches a strict new minimum, that is    \begin{eqnarray} \label{eq:jumpmin}  \underline{S}^+_{t} > S^+_{\xi^+(t)}.  \end{eqnarray} Using standard properties of the Skorokhod topology \cite[Chapter VI]{JS03}, we deduce from the above display, \eqref{eq:deftauplusn} and Theorem* \ref{thm:exploSLE} that for every $t >0$ we have the following convergence in distribution 
 \begin{eqnarray} \label{eq:easy} \frac{\tau^+([nt])}{n} &\xrightarrow[n\to\infty]{(d)}& \xi^+(t),  \end{eqnarray} and similarly when ``$+$'' is replaced by ``$-$''. 
We denote by  $$ \mathcal{R}^+ = \{ t \geq 0 : S^+_{t} = \underline{S}^+_{t}\} \quad  \mbox{and} \quad  \mathcal{R}^- = \{ t \geq 0 : S^-_{t} = \underline{S}^-_{t}\}$$ the set of times corresponding to minimal records of the processes $S^+$ and $S^-$. These two random closed sets are a.s.\,\,perfect (not isolated points), also it is known that we almost surely have $$  \mathcal{R}^+ \cap \, \mathcal{R}^- = \{0\},$$
see \cite[Chapter 5]{Ber99}. The last display entails that for every $k \in \{ 0,1,2,\ldots\}$ the $k$th alternate composition 
$$\xi^{(k)} = \underbrace{\xi^\pm \circ \cdots \circ \xi^{+}\circ\xi^{-} \circ \xi^{+}}_{k \ \mathrm{terms}},$$
is well-defined and that we have $ t < \xi^{(1)}(t)< \xi^{(2)}(t) < \cdots$ as well as  $\xi^{(n)}(t) \to \infty$ for every $t >0$ as $n$ goes to infinity.   Note that for any $t>0$, by the scaling property of the stable processes we have the following equality in distribution   \begin{eqnarray} \label{eq:iddist}  \left(\xi^{(k)}(t)\right)_{k \geq 0} = t \cdot \left(\xi^{(k)}(1)\right)_{k \geq 0}.  \end{eqnarray}  We will later study the behavior of $\xi^{(n)}(1)$ as $n \to \infty$, see Proposition \ref{prop:commute}. In the spirit of \eqref{eq:jumpmin} one can check that $  \xi^{(k)}(t)$ is a jump time of $S^\pm$ (depending on the parity of $k$) that reaches a strict new minimum a.s. We mimic the definition of $\xi^{(k)}$ and set $\tau^{(k)}$ to be the $k$th alternate composition  $$\tau^{(k)} = \underbrace{\tau^\pm \circ \cdots \circ \tau^{+}\circ\tau^{-} \circ \tau^{+}}_{k \ \mathrm{terms}}.$$ The above considerations show that Theorem* \ref{thm:exploSLE} actually leads to the following extension of \eqref{eq:easy}: for every  $t >0$  we have the following convergence in distribution
   \begin{eqnarray} \label{eq:compose}  \left(\frac{\tau^{(k)}([nt])}{n} \right)_{k \geq 0}  &\xrightarrow[n\to\infty]{(d)}& \Big(\xi^{(k)}(t)\Big)_{k \geq 0}, \end{eqnarray} for the topology of simple convergence.    \medskip

We now introduce similar notions in order to describe the alternative bouncings of the SLE on $ \mathbb{R}_{+}$ and $ \mathbb{R}_{-}$. In the following lines, it is important to parametrize the SLE$_{6}$ and we recall from Section \ref{sec:defSLE} that $(\gamma_{t})_{t \geq 0}$ is a standard chordal  SLE$_{6}$ on $ \mathbb{H}$ starting from $0$ and parametrized by its half-plane capacity. In accordance to the above notation, for every $t \geq 0$ we put 
$$ \theta^+(t) = \inf\{ s \geq t : \gamma_{s} \in \mathbb{R}_{+}\},$$ where an obvious definition holds for $ \theta^-$. Here also, for every $k \geq 0$ we denote by $\theta^{(k)}$ the $k$th alternated composition 
$ \theta^\pm \circ \ldots \circ \theta^- \circ \theta^+$. Again, the scaling property of the SLE process implies that   \begin{eqnarray} \label{eq:iddist2} \left(\theta^{(k)}(t)\right)_{ k \geq 0} =  t \cdot \left( \theta^{(k)}(1) \right)_{ k \geq 0}  \end{eqnarray} in distribution for every $t >0$. Proposition \ref{prop:commuteSLE} studies the behavior of $ \theta^{(n)}(1)$  as  $n \to \infty$.

 When the SLE$_{6}$ curve $\gamma$ is used to explore the half-planar triangulation $T_{\infty,\infty}$ we will need to tie the continuous parametrization of the curve $\gamma$ to the discrete exploration steps. If $t_{i}$ for $i = 0,1,2,\ldots$ are the continuous times at which the $i$th edge to peel is discovered by the SLE process then for every $ k,i \in \{0,1,2,\ldots\}$ we let 
   $$ \tilde{\theta}^{(k)}(i) = \inf\big\{ j \geq 0 : \theta^{(k)}(t_{i}) \leq t_{j}\big\}.$$

  As promised in the introduction of this section, we prove that the scaling limit of the alternative bouncing on $ \mathbb{R}_{+}$ and $ \mathbb{R}_{-}$ by the SLE are described by the alternative minimal records of the processes $S^+$ and $ S^-$. More precisely, we have

\begin{proposition2}[Connecting $\theta$ and $\tau$] \label{prop:limitbouncing} For every $k \geq 1$  we have 
 \begin{eqnarray*}\frac{\tilde{\theta}^{(k)}(n)}{\tau^{(k)}(n)} & \xrightarrow[n\to\infty]{(P)} & 1. \end{eqnarray*}
\end{proposition2} 

\proof 
\textsc{Lower bound.} Assume that at time $t \geq 0$ we have $\gamma_{t} \in \mathbb{R}_{+}$. Intuitively the next edge to peel will be close to the extreme-right edge of the explored part, that is with a minimal horodistance. Indeed, an easy adaptation of the proof of Proposition \ref{prop:boundetai} shows that the next edge to peel $a_{i}$ has a horodistance $ \mathcal{H}^+(i)$ close to $ \underline{ \mathcal{H}}^+(i-1)$  in the sense that asymptotically we have
$$  \frac{\mathcal{H}^+(i)-\underline{ \mathcal{H}}^+(i-1)}{ \log^{5+ \varepsilon} i} \leq 1.$$
Since by Theorem* \ref{thm:exploSLE},  the quantity $ \underline{\mathcal{H}}^+(i-1)$ is of order $ i^{2/3}$ that means that $ \mathcal{H}^+(i)$ is very close to its past infimum. Using the fact that the set of minimal records of a $3/2$-stable process has no isolated point and standard properties of stable processes,  Theorem* \ref{thm:exploSLE} implies that for any $ \varepsilon>0$ with high probability there exists $ (1- \varepsilon) \leq j \leq (1+ \varepsilon)i$ such that $ \underline{ \mathcal{H}}^+(j) = \mathcal{H}^+(j+1/2)$. Iterating this argument we get that  \begin{eqnarray*} P\Big( \tau^{(k)}(n) \leq (1+ \varepsilon)\tilde{\theta}^{(k)}(n)\Big) & \xrightarrow[n\to\infty]{} & 1, \end{eqnarray*} for any $ k \in \{ 1,2, \ldots\}$ and any $ \varepsilon >0$. 

\textsc{Upper bound.} Fix $n \geq 0$ (large). By Lemma \ref{lem:discetebouncing} if $\tau^+(n)$ is a good peeling time then there is a positive probability that $\gamma$ touches $ \mathbb{R}_{+}$ between the peeling steps $ \tau^+(n)$ and $ \tau^+(n)+2$. We claim that in fact, the SLE curve will hit $  \mathbb{R}_{+}$ between the peeling steps $ \tau^{(1)}(n) = \tau^{+} (n)$ and $ \tau^{(2)}(n) = \tau^-(\tau^+(n))$ with a probability tending to $1$ as $ n \to \infty$. 

Indeed, by standard properties of the stable process, the time  $\xi^+(1) \in \mathcal{R}^+$ is not isolated from the right  in $ \mathcal{R}^+$. Using \eqref{eq:jumpmin} and properties of the Skorokhod topology, it follows from Theorem*  \ref{thm:exploSLE} that for any $p \geq 0$  we have 
 \begin{eqnarray*} \frac{\tau^{+,(p)}(n)- \tau^+(n)}{n} & \xrightarrow[n\to\infty]{(P)} & 0,  \end{eqnarray*}
 where $\tau^{+,(p)}$ is the $p$-fold composition of $\tau^+$. Since $n^{-1}\tau^{(2)}(n)$ converges in distribution towards  $\xi^{(2)}(1) > \xi^+( 1)$,  we have $$\#\big\{  \tau^+(n) \leq i \leq \tau^-(\tau^+(n)) :  \mathcal{H}^+(i+1/2)= \underline{ \mathcal{H}}^+(i)\big\}  \xrightarrow[n\to\infty]{(P)} \infty.$$  We then claim that the last display remains true if we only restrict to good peeling times. A formal proof of this fact is tedious and we shall not enter these details since we anyway rely on $(*)$. Applying successively Lemma \ref{lem:discetebouncing} to these times,  we deduce that with high probability the SLE curve touches $ \mathbb{R}_{+}$ after the step ${\tau^+(n)}$ but before step $ \tau^{(2)}(n)$. An easy extension of the above argument then yields 
 \begin{eqnarray*} P\Big(\tilde{\theta}^{(k)}(n) \leq (1+ \varepsilon)\tau^{(k)}(n)\Big) & \xrightarrow[n\to\infty]{} & 1, \end{eqnarray*} for any $ k \in \{ 1,2,\ldots\}$ and any $ \varepsilon >0$.

 \endproof

The next two sections are devoted to two computations which investigate the behavior of $ \theta^{(n)}(1)$ and $\xi^{(n)}(1)$ as $n \to \infty$. These are technical propositions and their proofs can be skipped at first reading. This piece of information, combined with Proposition* \ref{prop:limitbouncing} is the heart of the proof of Theorem* \ref{thm:main}. 
 Since we will heavily deal with large deviations estimates we introduce a special notation for it. 

 \paragraph{A notation for large deviations. }Let $I = \mathbb{Z}_{+}, \mathbb{R}_{+}$ or $(0,1)$ and $\omega \in \{0, \infty\}$. If a real stochastic process $(X_{i})_{i \in I}$ indexed by $I$ satisfies a weak law of large numbers: $$ \lim_{i \to \omega} \frac{X_{i}}{ f(i)} = K,$$ in probability for some function $f$ such that $|f| \to \infty$ as $ i \to \omega$ (e.g. $f(i) = i$ or $f(i) = \log i$) and some constant $K \in \mathbb{R}$, we will say that large deviations hold if for every $\eta >0$ there exist $c_{1},c_{2}>0$ (which depend on $\eta$) such that for all $i \in I$ sufficiently close to $\omega$ we have
$$ P \left(\left| \frac{X_{i}}{f(i)} - K \right| >\eta \right) \leq c_{1}e^{-c_{2} |f(i)|},$$ and we write  $$  \displaystyle \limLD{X_{i}}{f(i)}{i \to \omega} = K.$$ 

 Let us give a few examples. The most basic one is to consider a sequence $\zeta_{1}, \ldots , \zeta_{n}$ of i.i.d.\,\,random variables such that $ E[\exp( \lambda |\zeta|)]< \infty$ for some $\lambda >0$.  Then by classical results on large deviations, their partial sums $S_{n} = \zeta_{1} + \cdots + \zeta_{n}$ satisfy   \begin{eqnarray} \label{eq:LDtrivial} \limLD{S_{n}}{n}{n \to \infty} &=& E[\zeta]. \end{eqnarray} 
 Various other examples will arise in this work and are based on scale invariance. E.g., consider the $ \frac{3}{2}$-stable process $S^+$ and its infimum process $ \underline{S}^+$. For any $t >0$, by the scaling property we have $\underline{S}^+_{t}= t^{2/3} \underline{S}^+_{1}$. Also, by standard properties, the law of $\underline{S}^+_{1}$ has a polynomial tail in $-\infty$ and a bounded density around $0$, thus we have  $ P( |\log(-\underline{S}^+_{1})| > x) \leq e^{-cx}$  for some $c>0$ as $x \to \infty$. For every $\eta >0$ and $t \geq 1$ we have
  \begin{eqnarray} 
  P \left( \left| \frac{\log - \underline{S}_{t}^+}{\log t} - \frac{2}{3} \right| > \eta \right) &=&   P \left( \left| \log - \underline{S}_{1}^+\right| > \eta \log t \right) \nonumber \\
  & \leq & c_{1} \exp(-c_{2} \log t)\nonumber \\ & \Rightarrow & \label{eq:LDstable} \underset{t \to \infty}{ \mathsf{limLD}}\,\left( \log |\underline{S}^+_{t}|, \log t\right) = \frac{2}{3}.  \end{eqnarray}
 The last display also holds if we replace $ t \to \infty$ by $t \to 0$. 
 Another useful example comes from the SLE$_{6}$ curve $(\gamma_{t})$ on $ \mathbb{H}$. For any $t>0$ we consider the random variable $ \underline{\gamma}^+_{t} = \sup \left\{\gamma_{[0,t]} \cap \mathbb{R}_{+}\right\}$. By the scaling property of the SLE process, for any $t >0$  we have $$ \underline{\gamma}^+_{ t} =t^{1/2} \underline{\gamma}^+_{1}$$ in distribution. Furthermore, by standard properties  \cite[Chapter 6]{Law05} there exists $c>0$ such that we have $P( \underline{\gamma}^+_{1} \leq \varepsilon) \leq \varepsilon^{c}$ as $ \varepsilon \to 0$ as well as $P( \underline{\gamma}^+_{1} \geq x) \leq x^{-c}$ as $x \to \infty$. Using the same proof as above we deduce that 
  \begin{eqnarray} \label{eq:SLEld} \underset{t \to \infty}{\mathsf{limLD}}\ \left(\log \underline{\gamma}^+_{t}, \log t \right) = \frac{1}{2}.\end{eqnarray}
  Again, the last display holds when $t \to 0$ instead of $t \to \infty$.

%%%%%%%%%%%%%%%%%%%%SUBSECTION
\subsection{A $ \frac{3}{2}$-stable calculation: estimates for $ \xi^{(n)}$}
Recall the definition of $\xi^{(n)}(t)$ from Section \ref{sec:discretebouncing}. In order to lighten notation, in this section we put $\xi^{(n)}:= \xi^{(n)}(1)$ for every $n \geq 0$. 
\label{sec:comStable}
 
 \begin{proposition} \label{prop:commute} We have $ \displaystyle \underset{n \to \infty}{ \mathsf{limLD}}\,\left(\log \xi^{(n)},n\right)  = \frac{\pi}{\sqrt{3}}$.
 \end{proposition}

 Before starting the proof, let us recall some useful facts about the $\frac{3}{2}$-stable process.  We refer to \cite{Ber96,Ber99} for the derivations of these classical identities. Let $S$ be a standard $\frac{3}{2}$-stable Lévy process with no positive jumps and let $ \underline{S}$ be its running infimum process. The reflected process $S- \underline{S}$ admits a local time at $0$ denoted by $(L_{t})_{t \geq 0}$. Its right-continuous inverse $L^{-1}$ is a $ \frac{1}{3}$-stable subordinator (\cite[Chap. VIII, Lemma 1]{Ber96}) and thus follows the generalized arcsine law (\cite[Chap. III, Theorem 6]{Ber96}): For every $x >0$   \begin{eqnarray} 
 \quad x^{-1}\sup \{ t \leq x : S_{t} = \underline{S}_{t}\} \quad \overset{(d)}{=}\quad \frac{ \sqrt{3}}{2\pi} s^{-2/3}(1-s)^{-1/3} ds 1_{(0,1)}(s).  \label{eq:arcsine} \end{eqnarray}
Recall that a random closed set $ \mathcal{S} \subset \mathbb{R}_{+}$ such that almost surely $ \mathcal{S}$ is not bounded, has no isolated point and such that $0
\in  \mathcal{S}$ is a \emph{regenerative} set if for any $t \geq 0$, conditionally on $Z_{t} =\min [t,\infty)\cap  \mathcal{S}$, the set $(  \mathcal{S}\cap [Z_{t},\infty))-Z_{t}$ is independent of $(\mathcal{S} \cap [0,Z_{t}])$ and is distributed as $ \mathcal{S}$. Any regenerative set can be seen as the range of a subordinator unique up to multiplicative constant, see \cite{Ber99}. A regenerative set is thus characterized by a drift parameter $d\geq 0$ and a positive Lévy measure $\pi$ (called the regenerative measure) unique up to multiplication by the same constant.  

In our case, the random closed set $ \mathcal{R} = \{ t \geq 0 : S_{t} = \underline{S}_{t}\}$ is a regenerative set (it corresponds to the range of the subordinator $L^{-1}$) with no drift and regenerative measure 
 \begin{eqnarray} \label{eq:regemes} x^{-4/3} \mathbf{1}_{x>0} dx. \end{eqnarray} 
For every $t >0$, almost surely $t \notin \mathcal{R}$ and $ \mathcal{R}$ has Hausdorff dimension $1/3$.  Also recall from \cite[Chap. 5]{Ber99} that the intersection of two independent copies of $ \mathcal{R}$ is almost surely reduced to $\{0\}$.

\begin{figure}[!ht]
 \begin{center}
 \includegraphics[width=14cm]{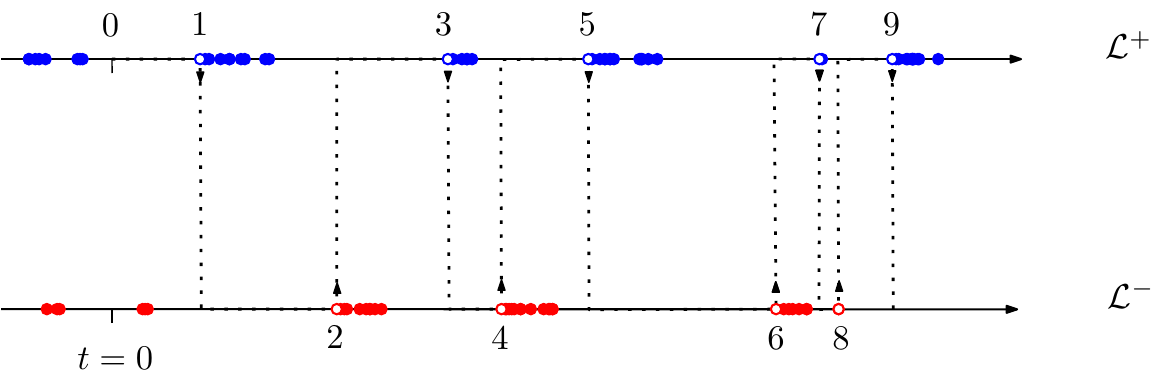}
 \caption{ \label{fig:commutelog} The sets $ \mathcal{L}^+$ and $\mathcal{L}^-$ and the steps of the chain $ (\log \xi^{(n)})_{n \geq 0}$.}
 \end{center}
 \end{figure}

 \proof[Proof of Proposition \ref{prop:commute}] Due to the logarithm in the statement of Proposition \ref{prop:commute} it is more convenient to deal with the logarithm of $ \mathcal{R}^+$ and $ \mathcal{R}^-$: we set $ \mathcal{L}^+= \log( \mathcal{R}^+ \backslash\{0\})$ and 
$ \mathcal{L}^-= \log( \mathcal{R}^-\backslash\{0\})$. Clearly we have $ \mathcal{L}^+ \cap \mathcal{L}^-=  \varnothing$ and $ \xi^{(n)}$ is measurable with respect to $ \mathcal{L}^+$ and $ \mathcal{L}^-$, see Fig.\,\ref{fig:commutelog}. It turns out that $ \mathcal{L}^+$ and $  \mathcal{L}^-$ are again regenerative sets, but not started at $0$: Let $ \mathsf{L}$ be a random set having the law of $ \mathcal{L}^{+}$ or $ \mathcal{L}^-$  translated at its first positive value
$$ \mathsf{L} \overset{(d)}{=} \big( \mathcal{L}^+ - \inf \mathcal{L}^+ \cap [0,\infty)\big) \cap [0, \infty).$$

 \begin{lemma} \label{lemma:log}The random set $ \mathsf{L}$ is a regenerative set with no drift and regenerative measure  $$ \nu( dx) = \frac{e^x}{(e^x-1)^{4/3}} dx.$$
 \end{lemma}
 \proof This comes from a straightforward calculation: For every $x_0 \in \mathbb{R}_+^*$ the push-forward of the measure $ \frac{ dx}{x^{4/3}}$ on $ \mathbb{R}_{+}$ given in \eqref{eq:regemes} by the map $ u \mapsto \log(u+x_0)-\log(x_0)$ is a multiple (depending of $x_0$) of the measure $\nu( dx)$. \endproof 
The main observation is the following.

\begin{lemma} \label{lemma:MC}
The process $ X_{n} := \log \xi^{(n+1)} - \log \xi^{(n)}$ for $n \geq 2$ is a Markov chain with transition kernel    \begin{eqnarray*} p(x, dy) &=& \frac{ \sqrt{3}}{2\pi} \left(\frac{e^x-1}{e^{x}(e^y-1)}\right)^{1/3}\frac{dy}{1-e^{-x-y}}.  \end{eqnarray*}
\end{lemma}
\proof For $x >0$ we denote by $G_{x} = \sup\{ s \leq x : s \in  \mathsf{L}\}$ and $D_{x} = \inf \{ s \geq x : s \in \mathsf{L}\}$. Since almost surely $x \notin \mathsf{L}$, we have  $ G_{x}< x< D_{x}$.  We denote by $ p_{x}( dy)$ the law of $D_{x}-x$ and will show that $(X_{i})_{i \geq 2}$ is a Markov chain with transition kernel $p(x, dy) = p_{x}( dy)$. 

 Let $i \geq 2$ be odd (say). The point $\log \xi^{(i)}$ thus belongs to $ \mathcal{L}^+$. Note that $X_{1},\ldots,X_{i-1}$ are measurable with respect to 
$$ \mathfrak{F}_{i} := \sigma(\mathcal{L}^+ \cap [0, \log \xi^{(i)}],\mathcal{L}^- \cap [0, \log \xi^{(i-1)}]).$$ We thus condition on $ \mathfrak{F}_{i}$ and look for the next point larger than or equal to $\log \xi^{(i)}$ belonging to $\mathcal{L}^-$.  By the regenerative property of $ \mathcal{L}^-$, the conditional distribution of $ \mathcal{L}^- \cap [ \log \xi^{(i-1)}, \infty)$ is that of $ \mathsf{L} +\log \xi^{(i-1)}$ (here we use that $i \geq 2$). The conditional law of $X_{i} = \log \xi^{(i+1)}- \log \xi^{(i)}$ is that of the law of $D_{X_{i-1}}-X_{i-1}$. Consequently, conditionally on $X_{1},\ldots,X_{i-1}$, the variable $X_{i}$ is distributed as $p_{X_{i-1}}( dy)$ as desired. 
 
 Let us now compute the distribution $p_{x}( dy)$ for $x > 0$. This is a pretty straightforward calculation but we provide the details for the reader's convenience. We first compute the distribution of $G_{x}$. For this we use the arcsine law \eqref{eq:arcsine} on the original regenerative set $ \mathcal{R}$. Indeed if $ \mathsf{R}$ is a version of $ \mathcal{R}$ started at time $1$ then $G_{x}$ has the same distribution as 
  \begin{eqnarray*} G_{x} & \overset{(d)}{=} & \log \big( \sup\{ s \leq e^x : s \in \mathsf{R}\}\big).  \end{eqnarray*}
  The law of the random variable inside the logarithm on the right-hand side minus $1$ and divided by $e^x-1$ is the arcsine law \eqref{eq:arcsine} of parameter $1/3$. In other words, for any positive measurable function $f$ we have 
  $$ E[f(G_{x})] = \int_{0}^1 ds\ \frac{ \sqrt{3}}{2 \pi} s^{-2/3}(1-s)^{-1/3} f\Big( \log \big( (e^x-1)s+1 \big)\Big).$$ Performing the change of variable $u =\log ( (e^x-1)s+1)$, the law of $G_{x}$ is given by   \begin{eqnarray} \label{lawGt} G_{x}& \overset{(d)}{=} & \frac {\sqrt{3}}{2 \pi} \frac{  e^u \mathbf{1}_{0 < u< x}}{(e^u-1)^{2/3}(e^x-e^u)^{1/3}}du .  \end{eqnarray}
Finally, conditionally on $G_{x}$, by the regenerative property of $ \mathsf{L}$,  the law of $D_{x}-x$ is given by $\nu( dy \mid y > x-G_{x})$. Using Lemma  \ref{lemma:log} (and the easy identity $ \nu [x, \infty) = 3/(e^x-1)^{1/3}$ for $x >0$) we get that for any positive measurable $f$
 \begin{eqnarray*}
 && \hspace{-1.5cm}E[f(D_{x}-x)] \\
 &=& \int_{0}^x du \ \frac {\sqrt{3}}{2 \pi} \frac{  e^u}{(e^u-1)^{2/3}(e^x-e^u)^{1/3}} \cdot \frac{(e^{x-u}-1)^{1/3}}{3} \int_{x-u}^\infty da\ \frac{e^a}{(e^a-1)^{4/3}}f\big(a-(x-u)\big) \\
 & \underset{\begin{subarray}{c}v = x-u\\ y= a-v \end{subarray}}{=} &  \frac {1}{2 \pi \sqrt{3}}\int_{0}^x dv \  \frac{  e^{x-v}(e^{v}-1)^{1/3}}{(e^{x-v}-1)^{2/3}(e^x-e^{x-v})^{1/3}} \int_{0}^\infty dy\ \frac{e^{y+v}}{(e^{y+v}-1)^{4/3}}f(y)   \\
 &=&\frac {1}{2 \pi \sqrt{3}} \int_{0}^\infty dy\  f(y) \int_{0}^x  dv\ \frac{e^{y+v}}{(e^{y+v}-1)^{4/3}} \frac{  e^{-2v/3}}{(e^{-v}-e^{-x})^{2/3}}.\end{eqnarray*} The last integral has been computed using Mathematica$^\copyright$, however it is easy (but tedious) to check \emph{a posteriori} that it is equal to the formula provided in the statement of the lemma.
 \endproof	 
Using the exact form of the probability transitions of the chain $(X_{i})$ it is easy to see that this chain is aperiodic, recurrent and ergodic. Furthermore, its  unique invariant and reversible probability measure is given by   \begin{eqnarray*} \varpi( dx) &=&  \frac{2^{2/3}\sqrt{\pi}}{ \Gamma(1/3) \Gamma(1/6)}\frac{e^x}{ \big( e^x (e^x-1)\big)^{2/3}}.  \end{eqnarray*} An application of the ergodic theorem  implies that   $n^{-1}\log \xi^{(n)} = n^{-1}(X_{0}+X_{1} + \cdots +X_{n-1})$ converges almost surely and in $ \mathbb{L}^1$ towards\footnote{Here and later, unexplained integrations have been realized using Mathematica$^\copyright$} $$\int_{ \mathbb{R}_{+}} \varpi( dx) x = \frac{\pi}{ \sqrt{3}},$$ which is the constant appearing in the statement of the Proposition \ref{prop:commute}.

We now set up large deviations estimates. Recall that a set $C \subset \mathbb{R}_{+}$ is small \cite[p 102]{MT09}  if there exists a probability measure $ \nu$ and $ \varepsilon>0$ such that for some $k \geq 1$ the $k$-steps transition kernel satisfies $$ p^k(x,A) \geq \varepsilon \nu(A), \quad \forall x \in C, \quad \forall A \mbox{\ Borel}.$$
The chain is uniformly ergodic  (see \cite[Theorem 16.0.2]{MT09}) if the full space is small. Unfortunately for us, it is easy to see that $p(x, dy)$ is concentrated around $0$ when $x$ is close to $0$. Hence the chain is not uniformly ergodic and some care is needed. It is however easy to see from the exact form of the transition kernels that any set $[b, \infty)$ with $b >0$ is a small set. 
We now establish that the chain $(X_{i})$ is $V$-geometrically ergodic (see \cite[Theorem 16.0.1]{MT09}) with the function $$V : x \in \mathbb{R}_{+} \mapsto (x\vee x^{-1/4}) \in \mathbb{R}_{+}.$$ This will allow us to apply the powerful machinery developed in  \cite[Chapter 16]{MT09}. For this, we compute the variation of $V$ after applying a one step transition of  the chain:
 \begin{align*} 
  pV(x) := \int p(x,dy) V(y) &= \frac{ \sqrt{3}}{2\pi} \left(\frac{e^x-1}{e^{x}}\right)^{1/3}\int_{0}^\infty  \frac{dy\ (y \vee y^{-1/4})}{(e^y-1)^{1/3}1-e^{-x-y}}. 
  \end{align*}
It is easy to see that $pV(x) \leq K$ for some constant $K>0$ uniformly in $x \geq 1$. On the other hand, when $x \to 0$ we have 
$$pV(x) \sim \frac{\sqrt{3}}{2\pi} x^{1/3} \int_{0}^\infty \frac{ dy \ y^{-1/4}}{y^{1/3}(x+y)} = x^{-1/4} \cdot \frac{\sqrt{3}}{2\pi} \int_{0}^\infty \frac{ dz}{z^{7/12}(1+z)}$$
where we have performed the change of variable $y=xz$. The right-hand side can be computed exactly and is equal to $V(x) \cdot \frac{ \sqrt{6}}{1+ \sqrt{3}}$ for $x <1$. Since $\frac{ \sqrt{6}}{1+ \sqrt{3}} <1$ we deduce that the condition $(V4)$ of \cite[p 376]{MT09} is indeed satisfied and thus the chain is $V$-geometrically ergodic. 

We first establish upper large deviations for the partial sums of the $X_{i}$. Fix $\eta >0$ and find $a\geq 0$  such that 
  \begin{eqnarray}\label{eq:eta}\int_{a}^\infty \frac{ dy\ y}{(1- e^{-y})(e^y-1)^{1/3}} \leq  \frac{\eta}{10}. \end{eqnarray} 
We have 
 \begin{align} P \left( n^{-1}\sum_{i=0}^{n-1} X_{i} - \frac{\pi}{ \sqrt{3}} > \eta \right)  & \leq P \left( n^{-1}\sum_{i=0}^{n-1} X_{i} \mathbf{1}_{X_{i} \geq a}  > \eta/2 \right) \label{eq:trois} + P \left( n^{-1}\sum_{i=0}^{n-1} X_{i}\mathbf{1}_{X_{i} < a} - \frac{\pi}{ \sqrt{3}} > \eta/2 \right). \end{align}
 Note that $F:x \mapsto x \mathbf{1}_{x<a}$ is a bounded function, so we can apply the results of \cite{KM03} and get that $ \limLD{\sum_{i=0}^n F(X_{i})}{n}{n \to \infty} = E_{\pi}[F(X)]$. In particular since $E_{\pi}[F(X)] \leq E_{\pi}[X]= \pi/ \sqrt{3}$ we deduce that for some $c_{1},c_{2}>0$ we have 
 $$ P \left( n^{-1}\sum_{i=0}^{n-1} X_{i}\mathbf{1}_{X_{i} < a} > \eta/2  + \frac{\pi}{ \sqrt{3}} \right) \leq c_{1} e^{-c_{2}n}.$$
To control the other term of the right-hand side of \eqref{eq:trois}, we remark that the probability transitions of the chain $(X)$ are bounded from above by
$$ p(x,dy) \leq \frac{ dy}{(1- e^{-y})(e^y-1)^{1/3}}.$$
And so 
$$ P \left( n^{-1}\sum_{i=0}^{n-1} X_{i} \mathbf{1}_{X_{i} \geq a}  > \eta/2 \right) \leq P \left( n^{-1}\sum_{i=0}^{n-1} Z_{i}  > \eta/2 \right),$$ where $(Z_{i})$ are i.i.d. random variables of law given by 
$$ \mathbf{1}_{y>a}\frac{ dy}{(1- e^{-y})(e^y-1)^{1/3}} + \left(1- \int_{a}^\infty \frac{ dx}{(1- e^{-x})(e^x-1)^{1/3}}\right) \delta_{0}.$$
Since $Z_{i}$ has mean less than $ \eta/10$ by \eqref{eq:eta} and has exponential moments, large deviations estimates \eqref{eq:LDtrivial} show that the last term is bounded by $c_{1} e^{-c_{2}n}$ for some $c_{1},c_{2}>0$. This completes the upper large deviations for the partial sums of the chain $X$,  the lower large deviations are similar and left to the reader. \qedhere  \medskip 

As a corollary of the last proposition, we study  the number of alternative minimal records of two independent stable processes $(S^+,S^-)$ between the times when $\underline{S}^+$ is between two fixed values. More precisely, for any $ x>0 $ set $\vartheta_{ x} = \inf\{ {t\geq 0} : \underline{S}^+_{t} \leq - x \}$ and for $0 <x<y$ put
$$  \mathsf{ComStable}(x,y) = \inf \Big\{ k \in \{ 1,3,5,\ldots\} : \xi^{(k)}(  \vartheta_{x}) \geq \vartheta_{y} \Big\}.$$ 
Remark that by scale invariance of the stable processes we have $ \mathsf{ComStable}(x,y) = \mathsf{ComStable}(1,y/x)$ in distribution. 
\begin{corollary} \label{cor:commuteStable} We have $\displaystyle \limLD{\mathsf{ComStable}(1,x)}{\log x}{ x \to \infty} = { \frac{ 3 \sqrt{3}}{  2\pi}}.$
\end{corollary}
\proof  By monotonicity of $ t \mapsto \xi^{(k)}(t)$ and of $ k \mapsto \xi^{(k)}(t)$, note that if $[b \log x]$ is odd and if we have simultaneously $ \vartheta_{x} \leq x^{a}$, $ \vartheta_{1} \geq x^{- \varepsilon}$ and  $ \xi^{([b \log x])} (x^{- \varepsilon}) \geq x^{a}$ then we have $\mathsf{ComStable}(1,x) \leq [b \log x]$. Taking $ \frac{\pi}{\sqrt{3}} b > (a + \varepsilon)> a > 3/2$ and using the last proposition together with \eqref{eq:LDstable} we get that for large $x$  so that $[b\log x]$ is odd
 \begin{eqnarray*} && \hspace{-3cm}P ( \mathsf{ComStable}(1,x) >  [b \log x])\\ &\leq& P( \xi^{([b \log x])} (x^{- \varepsilon})\leq x^a) + P( \vartheta_{x} \geq x^a)  + P( \vartheta_{1} \leq x^{- \varepsilon})\\
 &=& P\left( x^{ \varepsilon} \cdot \xi^{([b \log x])}(x^{-  \varepsilon}) \leq x^{a+ \varepsilon}\right) + P( \underline{S}^+_{x ^a}\geq -x ) + P( \underline{S}^+_{x^{- \varepsilon}} \leq -1) \\
 & \underset{  \mathrm{scaling}}{=} & P \left( \frac{\log \xi^{([b \log x])}}{\log x} \leq (a+ \varepsilon)\right) + P(\underline{S}^+_{1} \geq -x^{1-2a/3}) + P(\underline{S}^+_{1} \leq -x^{ 2\varepsilon/3}) \\
 & \underset{\begin{subarray}{c} \mathrm{Prop.} \ref{prop:commute} \\ \mathrm{ and }\ \eqref{eq:LDstable} \end{subarray}}{ \leq }&  c_{1} e^{-c_{2} \log x} \end{eqnarray*}
 for some $c_{1},c_{2}>0$ (depending on $a,b$ and $ \varepsilon$). This thus holds for any $b > \frac{3 \sqrt{3}}{2 \pi}$. The other inequality is similar. This  proves the corollary. \qedhere

%%%%%%%%%%%%%%%%%%%%SUBSECTION
\subsection{An SLE$_{6}$ calculation: estimates for $ \theta^{(n)}$} 
 In order to lighten notation, in this section we put $\theta^{(n)}:= \theta^{(n)}(1)$.

\begin{proposition} \label{prop:commuteSLE} We have $ \displaystyle \limLD{\log \theta^{(n)}}{n}{n \to \infty}  = \frac{4\pi}{\sqrt{3}}$.
 \end{proposition}
 \proof  Recall the notation of Section \ref{sec:defSLE}.  We start by a few classical facts on the SLE processes, see e.g. \cite[Section 8.3]{LSW03} for details. The image of the boundary of the hull $ \mathbb{H} \backslash H_{t}$ inside $ \mathbb{H}$ is sent by the uniformization mapping $g_{t}$ to a segment $[L_{t},R_{t}]$ and $g_{t}(\gamma_{t})$, denoted by $U_{t}$, lies inside $[L_{t},R_{t}]$. In particular,  the times when $\gamma$ touches $ \mathbb{R}_{+}$ correspond to the times when $ U_{t} = R_{t}$ and similarly for the left part. That is $ \theta^{(0)}=1$ and  $\theta^{(i+1)}= \inf \{ t \geq \theta^{(i)} : U_{t}=R_{t}\}$ for $i$ even and $\theta^{(i+1)}= \inf \{ t \geq \theta^{(i)} : U_{t}=L_{t}\}$ for $i$ odd. We now derive the equations driving these processes, we refer to \cite[Section 8.3]{LSW03} for more details. The Loewner equation \eqref{def:SLE6} tells us that $$dU_{t} = \sqrt{6} dB_{t},$$  where $B$ is a standard Brownian motion. Also an easy calculation using \eqref{def:SLE6} shows that as long as $ U_{t} \ne R_{t}$ and $ U_{t} \ne L_{t}$ we have 
 \begin{eqnarray} \label{eq:LtRt}
 \displaystyle dL_{t} = \displaystyle\frac{2 dt}{L_{t}-U_{t}},\quad \mbox{ and }\quad 
 \displaystyle dR_{t} = \displaystyle \frac{2 dt}{R_{t}-U_{t}}.  \end{eqnarray} However, since the parameter of the SLE is $ \kappa=6$, we will have infinitely many times at which $L_{t}= U_{t}$ or $R_{t}=U_{t}$ and the meaning of the last display is not clear anymore. One way to cope with this to first define simultaneously the processes $G_{t} = (U_{t}-L_{t})^2$ and $D_{t}=(R_{t}-U_{t})^2$ as the solutions of
 $$ 
 dG_{t} = 2 \sqrt{G_{t}} \, dU_{t} + 10 \,dt \quad \mbox{and} \quad  dD_{t} = - 2 \sqrt{D_{t}}\,dU_{t} + 10 \,dt,$$ starting from $0$ (with the same Brownian motion). Consequently, both $G_{t}$ and $D_{t}$ are distributed as $ 1/{6}$ times a squared Bessel process of dimension $5/3$ and are defined for all $t \geq 0$, see \cite[Chap. XI]{RY99}. From the triplet $(G_{t},U_{t},D_{t})$ we can then construct $(L_{t},U_{t},R_{t})=(U_{t}- \sqrt{G_{t}},U_{t}, \sqrt{D_{t}}-U_{t})$ for all times $t \geq 0$, see \cite[Chapter 8.3]{LSW03}. If we put $$X_{t} = \frac{U_{t}-L_{t}}{R_{t}-L_{t}}$$ and $ \Delta_{t} =  R_{t}-L_{t}$ applying Ito's formula we get
 \begin{eqnarray*} dX_{t} &=& \frac{ dU_{t}}{ \Delta_{t}} + \frac{2 dt}{ \Delta_{t}^2} \Big( \frac{1}{X_t}- \frac{1}{1-X_{t}}\Big),  \end{eqnarray*}
 which can be defined for all $t \geq 0$ using the above device. Performing the following time-change  $$r(t) = \int_{1}^{t} \frac{ ds}{ \Delta_{s}^2}, \qquad Z_{r(t)}= X_{t} $$ we obtain that $Z$ satisfies 
$$ \left\{ \begin{array}{ccl} Z_{0} & = &X_{1}\\ dZ_{t}&=& dU_{t} + 2 dt \Big( \displaystyle \frac{1}{Z_{t}}- \frac{1}{1-Z_{t}} \Big).  \end{array} \right.$$ 
After these transformations the alternative hitting times  of $1$ and $0$ by the process $Z$ are given by $r(\theta^{(i)})$ for $i \geq 1$. We now have two tasks. Firstly, understand the number of commutings between $0$ and $1$ for the process $Z$ as time goes to infinity, and secondly understand the asymptotics of the time change $r(t)$ in order to translate these results back to the $\theta^{(i)}$. 

\textsc{Commutings of $Z$.} The process $Z$ is strong Markov and symmetric with respect to $1/2$.  By looking at the SDE governing $Z$ we see that it evolves like a Bessel of dimension $5/3$ around $0$ and symmetrically around $1$. In particular, starting from $0$ the process $Z$  will eventually hit $1$ in finite time a.s.\,\,and vice versa. For $x \in [0,1]$, under $E_{x}$ the process $Z$ starts from $x$. We let $ \mathfrak{t}_{0}$ and $ \mathfrak{t}_{1}$ be the hitting times of $0$ and $1$ respectively by the process $Z$. We now state a technical lemma:

\begin{lemma} \label{lem:estim} For some $\lambda >0$ we have  \begin{eqnarray*}  E_{0} \left[ \exp\left(  \lambda \int_{0}^{\mathfrak{t}_{1}} \frac{du}{Z_{u}(1-Z_{u})} \right) \right] < \infty.  \end{eqnarray*} 
In particular $E_{0}[\exp( \lambda   \mathfrak{t}_{1})] < \infty$. Furthermore we have 
$$E_{0}[ \mathfrak{t}_{1}] = E_{1}[ \mathfrak{t}_{0}] = \frac{\pi}{7 \sqrt{3}}.$$
 \end{lemma}

\proof[Proof of Lemma \ref{lem:estim}] By symmetry in space and time of the process $Z$, to prove the first assertion of the lemma, it is sufficient to prove that for some $\lambda >0$ we have 
 \begin{eqnarray} \label{eq:g1}  E_{1/2} \left[ \exp\left(  \lambda \int_{0}^{\mathfrak{t}_{0} \wedge \mathfrak{t}_{1}} \frac{du}{Z_{u}(1-Z_{u})} \right) \right] < \infty.  \end{eqnarray} 
 For this we introduce the scale function $\phi$ of the process $Z$ which is defined for $x \in [0,1]$ by  
 $$ \phi(x) = \int_{0}^x \frac{du}{\big(u(1-u)\big)^{2/3}}.$$
 In particular we have $\Lambda := \phi(1) =  \frac{\Gamma(1/6) \Gamma(1/3)}{2^{2/3} \sqrt{\pi}}$ and $\phi$ satisfies $2 \phi'(x)( \frac{1}{x}- \frac{1}{1-x}) + 3 \phi''(x)=0$. Applying Ito's formula, it comes as no surprise that $Y_{t} = \phi(Z_{t \wedge \mathfrak{t}_{0} \wedge \mathfrak{t}_{1}})$ is a local martingale under $E_{1/2}$. Since the later is bounded it is even a true martingale. By the Dubins-Schwarz theorem, $Y$ is a time change of a Brownian motion. Specifically, we can write $ Y_{t} = \beta_{<Y>t}$ where $\beta$ is a standard Brownian motion started from $ \Lambda/2$. Stochastic calculus shows that $ d\hspace{-1mm}<\hspace{-1mm}Y\hspace{-1mm}>_{u} = 6 ( \phi' \circ \phi^{-1})^2(\beta_{<Y>_{u}})du$, consequently after the change of variable $v = <\hspace{-1mm}Y\hspace{-1mm}>_{u}$ we have
  \begin{eqnarray*}\int_{0}^{ \mathfrak{t}_{0} \wedge \mathfrak{t}_{1}} \frac{du}{Z_{u}(1-Z_{u})} &=& \int_{0}^{ \mathfrak{t}_{0} \wedge \mathfrak{t}_{1}} \frac{du}{ \phi^{-1}( \beta_{<Y>_{u}})(1- \phi^{-1}( \beta_{<Y>_{u}})} \\ &=& \int_{0}^{\tau_{0}\wedge \tau_{\Lambda}} \frac{dv}{ 6 (\phi' \circ \phi^{-1})^2 \cdot \phi^{-1} (1- \phi^{-1}) (\beta_{v})}\\
  &=& \int_{0}^{\tau_{0}\wedge \tau_{\Lambda}} dv \ \psi \circ \phi^{-1}( \beta_{v})   \end{eqnarray*}
 where $\psi : x \mapsto 6(x(1-x))^{1/3}$ and $\beta$ is a Brownian motion started from $ \Lambda/2$ and stopped at $\tau_{0}\wedge \tau_{\Lambda}$, the first hitting time of $0$ or of $\Lambda$ by $\beta$. Since the function $\psi$ is bounded by $6$ over $[0,1]$ we deduce that  
 $$ \int_{0}^{ \mathfrak{t}_{0} \wedge \mathfrak{t}_{1}} \frac{du}{Z_{u}(1-Z_{u})} \leq 6 \cdot \tau_{0} \wedge \tau_{ \Lambda}. $$ It is classical that $\tau_{0}\wedge \tau_{ \Lambda}$ has some exponential moment and so \eqref{eq:g1} follows.
 
 We easily deduce from the first point that under $E_{0}$ the variable $ \mathfrak{t}_{1}$ possesses some exponential moments and is in particular integrable. To compute its expectation consider now the function $f : [0,1] \to \mathbb{R}_{}$ which is $ \mathcal{C}^2$ over $[0,1)$, with  $f(0)=0$, $f'(0)=0$ and which satisfies the differential equation $$ 2 f'(x) \Big( \frac{1}{x}- \frac{1}{1-x}\Big) + 3 f''(x) =1.$$ Such a function exists and an can be expressed using hypergeometric functions\footnote{A computation with Mathematica gives $f(x) = 1/14 (-x + x^2 + x \mathrm{HypergeometricPFQ}[\{1, 1, 4/3\}, \{5/3, 2\}, x])$}. This function is positive, continuous over $[0,1]$ and $f(1)-f(0)= \frac{\pi}{7 \sqrt{3}}$. Another application of Ito's formula shows that $$( f(Z_{t})-t)_{t < \mathfrak{t}_{1}}$$ is a local martingale. Since $f$ is bounded it is even a true martingale. Applying the optional sampling theorem we deduce that  $E_{0}[f(Z_{t \wedge \mathfrak{t}_{1}})] = E_{0}[ \mathfrak{t}_{1} \wedge t]$ for every $t \geq 0$. Letting $t \to \infty$ we get by the dominated and monotone convergences theorems that \begin{eqnarray*} E_{0}[ \mathfrak{t}_{1}] = f(1)-f(0) = \frac{ \pi}{7 \sqrt{3}}.   \end{eqnarray*} \endproof

Let us now come back to the proof of Proposition \ref{prop:commuteSLE}. By applying the strong Markov property at the successive and alternate hitting times of $1$ and $0$ by the process $Z$, we deduce that the $n$th interlaced hitting time $r(\theta^{(n)})$ of $\{0,1\}$ by the process $Z$ is given by $ r(\theta^{(1)})+ \mathfrak{t}^{(2)}+ \cdots + \mathfrak{t}^{(n)}$ where $ \mathfrak{t}^{(i)}$ are i.i.d.\,\,copies of $ \mathfrak{t}_{1}$ under $E_{0}$. We deduce from Lemma \ref{lem:estim} and \eqref{eq:LDtrivial} that  \begin{eqnarray} \label{eq:ghjf} \limLD{r(\theta^{(n)})}{n}{n\to \infty} &=& \frac{\pi}{7 \sqrt{3}}.  \end{eqnarray}  

\textsc{Asymptotics of the time-change.} We now prove that
 \begin{eqnarray} \limLD{r(t)}{\log t}{t\to \infty}&=& \frac{1}{28}, \label{eq:tich} \end{eqnarray} which will together with the last display imply the proposition. Indeed, by monotonicity of $t \mapsto r(t)$ and $k \mapsto \theta^{(k)}$, if for some number $c\geq 0$ we have both $ r(\theta^{(n)}) \geq c$ and $r(t) \leq c$ then $\theta^{(n)} \geq t$. Choosing $ \frac{\pi}{7 \sqrt{3}}>a$ and $b> \frac{1}{28}$ and setting $c=b \log t =  a n$ we have 
 \begin{eqnarray*} P\left( \theta^{(n)} < t \right)=P\left( \frac{\log \theta^{(n)}}{n} < \frac{a}{b} \right) &\leq& P\left( r(\theta^{(n)}) \leq  an\right) + P\left( r(t) \geq b \log t \right) \\ & \underset{\begin{subarray}{c} \eqref{eq:ghjf} \ \mathrm{and }\ \eqref{eq:tich} \end{subarray}}{\leq}& c_{1}\exp(-c_{2} n), \end{eqnarray*} for some constants $c_{1},c_{2}>0$. Since $a/b$ can be made arbitrarily close to $ {4\pi}/{ \sqrt{3}}$ this proves one side of the proposition, the other inequality is similar.
 
From the SDE satisfied by $ \Delta_{s}$ we get  that \begin{eqnarray*} d \Delta_{s} &=& \frac{2 ds}{ \Delta_{s}X_{s}(1-X_{s})}\\
d (\log \Delta_{s}) &=& \frac{ds}{ \Delta_{s}^2} \frac{2}{ X_{s}(1-X_{s})}.  \end{eqnarray*}
Integrating over $[1,t]$ and performing the change of variable $u = r(t)$ with $ du = dt/ \Delta_{t}^2$ we get  \begin{eqnarray} \label{eq:integ} \log( \Delta_{t})- \log( \Delta_{1}) &=& \int_{0}^{r(t)} \frac{2 du}{Z_{u}(1-Z_{u})}. \end{eqnarray}
Recall that $\Delta_{t} = R_{t}-L_{t}$ is the sum of two (depend) multiples of Bessel processes of dimension $5/3$. The scaling property of these then imply that $\Delta_{t} = \sqrt{t} \Delta_{1}$ in distribution and easy estimates actually show that 
 \begin{eqnarray} \label{eq:scaling2} \limLD{\log \Delta_{t}}{\log t}{t \to \infty} = \frac{1}{2}.  \end{eqnarray} On the other hand, recall that the invariant measure  of a  diffusion $ d Z_{t} = - \nabla \psi(t) dt + \sqrt{2 \beta^{-1}} dB_{t}$ is proportional to  $ \rho ( dx) \propto \exp(-\beta \psi(x)) dx$. In the case of $Z$, the invariant probability measure is thus \begin{eqnarray*} \rho( dx) &=& \frac{\Gamma({10/3})}{\Gamma({5/3})^2}\big(x(1-x)\big)^{2/3}.  \end{eqnarray*} In particular an application of the ergodic theorem shows that  \begin{eqnarray*}  \lim_{t \to \infty} t^{-1}\int_{0}^t \frac{ du}{ Z_{u}(1-Z_{u})} &=& \int_{0}^1  \frac{\rho ( dx)}{x(1-x)} =  7,  \end{eqnarray*} almost surely and in $L^1$. We can strengthen the last display. Indeed, by decomposing the process $Z$ into independent excursions between $0$ and $1$, and using Lemma \ref{lem:estim} and \eqref{eq:LDtrivial} one deduces that large deviations hold for the last display, that is 
 $$ \limLD{ \int_{0}^t \frac{du}{Z_{u}(1-Z_{u})}}{t}{t \to \infty} = 7.$$ It is now easy to combine the last display with \eqref{eq:scaling2} and \eqref{eq:integ} to complete the proof of \eqref{eq:tich}. \hfill $\square$  \medskip

As a corollary of the last proposition, we study  the number of alternative bouncings on $ \mathbb{R}_{+}$ and $ \mathbb{R}_{-}$ that the curve $\gamma$ is doing between two fixed points. Recall the notation introduced before \eqref{eq:SLEld}. For any $ x>0$, set $\varkappa_{x} = \inf\{ {t\geq 0} : \underline{\gamma}^+_{t} \geq  x \}$ and for $0 <x<y$ put
$$ \mathsf{ComSLE}(x,y) = \inf \Big\{ k \in \{ 1,3,5,\ldots\} : \theta^{(k)}(  \varkappa_{x}) \geq \varkappa_{y} \Big\}.$$ 
Remark that by scale invariance of the stable processes we have $ \mathsf{ComSLE}(x,y) = \mathsf{ComSLE}(1,y/x)$ in distribution. 
\begin{corollary} \label{cor:commuteSLE} We have $\displaystyle \limLD{\mathsf{ComSLE}(1,x)}{\log x}{ x \to \infty} = { \frac{ \sqrt{3}}{  2\pi}}.$
\end{corollary}

\proof The proof is similar to that of Corollary \ref{cor:commuteStable} and follows from the last proposition together with the square-root scaling property of the SLE$_{6}$ process. Let us repeat the argument. By monotonicity of $ t \mapsto \theta^{(k)}(t)$ and of $ k \mapsto \theta^{(k)}(t)$, note that if $[b \log x]$ is odd and if we have simultaneously $ \varkappa_{x} \leq x^{a}$, $ \varkappa_{1} \geq x^{- \varepsilon}$ and  $ \theta^{([b \log x])} (x^{- \varepsilon}) \geq x^{a}$ then we have $\mathsf{ComSLE}(1,x) \leq [b \log x]$. Taking $ \frac{4\pi}{\sqrt{3}} b > (a + \varepsilon)> a > 2$ and using the last proposition together with \eqref{eq:SLEld} we get that for large $x$  so that $[b\log x]$ is odd
 \begin{eqnarray*} &&\hspace{-3cm}P ( \mathsf{ComSLE}(1,x) >  [b \log x])\\  &\leq& P( \theta^{([b \log x])} (x^{- \varepsilon})\leq x^a) + P( \varkappa_{x} \geq x^a)  + P( \varkappa_{1} \leq x^{- \varepsilon})\\
 &=& P\left( x^{ \varepsilon} \cdot \theta^{([b \log x])}(x^{-  \varepsilon}) \leq x^{a+ \varepsilon}\right) + P( \underline{\gamma}^+_{x ^a}\leq x ) + P( \underline{\gamma}^+_{x^{- \varepsilon}} \geq 1) \\
 & \underset{  \mathrm{scaling}}{=} & P \left( \frac{\log \theta^{([b \log x])}}{\log x} \leq (a+ \varepsilon)\right) + P(\underline{\gamma}^+_{1} \leq x^{1-a/2}) + P(\underline{\gamma}^+_{1} \geq x^{ \varepsilon/2}) \\
 & \underset{\begin{subarray}{c} \mathrm{Prop.} \ref{prop:commuteSLE} \\ \mathrm{ and }\ \eqref{eq:SLEld} \end{subarray}}{ \leq }&  c_{1} e^{-c_{2} \log x}
\end{eqnarray*}
 for some $c_{1},c_{2}>0$ (depending on $a,b,\varepsilon$). This thus holds for any $b > \frac{\sqrt{3}}{2 \pi}$. The other inequality is similar.  \qedhere

\begin{remark} These commuting estimates for the SLE$_{6}$ are closely related to the work of Hongler and Smirnov \cite{HS11}. Indeed, these authors computed the limit of the expected number of clusters for critical site percolation on the triangular lattice in a rectangle of fixed aspect ratio as the mesh goes to $0$. In terms of SLE$_{6}$ (the limit of the percolation interface), this boils down to computing the expectation of the number of commutings the latter is doing between the top and bottom boundaries of the rectangle or equivalently (by conformal invariance) 
 the expected number of times an SLE$_{6}$ bounces off $ \mathbb{R}_{+}$ and $ \mathbb{R}_{-}$ in a semi-ring region, see Fig.\,\ref{fig:HS}. % We preferred not to use this result (which would anyway require additional work to get large deviations estimates) and rather stayed on the continuous setup. 
 
 \begin{figure}[!ht]
  \begin{center}
  \includegraphics[width=16cm]{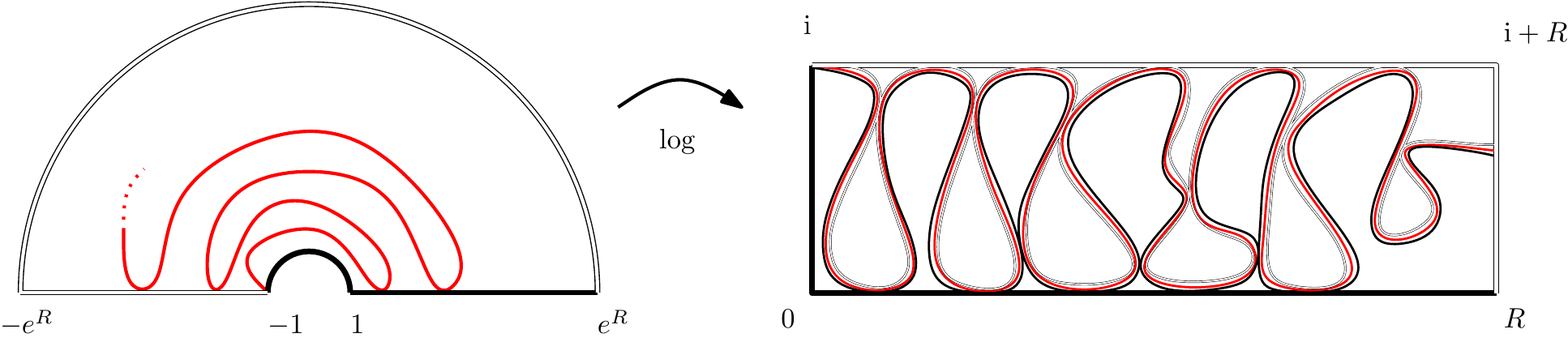}
  \caption{ \label{fig:HS}Connecting Proposition \ref{prop:commuteSLE} with the work \cite{HS11}.}
  \end{center}
  \end{figure}
 \end{remark}

%%%%%%%%%%%%%%%%SECTION
\section{Conformal measure on the boundary}
\label{sec:proof}
With all the ingredients that we have gathered we can now proceed to the proof of Theorem* \ref{thm:main}.
Consider the uniformization of a UIHPT onto $ \mathbb{H}$ such that the origin and target of the root edge are sent to $-\frac{1}{2}$ and $ \frac{1}{2}$ and $\infty$ to $\infty$. Recall also that the $k$th vertex on the right of the origin of $T_{\infty,\infty}$ has image $ \mathcal{X}_{k} \in  \mathbb{R}_{+}$. The sketch of the proof of Theorem* \ref{thm:main} $(iii)$ can be found in the introduction, however the following lines are a bit more technical since we will need precise estimates to rigorously derive the Hausdorff dimension of $\mu$.

%%%%%%%%%%%%%%%%%%%%SUBSECTION
\subsection{Discrete estimates}

\begin{proposition2} \label{prop:discretemass} Let  $ \eta, \varepsilon \in (0,1)$. Then there exist $c_{1}, c_{2} >0$ such that 
$$ \limsup_{n \to \infty} P \left( \left|\log \frac{\mathcal{X}_{ [ \varepsilon n]}}{\mathcal{X}_{n}}  - 3 {\log \varepsilon} \right|> \eta |\log \varepsilon|  \right) \leq c_{1} \exp(-c_{2} |\log \varepsilon|).$$
\end{proposition2}

\proof Fix  $  \varepsilon \in (0,1)$. Denote by $t_{ \varepsilon,n}$ the first time the exploration process triggers a peeling step that ``swallows'' or touches the $[ \varepsilon n]$th vertex on the right of the root edge, recalling \eqref{eq:hithor} we thus have
 $$t_{ \varepsilon,n}= \inf \big\{ k \geq 0 : \underline{\mathcal{H}}^+(k) \leq -[ \varepsilon n]\big\}.$$ (Note that this time is almost surely finite by Theorem* \ref{thm:exploSLE}.) Introduce then the number of discrete remaining commutings necessary to discover the $n$th point on the right boundary, that is
$$  {C}^{dis}( \varepsilon,n) = \inf \left\{ k \in \{1,3,5,\ldots\} :  \underline{\mathcal{H}}^+ \left( \tau^{(k)} ( t_{ \varepsilon,n})\right) \leq -n \right\}.$$
By  Theorem* \ref{thm:exploSLE} and using standard arguments as those developed in Section \ref{sec:discretebouncing}, the random variable $C^{dis}( \varepsilon,n)$ converges in distribution as $n \to \infty$ towards the random variable $\mathsf{ComStable}( \varepsilon,1)$ defined just before Corollary \ref{cor:commuteStable}. The same Corollary \ref{cor:commuteStable} thus entails that for every $ \eta >0$ there exist $c_{1},c_{2}>0$ such that
 \begin{eqnarray} \label{eq:limsup1} \limsup_{n \to \infty} P \left( \left|\frac{\log C^{dis}( \varepsilon,n)}{|\log \varepsilon|} - \frac{3 \sqrt{3}}{2 \pi} \right|> \eta \right) \leq c_{1} \exp(-c_{2} |\log \varepsilon|).  \end{eqnarray}
Let us now focus on the SLE exploration.   An easy adaptation of Proposition* \ref{prop:limitbouncing} implies that the number $C^{dis}( \varepsilon,n)$ is asymptotically equal as $n \to \infty$ to the number of commutings the SLE$_{6}$ is doing after having swallowed the point $  \mathcal{X}_{[ \varepsilon n]}$ until it swallows $ \mathcal{X}_{n}$, i.e. 
 \begin{eqnarray*} \big| C^{dis}( \varepsilon,n) - \mathsf{ComSLE}( \mathcal{X}_{[ \varepsilon n]}, \mathcal{X}_{n})  \big| &\xrightarrow[n\to\infty]{(P)} &0.  \end{eqnarray*}
Consequently, by the last display and \eqref{eq:limsup1} we have 
 \begin{eqnarray} \label{eq:lismup2} \limsup_{n \to \infty} P \left( \left|\frac{\log  \mathsf{ComSLE}(  \mathcal{X}_{[ \varepsilon n]}, \mathcal{X}_{n})}{|\log \varepsilon|} - \frac{3 \sqrt{3}}{2 \pi} \right|> \eta \right) \leq c_{1} \exp(-c_{2} |\log \varepsilon|).    \end{eqnarray}
Since the SLE$_{6}$ is independent of the map, if we condition on $  \mathscr{T}_{\infty,\infty}$, using  Corollary \ref{cor:commuteSLE} we get that for small enough $ \varepsilon>0$ we have (note that $ 4 \cdot \frac{ \sqrt{3}}{2 \pi} >1$)
 \begin{eqnarray} P \left(\left. \left|\frac{\log  \mathsf{ComSLE}(  \mathcal{X}_{[ \varepsilon n]}, \mathcal{X}_{n})}{|\log \varepsilon|} - \frac{3 \sqrt{3}}{2 \pi} \right|> \eta \ \right| \left\{\left|\log \frac{\mathcal{X}_{ [ \varepsilon n]}}{\mathcal{X}_{n}}  - 3 {\log \varepsilon} \right|> 4\eta  |\log \varepsilon|\right\}\right) &\geq& \frac{1}{2}.   \label{eq:lastdiaplsy}\end{eqnarray}
So that for small enough $ \varepsilon>0$:
 \begin{eqnarray*} &&\hspace{-3cm}\limsup_{n \to \infty} P \left( \left|\log \frac{\mathcal{X}_{ [ \varepsilon n]}}{\mathcal{X}_{n}}  - 3 {\log \varepsilon} \right|> 4\eta |\log \varepsilon|\right) \cdot \frac{1}{2}\\ &\underset{ \eqref{eq:lastdiaplsy}}{\leq}&  \limsup_{n \to \infty} P \left( \left|\frac{\log  \mathsf{ComSLE}(  \mathcal{X}_{[ \varepsilon n]}, \mathcal{X}_{n})}{|\log \varepsilon|} - \frac{3 \sqrt{3}}{2 \pi} \right|> \eta \right)\\ &\underset {\eqref{eq:lismup2}}{\leq}& c_{1} \exp(-c_{2} |\log \varepsilon|).  \end{eqnarray*} This completes the proof of the proposition*.
\qedhere \bigskip 

We will also rely on an adaptation of the last proposition* in order to compare the relative positions of $  \mathcal{X}_{k}$ and $ \mathcal{X}_{k'}$  when $k$ and $k'$ are of the same order.

\begin{proposition2} \label{prop:small} For every $u >0$ and for every $\eta >0$, there exist $c_{1},c_{2}>0$ such that we have 
$$ \limsup_{n \to \infty} P \left(	 \left|\log \frac{ \mathcal{X}_{[un]}}{ \mathcal{X}_{n}}\right| \geq \eta |\log \varepsilon| \right) \leq c_{1} e^{-c_{2}|\log \varepsilon|}.$$
\end{proposition2}
\proof[Sketch of the proof.] The proof uses the same arguments as in Proposition* \ref{prop:discretemass} so we only sketch it. Fix $ u <1$ for definiteness and consider $ C^{dis}(u,n)$ to be the number of commutings realized by the horodistance process between the discovery of the $[un]$th vertex on the right of $ \vec{e}$ and the $n$th one. On the one hand, as in Proposition* \ref{prop:discretemass}, $C^{dis}( u,n)$ converges as $n \to \infty$ towards $\mathsf{ComStable}( u,1)$ which is of order $1$. On the other hand, using the SLE$_{6}$ interpretation of the exploration (and Proposition* \ref{prop:limitbouncing}) we also get that $C^{dis}(u,n)- \mathsf{ComSLE}( \mathcal{X}_{[un]}, \mathcal{X}_{n}) $ converges towards $0$ in probability as $n \to \infty$. Using  Corollary \ref{cor:commuteStable} and \ref{cor:commuteSLE} we thus get that 
$$ 1 \approx \mathsf{ComStable}( u,1) \approx \mathsf{ComSLE}( \mathcal{X}_{[un]}, \mathcal{X}_{n}) \approx  \frac{ \sqrt{3}}{2\pi} \log \frac{ \mathcal{X}_{n}}{ \mathcal{X}_{[ u n]}},$$ so that $ \mathcal{X}_{[ u n]}$ and $ \mathcal{X}_{n}$ are of the same order of magnitude. Details are left to the reader.\endproof

We extend the definition of $ \mathcal{X}_{k}$ to every integer $ k \in \mathbb{Z}$ in a straightforward manner. 

\begin{proposition}[Re-rooting and symmetry] \label{prop:reroot2} For every $n \geq 0$ and every integer $u_{n}$ we have the following identity in distribution 
 \begin{eqnarray*} \left(\frac{ \mathcal{X}_{k+u_{n}}- \mathcal{X}_{u_{n}}}{ \mathcal{X}_{n+u_{n}}- \mathcal{ X}_{u_{n}}}\right)_{k  \in \mathbb{Z}} \overset{(d)}{=} \quad \left(\frac{ \mathcal{X}_{k}}{ \mathcal{X}_{n}}\right)_{k \in \mathbb{Z}} \overset{(d)}{=}\quad  \left(\frac{ \mathcal{X}_{-k}}{ \mathcal{X}_{-n}}\right)_{k \in \mathbb{Z}}.  \end{eqnarray*}
\end{proposition}
\proof For $n \geq 0$, the lattice $ \tilde{T}_{\infty,\infty}$ obtained from $T_{\infty,\infty}$ after re-rooting at the $u_{n}$th edge on the right of $ \vec{e}$ is still distributed as the UIHPT. The uniformization $ \tilde{ \mathscr{T}}_{\infty,\infty}$ of $ \tilde{T}_{\infty,\infty}$ (with the root edge sent to $[-1/2,1/2]$ and infinity to infinity) is obtained from $ \mathscr{T}_{\infty,\infty}$ by translation and dilation. Thus if $ (\tilde{ \mathcal{X}}_{k})$ denotes the positions of the vertices on the right of the root edge of the uniformization of $ \tilde{T}_{\infty,\infty}$ we get that 
$$ \left(\frac{ \mathcal{X}_{ k+u_{n}}- \mathcal{X}_{u_{n}}}{  \mathcal{X}_{n+u_{n}}-\mathcal{X}_{u_{n}}}\right)_{k \geq 0} = \left(\frac{\tilde{ \mathcal{X}}_{k}}{ \tilde{ \mathcal{X}}_{n}}\right)_{k \geq 0}.$$
Since $ ( \tilde{ \mathcal{X}}_{k})_{ k \geq 0}$ has the same law as $ ( \mathcal{X}_{k})_{ k \geq 0}$  the first identity in distribution follows. The second one is obtained by flipping $T_{\infty,\infty}$ horizontally, operation which leaves its distribution unchanged. \qedhere \medskip

Combining Proposition* \ref{prop:small} with Proposition \ref{prop:reroot2} we deduce that  
 \begin{eqnarray*} \limsup_{n \to \infty} P \Big(	 \star \geq \eta |\log \varepsilon| \Big) \leq c_{1} e^{-c_{2}|\log \varepsilon|} \end{eqnarray*} where $\star$ can be replaced by 
 $$ \star =   \frac{ \mathcal{X}_{[n/2]}}{ \mathcal{X}_{n}} \mbox{ or }  \frac{ \mathcal{X}_{n}- \mathcal{X}_{[n/2]}}{ \mathcal{X}_{n}} \mbox{ or } \frac{ \mathcal{X}_{[3n/2]}}{ \mathcal{X}_{n}} \mbox{ or } \frac{ \mathcal{X}_{[3n/2]} - \mathcal{X}_{n}}{ \mathcal{X}_{[3n/2]}- \mathcal{X}_{[n/2]}}.$$
After some manipulations this eventually implies
 \begin{eqnarray} \label{eq:corcor}\limsup_{n \to \infty} P \left(	 \left|\log \frac{ (\mathcal{X}_{[3/2n]}- \mathcal{X}_{n}) \wedge ( \mathcal{X}_{n}- \mathcal{X}_{[n/2]})}{ \mathcal{X}_{n}}\right| \geq \eta |\log \varepsilon| \right) \leq c_{1} e^{-c_{2}|\log \varepsilon|}.  \end{eqnarray}

\subsection{Dimension of the random measure}
Recall that we consider the random  measure $\mu_{n}$ defined by
$$ \mu_{n}= \frac{1}{n}\sum_{k=1}^n \delta_{ {\mathcal{X}_{k}}/{ \mathcal{X}_{n}}}.$$
 Hence $\mu_{n}$ is a random probability measure on $[0,1]$. We briefly remind the reader about the basics of convergence in distribution for random measures on $ \mathbb{R}$ (the interested should consult the authoritative reference \cite{Kal76} for proofs and more general statements and \cite{Gra77} for a smooth introduction). We endow the set $\mathcal{M}$ of all positive Radon measures on $\mathbb{R}$ with the topology  $\mathcal{T}$ of vague convergence, that is, the weakest topology 
which makes the mappings 
\[
\mu \in   \mathcal{ M} \mapsto \mu f:= \int_{\mathbb R} \mathrm{d}\mu f ,
\quad  f \in C_K,
\]
continuous. (Here, $C_K$ is the set of continuous functions 
$f: \mathbb R \to \mathbb R$ with compact support.)
A \emph{random measure} is a random element of the space
$( \mathcal{M}, \mathcal{T})$, viewed as a measurable space with
$\sigma$-algebra generated by the sets in $\mathcal{T}$. 
A sequence $\lambda_{1}, \lambda_{2}, \ldots$ of random measures 
converges in distribution towards a random measure $\lambda$ 
if for any bounded continuous mapping $F : ( \mathcal{M}, \mathcal{T}) 
\to \mathbb{R}$ we have $E[F(\lambda_{i})] \to E[F(\lambda)]$ as $i \to \infty$.

Actually, convergence of $\lambda_n$ to $\lambda$ in distribution
is equivalent to:
$\lambda_{n} f \xrightarrow[]{} \lambda f$,
for any continuous $f \in C_K$ (see Theorem 4.2 in \cite{Kal76}). The
latter convergence is convergence in distribution of real-valued 
random variables.
The set of all random probability measures on $[0,1]$ is tight for this convergence in distribution. Hence, from any sequence of integers going to $\infty$ we can extract a subsequence $n_{k} \to \infty$ such that there exists a random probability measure $\mu$ satisfying 
 \begin{eqnarray*} \mu_{n_{k}} & \xrightarrow[k\to\infty]{(d)} & \mu.  \end{eqnarray*}

%%%%%%%%%%%%%%%%%%%%SUBSECTION

\begin{center} \textit{To lighten notation, we suppose in the rest of this section that above extraction has been realized and that all the statements $ n \to \infty$ have to be interpreted along this subsequence.} \end{center}
% When $\mu$ has no atoms, the last convergence is equivalent to the fact that $ \mu_{n}[a,b] \to \mu[a,b]$ in distribution for every $0 \leq a<b \leq 1$. In the general case an application of Portmanteau theorem gives   \begin{eqnarray} \label{eq:convmesmes}
%\limsup_{n\to \infty} P( \mu_{n}[a, b] \geq x) \leq P ( \mu[a, b] \geq x) \qquad \mbox{and} \qquad   P ( \mu(a, b) > x) \leq \liminf_{n \to \infty} P( \mu_{n}(a, b) > x),  \end{eqnarray} for every $ x \geq 0$.  
To get the third part of  Theorem* \ref{thm:main} we will prove that balls of radius $r$ around typical points of $\mu$ roughly have volume $ r ^{1/3}$ when $r \to 0$. Indeed, Theorem* \ref{thm:main} $(iii)$ is a standard consequence of the following result* (see for example \cite[Lemma 4.1]{LPP95}):

  \begin{corollary2}[H\"older exponent] \label{cor:Holder} Almost surely, for $\mu$-almost all $x \in [0,1]$  we have 
  $$ \lim_{ r \downarrow 0} \frac{ \log \mu (B_{r}(x))}{ \log r} = \frac{1}{3},$$
  where $B_{r}(x) = [x- r, x+r]$ is the ball of radius $r$ around $x$.
 \end{corollary2}
 
 \proof Conditionally on $\mu_{n}$, let $X_{n}$ be a random point sampled accord to $\mu_{n}$ and similarly conditionally on $\mu$, let $X$ be sampled according to $\mu$. By definition of $\mu_{n}$ notice that we can write $X_{n} = \mathcal{X}_{ \lceil Un \rceil}/ \mathcal{X}_{n}$ where $U \in (0,1)$ is a uniform random variable independent of $T_{\infty,\infty}$ and $\lceil a \rceil$ is the lowest integer larger than $a$. Now fix $a < 1/3$ and write $x = r^{a}$. Since $\mu$ is the (subsequential limit) of the $\mu_{n}$'s we have 
  \begin{eqnarray} &&\hspace{-2cm}P\left( \mu B_{r/2}(X) \geq 8 x\right) \nonumber \\ & \leq & \limsup_{n \to \infty} P \left(\mu_{n} B_{r}(X_{n}) \geq 4x\right) \nonumber \\ & \leq & \limsup_{n \to \infty} P \left( \mu_{n}[X_{n},X_{n}+r] \geq 2x \right) + \limsup_{n \to \infty}P \left( \mu_{n}[X_{n}-r,X_{n}] \geq 2x \right). \label{eq:preuve2} \end{eqnarray}
By definition of $\mu_{n}$, note that the event $\{\mu_{n}[X_{n},X_{n}+r] \geq 2x\}$ can be written as 
$$\{\mu_{n}[X_{n},X_{n}+r] \geq 2x\} = \left\{ \frac{ \mathcal{X}_{ \lceil Un\rceil +\lceil 2x n \rceil -1}- \mathcal{X}_{\lceil Un\rceil}}{ \mathcal{X}_{n}} \leq r\right\}   \underset{ \mathrm{large \ } n's}{\subset} \left\{ \frac{ \mathcal{X}_{ [xn]+\lceil Un\rceil }- \mathcal{X}_{\lceil Un\rceil}}{ \mathcal{X}_{n}} \leq r\right\}. $$
We now write
 \begin{eqnarray*} && \hspace{-1cm}P \left( \frac{ \mathcal{X}_{ [xn]+\lceil Un\rceil }- \mathcal{X}_{\lceil Un\rceil}}{ \mathcal{X}_{n}} \leq r \right)\\ 
 & = & P \left( \frac{ \mathcal{X}_{ [xn]+\lceil Un\rceil }- \mathcal{X}_{\lceil Un\rceil}}{ \mathcal{X}_{n + \lceil Un\rceil}- \mathcal{X}_{\lceil Un\rceil}} \cdot \frac{ \mathcal{X}_{n+\lceil Un\rceil}- \mathcal{X}_{\lceil Un\rceil}}{ \mathcal{X}_{n}} \leq r \right)   \\ 
 & \leq & P \left( \log \frac{ \mathcal{X}_{ [xn]+\lceil Un\rceil }- \mathcal{X}_{\lceil Un\rceil}}{ \mathcal{X}_{n + \lceil Un\rceil}- \mathcal{X}_{\lceil Un\rceil}} \leq (1- \varepsilon) \log r \right) + P \left( \log \frac{ \mathcal{X}_{n+ \lceil Un \rceil}- \mathcal{X}_{\lceil Un\rceil}}{ \mathcal{X}_{n}} \leq  \varepsilon \log r \right) \\
  & \leq & P \left( \log \frac{ \mathcal{X}_{ [xn]+u_{n} }- \mathcal{X}_{u_{n}}}{ \mathcal{X}_{n +u_{n}}- \mathcal{X}_{u_{n}}} \leq (1- \varepsilon) \log r \right) + P \left( \log \frac{ (\mathcal{X}_{[3n/2]}- \mathcal{X}_{n}) \wedge ( \mathcal{X}_{n} - \mathcal{X}_{[n/2]})}{ \mathcal{X}_{n}} \leq  \varepsilon \log r \right) ,\end{eqnarray*} where we have chosen $ \varepsilon >0$ so that $(1- \varepsilon)/a >3$ and put $u_{n} = \lceil Un \rceil$. Indeed notice that when $U \in (0,1)$ then $[\lceil Un \rceil,n+ \lceil Un \rceil]$ contains either $[[n/2],n]$ or $[n,[3n/2]]$. We can take $\limsup$ and apply Proposition \ref{prop:reroot2} together with Proposition* \ref{prop:discretemass} to the first member of the right-hand side and \eqref{eq:corcor} to the second member of the right-hand side to  deduce that there exist constants $c_{1},c_{2}>0$ so that 
 $$\limsup_{n \to \infty} P \left( \mu_{n}[X_{n},X_{n}+r] \geq 2x\right) \leq c_{1} \exp(-c_{2} |\log r|).$$
 A similar reasoning holds for $ \limsup_{n \to \infty} P \left( \mu_{n}[X_{n}-r,X_{n}] \geq 2x \right)$. Gather-up the pieces of \eqref{eq:preuve2} and establishing the corresponding lower bound (left to the reader) we finally get that 
 $$ \limLD{\log \mu B_{r}(X)}{\log r}{ r \to 0} = \frac{1}{3}.$$
Taking $r = 2^{-k}$, the Borel--Cantelli lemma shows that 
$ \frac{\log \mu B_{2^{-k}}(X)}{ \log 2^{-k}} \to \frac{1}{3}$ almost surely as $k \to \infty$ which easily implies the statement of the corollary*. \endproof 
 
\begin{corollary2} The random probability measure $\mu$ is almost surely non-atomic. 
\end{corollary2}
\proof	 This is a straightforward consequence of the last corollary*. Indeed, if $\mu$ had a probability at least $ \varepsilon>0$ of having an atom of mass at least $ \varepsilon$ then $X$ would be located on a point of $\mu$-mass $ \varepsilon$ with probability at least $ \varepsilon^2$ and Corollary* \ref{cor:Holder} would not hold. \endproof

%%%%%%%%%%%%%%%%%%%%SUBSECTION
\subsection{Full support}

\begin{proposition2} The random probability measure $\mu$ has topological support equal to $[0,1]$ a.s. \end{proposition2}

\proof Let us argue by contradiction and suppose that with positive probability $\mu$ has not full support i.e. $P(\exists x \in [0,1] \mbox{ and } \varepsilon >0 : \mu B_{ \varepsilon}(x) =0 ) >0$. By compactness, we can thus suppose that for some some $ x \in [0,1]$ and $ \varepsilon \in (0,1)$ we have $ P( \mu B_{ 2\varepsilon}(x) = 0 ) \geq 2\varepsilon$. We will further assume that $x > 2 \varepsilon$ and 
$$ P( \mu B_{ 2\varepsilon}(x) = 0 \mbox{ and } \mu B_{ \varepsilon/2}(0) \geq 2 \varepsilon ) \geq 2\varepsilon.$$ (The boundary case $x=0$ is similar and left to the reader). Going back to a discrete level we deduce that we have a sequence of integers $n_{k}$ such that for every $ \delta >0$ the event $$ E_{n_{k}} = \left\{\mu_{n_{k}} B_{ \varepsilon}(x) \leq \delta \mbox{ and } \mu_{n_{k}} B_{ \varepsilon}(0) \geq  \varepsilon \right\}$$ is asymptotically of probability larger than or equal to $\varepsilon$. As before, we will lighten notation and assume that all the statements involving $n$ in the following lines have to be restricted to this subsequence. 
Unsurprisingly, we consider the exploration of the UIHPT by an  SLE$_{6}$. Since the SLE$_{6}$ process is independent of the map (and thus of its uniformization) conditionally on $ \mathscr{T}_{\infty,\infty}$ (and a fortiori on $E_{n}$) there is a positive probability $c>0$ (depending on $x$ and $ \varepsilon$, but not on $\delta$) that the curve $\gamma$ touches the interval $[x- \varepsilon, x+ \varepsilon]$ then touches $ \mathbb{R}_{-}$ and finally touches the interval $[ x- \varepsilon, x+ \varepsilon]$ again, that is with our notation 
$$ P \big( \mathsf{ComSLE}(x- \varepsilon, x+ \varepsilon) \geq 5\big) >c.$$  Recall now the notation $ t_{ \varepsilon,n}$ and $ C^{dis}( \varepsilon, n)$ from the proof of Proposition* \ref{prop:discretemass}. In terms of horodistance process, using (a variant of)  Proposition* \ref{prop:limitbouncing}, with high probability these visits in $[ x - \varepsilon, x + \varepsilon]$ by the SLE process can be associated with some peeling time $k$ such that $ \underline{ \mathcal{H}}^+(k) = \mathcal{H}^+(k+1/2)$ where $ \underline{ \mathcal{H}}^+(k) \in ( \varepsilon n , n)$. That is with some time  $\tau^{(k)} ( t_{ \varepsilon,n})$ and $ \tau^{(k+2)}( t_{ \varepsilon,n})$ for $k \leq C^{dis}( \varepsilon, n)$. By definition of the event $E_{n}$, for large $n \geq 0$ we have 
$$ P \left( \exists k \leq C^{dis}( \varepsilon,n) :  \Big|\underline{ \mathcal{H}}^+\big( \tau^{(k)} ( t_{ \varepsilon,n})\big) - \underline{ \mathcal{H}}^+\big( \tau^{(k+2)} ( t_{ \varepsilon,n})\big) \Big| \leq \delta n \right) > c',$$  for some positive constant $c'>0$ (depending on $x$ and $ \varepsilon$ but not on $\delta$). Taking the scaling limit of the horodistance processes using Theorem* \ref{thm:exploSLE}, a similar statement must  hold for the stable processes $(S^+,S^-)$ more precisely: There exists some constant $ c''>0$ such that for any $ \delta>0$ we have 
$$ P \left( \exists k \leq \mathsf{ComStable}( \varepsilon,1) :  \Big|\underline{S}^+\big( \xi^{(k)} ( \vartheta_{ \varepsilon})\big) - \underline{S}^+\big( \xi^{(k+2)} ( \vartheta_{ \varepsilon})\big) \Big| \leq \delta \right) > c''.$$
Letting $\delta \to 0$ we reach a contradiction since $ k \mapsto  \underline{S}^+(\xi^{(k)}(\vartheta_{ \varepsilon}))$ is strictly decreasing a.s.\qedhere

%%%%%%%%%%%%%%%%SECTION
\section{Discussion and comments}
 \label{sec:comments}
First of all, let us mention that although this paper was focused with the case of triangulations, we do not perceive any major conceptual obstacle in deriving the same results (provided that a variant of $(*)$ holds) for other classes of maps like quadrangulations or general planar maps. Indeed, though the peeling transitions of Section \ref{sec:peeling} are more complicated, they exhibit the same large-scale property and a variant of Proposition \ref{prop:deuxstables} should hold (see \cite{ACpercopeel}). 

A first natural question is to sharpen Proposition* \ref{prop:discretemass} to get a (more) precise result on the location $ \mathcal{X}_{n}$ of the $n$th vertex on the right of the root in the uniformized UIHPT:

\begin{question}  Prove that $ \displaystyle \frac{\log  \mathcal{X}_{n}}{\log n} \xrightarrow[n\to\infty]{(P)} 3$. Do we actually have $ (n^{-3}\mathcal{X}_{n})_{n\geq0}$ tight? 
 
\end{question}

%%%%%%%%%%%%%%%%%%%%SUBSECTION
\subsection{Discussion on the $(*)$ property}
We give here some elements supporting  $(*)$.
\begin{itemize}
\item First of all, we have seen in Proposition \ref{prop:boundetai} that the tail of $\eta_{i}^+$ is very light and that $\eta_{i}^+ \leq \log^{5+ \varepsilon} i$ eventually. We believe that the exponent $5$ could be brought down to $1$ with some work. Thus only a collective behavior of the $\eta_{i}^+$ could violate $(*)$
\item Although not independent, the $\eta^+_{i}$ decorrelate. Quantifying the speed of mixing is a path towards a proof of $(*)$. In particular, it happens that during the exploration the SLE$_{6}$ creates a ``bubble'': it explores for some time a new region connected to the past by a single triangle. This should  correspond on a continuous level to the pinch-points of the SLE$_{6}$. On an intuitive level the $\eta_{i}^+$ in such a region can be thought as independent of the past. 
\item Another heuristic for the decorrelation of the $\eta_{i}^+$ is the following. On a rough level, one can imagine that $\eta_{i}^+$ is correlated with  $\eta_{j}^+$ with $i <j$ only if during the $j$th peeling step the curve $\gamma$ comes back close to its location at time $i$ and touches a part of the continuous hull which was close (say within horodistance $10$) from the edge $a_{i}$:
$$ \eta_{i}^+ \mbox{ correlated with } \eta_{j}^+ \quad \iff \quad {\left\{\begin{array}{l} \mbox{The edge $a_{i}$ is still on the boundary of $K_{j}$}\\
\mbox{and $ \mathcal{H}^+(j) \approx \mathcal{H}^+(i)$} \end{array}\right\} },$$ we call this event $ \mathrm{Corr}(i,j)$. 
However, it is easy to see that conditionally on $ \mathrm{Corr}(i,j)$ and on the past before $j$, there is probability bounded away from $0$ that the next peeling step  swallows $a_{i}$. Hence, the number of $j \geq i$ such that $ \mathrm{Corr}(i,j)$  holds has an exponential tail. So if we believe in the last display we would have $E[( \sum_{i \leq n} \eta^+_{i})^2] \leq C n$ for some $C>0$ and by Markov inequality $ \sum_{i \leq n} \eta_{i}^+$ could not be much larger than $\sqrt{n}$ (condition $(*)$ just needs $o(n^{2/3})$).
\end{itemize}

\paragraph{Percolation interface.} Recall from Remark \ref{rem:percoexplo} that Theorem* \ref{thm:exploSLE} and Theorem \ref{thm:exploperco} show that the interfaces in critical site percolation on the UIHPT (see Theorem \ref{thm:exploperco}) converge, in terms of horodistances, towards the SLE$_{6}$ exploration. Obviously, one can wonder if a geometric statement holds: Conditionally on $T_{\infty,\infty}$ sample a critical site percolation (with parameter $1/2$) with boundary condition  black-white as in Section \ref{sec:perco}. This naturally defines a curve on the Riemann surface $ \mathbf{T}_{\infty,\infty}$ (join with straight lines the middle of the edges $\bullet-\circ$ in each peeled triangle), see Fig.\ref{fig:exploration}. Denoting by $ \mathcal{P}$ the image of this interface on the uniformization $ \mathscr{T}_{\infty,\infty}$ one would like to show that $\lambda \cdot \mathcal{P}$ converges as $\lambda \to 0$ in distribution (for the Hausdorff distance on any compact sets of $ \mathbb{H}$) towards a standard chordal SLE$_{6}$.

%%%%%%%%%%%%%%%%%%%%SUBSECTION
\subsection{Towards a characterization of $\mu$}
In this work we used Theorem* \ref{thm:exploSLE} on exploration along an SLE$_{6}$ in order to derive a few properties of the sequential limits of the $\mu_{n}$'s (Theorem* \ref{thm:main}). Theorem* \ref{thm:exploSLE} is actually much stronger and \emph{it may be the case that it implies alone the convergence of the $\mu_{n}$'s towards the random measure $\tilde{\mu}_{\rho}$ with $ \rho = \sqrt{8/3}$ of Question \ref{ques:limit}}. Let us comment on this.\medskip

First, we can extend the definition of $\mu_{n}$ to a infinite random measure on $ \mathbb{R}$: for $n \geq 1$ in this section we let 
$$ \boldsymbol{\mu}_{n} = \frac{1}{n} \sum_{k \in \mathbb{Z}} \delta_{ \mathcal{X}_{k}/ \mathcal{X}_{n}},$$ where $ \mathcal{X}_{k}$ are the location of the vertices of the boundary of the UIHPT in the uniformization $ \mathscr{T}_{\infty,\infty}$. Clearly, $ \boldsymbol{\mu}_{n}$ is an infinite measure and $ \boldsymbol{\mu}_{n}( \cdot \cap [0,1]) = \mu_{n}$. Using similar arguments as those developed in Section 	\ref{sec:proof} one can show that $ (\boldsymbol{ \mu}_{n})_{n}$ is tight and we denote $ \boldsymbol{\mu}_{\infty}$ one possible limiting random measure on $ \mathbb{R}$. Let us now present a corollary of Theorem* \ref{thm:exploSLE} just in terms of $ \boldsymbol{\mu}_{\infty}$:

If $\nu$ is a random measure on $ \mathbb{R}$, independently of it let $(\gamma_{t})_{t \geq 0}$ be a chordal SLE$_{6}$ on $ \mathbb{H}$ starting from $0$ and denote by $ \Gamma_{ \nu}^{+} = \{(t,\nu[0, \gamma_{t}] ) : \gamma_{t} \in \mathbb{R}_{+}\}$ and similarly for $ \Gamma^-_{ \nu}$. On the other hand, let $S^+$ and $S^-$ be two independent stable processes and recall the notation $ \underline{S}^+$ and $ \underline{S}^-$ for their running infimum processes. We denote $ \mathcal{M}^+ = \{ (t, -S^+_{t}) : S^+_{t} = \underline{S}^+_{t}\}$ and define $ \mathcal{M}^-$ similarly from $S^-$. 

   Then the analog of Theorem* \ref{thm:exploSLE} gives the following identity in distribution up to time reparametrization of the first coordinate
  \begin{eqnarray}  \label{eq:caracterization} \big( \Gamma^{+}_{ \boldsymbol{\mu}_{\infty}}, \Gamma^{-}_{ \boldsymbol{\mu}_{\infty}}\big) \overset{(d)}{=} \big( M^+,M^-\big).  \end{eqnarray}
  
\begin{question}[Characterization] Does the last display characterize the law of $ \boldsymbol{\mu}_{\infty}$? \label{ques:characterization}
\end{question}  

A positive answer the last question would imply a positive answer to Question \ref{ques:limit}. Indeed, in the recent work \cite{DMS14} the authors show (among other things) that we have the identity in distribution up to time reparametrization of the first coordinate
 \begin{eqnarray}  \label{eq:caracterization?} \big( \Gamma^{+}_{\mu_{\varrho}}, \Gamma^{-}_{ \mu_{\varrho}}\big) \overset{(d)}{=} \big( M^+,M^-\big), \end{eqnarray}
where $\mu_{\varrho}$ is the random measure defined in Question 1 (see Introduction).

\subsection{Full-plane}
To conclude, we would like to mention that our tools might be extended to investigate the conformal structure of random maps ``in the bulk'' . Similarly as in this work, the limiting random measures in the plane induced by the conformal uniformization of random planar maps is conjectured to be described by the exponential of the GFF with the same parameter $ \varrho= \sqrt{8/3}$ and should be of fractal dimension $2/3$.

One should this time use a full-plane version of the SLE$_{6}$ instead of the chordal version. Again, the locality property of the SLE$_{6}$ will imply that the exploration is Markovian (see \cite{Ang03,BCsubdiffusive} for Markovian exploration of the UIPT/UIPQ) and thus under a similar $(*)$ condition one would be able to understand the change in the boundary length on the left and right of the peeling point. Two difficulties will then arise: first the scaling limit of the variations of the boundary length in a full-plane exploration (of the UIPT say) now involves $3/2$-stable processes \emph{conditioned to survive}, second the geometric property  (bouncing on the real line) used in this work has to be replaced by another geometric property computable in terms of the ``horodistances'': the \emph{winding number}. We cherish the hope of pursuing these ideas in future works.

  \end{document}